\documentclass[11pt]{article}
\usepackage[english]{babel}
\usepackage[latin1]{inputenc}
\usepackage{exscale}
\usepackage[centertags]{amsmath}
\usepackage{amsfonts,amstext,amssymb,amsthm}
\usepackage{mathrsfs}
\usepackage{comment}
\usepackage{newlfont}
\usepackage{enumerate}
\usepackage{graphicx, graphics}
\usepackage{these}
\usepackage{color}
\usepackage[usenames,dvipsnames]{xcolor}
\definecolor{BlueFonse}{rgb}{0,0,1}
\definecolor{BlueFonse1}{cmyk}{1,0,0,0.7}
\usepackage[colorlinks=true,linkcolor=blue,citecolor=ForestGreen]{hyperref}
\usepackage{hyperref}
\usepackage{url}
\usepackage{tikz}
\usetikzlibrary{calc}
\usetikzlibrary{matrix,arrows,decorations.pathmorphing,decorations.pathreplacing,shapes.misc}
\usetikzlibrary{decorations.markings}
\usetikzlibrary{patterns}
\usetikzlibrary{snakes}
\title{Regularity of isoperimetric regions that are close to a smooth manifold\footnote{Work supported by the grant "Borse per l'Estero" of INDAM}{}} 
\author{Stefano Nardulli}
\begin{document}
      \maketitle
      \tableofcontents
\section{Introduction}
\subsection{A regularity theorem}

In this article we give some regularity results for the boundary of isoperimetric regions in a smooth complete Riemannian manifold with variable metric using the theory of regularity of Allard. In the remaining part of this paper we always assume that all the Riemannian manifolds $(M, g)$ considered are smooth with smooth Riemannian metric $g$. $\mathbb{M}^n_k$ will denote the simply connected space form of constant sectional curvature $k$. For every $m\in\mathbb{N}$, we denote by $\mathcal{H}_g^m$, the $m$-dimensional Hausdorff  measure associated to the metric space $(M, d_g)$, where $d_g$ is the canonical length space metric associated to $g$, i.e., $d_g(x,y):=\inf\{l_g(\gamma): \gamma\:\text{is a piecewise}\:C^1\text{curve joining}\:x\:\text{to}\:y\}$, $l_g(\gamma)$ represents the length of the curve $\gamma$ with respect to the Riemannian metric $g$. $Vol_g=\mathcal{H}_g^n$ will denote the canonical Riemannian measure induced on $M$ by $g$, by $A_g$ we will denote often $\mathcal{H}_g^{n-1}$, $\M_g$ indicates the mass of a current, the notation here is the standard one of Federer's book \cite{Fed}. When it is already clear from the context, explicit mention of the metric $g$ will be suppressed. 
All along this text we will encounter a lot of constants that depends on various geometrical quantities and the metric $g$, when a constant $c$ depends on the metric g and on its first and/or second derivatives continuously, we will denote this fact by writing $c=c(g, \partial g, \partial^2 g)$. In what follow we will also be concerned with isoperimetric regions that are close to a fixed open relatively compact set $B\subseteq M$ with smooth boundary $\partial B$. As a consequence of this fact we need to consider also the dependence on $B\subseteq M$ and a fixed once at all $\xi\in\mathfrak{X}(M)$, where $\mathfrak{X}(M)$ is the set of smooth vector fields defined on $M^n$, these objects depending only on the differentiable structure of $M^n$ (that is considered fixed in this work) without any reference to a Riemannian metric. To deal with the case of variable metrics we will use the celebrated Nash's isometric embedding Theorem. So in general our constants $c$ will be of two kinds. The constants of the first kind are of the form $c=c(B, \xi, g, \partial g, \partial^2 g)$ when they do not depend on the Nash's isometric imbedding in some higher dimensional Euclidean space $i_g:(M^n, g)\to (\R^N, \delta)$ for some $N>n$, where $\delta$ is the Euclidean canonical metric. The constants are of the second kind when they depend on a Nash's isometric imbedding. In this latter case we write $c=c(B, \xi, i_g, g, \partial g, \partial^2 g)$. However, as we will see later in Section \ref{Section:VariableMetrics}, $C^2$ dependence on $i_g$ means to depend $C^4$ on the metric so in general another kind of typical constants that we will encounter is $c=c(B, \xi, i_g, g, \partial g, \partial^2 g)=c(B, \xi, g, \partial g, \partial^2 g, \partial^3 g,\partial^4 g)$. 
\begin{Def}
     Let $(M, g)$ be a smooth (possibly non-complete) Riemannian manifold of dimension $n$. 
     We denote by $\tau_{M}$ the class of relatively compact open sets of $M$ with          
     $C^{\infty}$ boundary.
    The function $I_{(M,g)}:[0,Vol_g(M)[\rightarrow [0,+\infty [$ defined by
    \begin{equation}
               I_{(M,g)}(v):=\inf\left\{A_g(\partial \Omega )\:|\:\Omega\in \tau_{M}, Vol_g(\Omega)=v\right\},
    \end{equation}
is called the \textbf{isoperimetric profile function} (shortly the \textbf{isoperimetric profile}) of the manifold $M$.
We define an \textbf{isoperimetric region for volume $v$} as an $n$-dimensional integral normal current $T$, such that $\M_g(T)=v$ and $\M_g(\partial T)=I_{(M,g)}(v)$.
\end{Def}
The regularity theory for minimizing currents, inaugurated by Ennio De Giorgi in codimension $1$ (see for example \cite{DeGiorgi}), Federer and Fleming in any codimension, and fully developed in the work of Almgren and Allard, shows that isoperimetric regions are almost smooth.  Precisely, they are submanifolds with smooth boundary on the complement of a singular set of codimension at least equal to $7$ \cite{Alm}. In codimension $1$ one can compare with \cite{GMTamanini} in which the theory of finite perimeter sets is adopted.  On the other hand, for manifolds $M^n$ of dimension $n\ge8$ there can be minimizing currents whith non-smooth boundary (see \cite{Alm}, \cite{Morg1}, \cite{BDeGG}). The first result along these lines, due to Bombieri, De Giorgi, and Giusti \cite{BDeGG}, shows that the cone $C:=\{ (x,y)\in \mathbb{R}^4\times \mathbb{R}^4 : |x|=|y| \}$, as conjectured by James Simons is singular at the origin and has minimal area in $\mathbb{R}^{8}$. In every ball of $\mathbb{R}^{n}$, such a current is a minimal hypersurface. Coming back to the isoperimetric problem, consider a point $p$ belonging to the support of the boundary of the current $T$, $p\in spt(\partial T)$, for some isoperimetric region $T$, then the tangent cone of $\partial T$ at $p$ have to be area minimizing in $T_pM^n$. If the point $p$ is regular then the tangent cone of $\partial T$ at $p$ is an hyperplane. If $p$ is a singular point the tangent cone at $p$ could be a genuine cone. In fact, there are examples with non void singular part, for more details about this matter we recommend the lecture of Proposition $3.5$ of \cite{Morg1} in which is proved that if $T$ is an isoperimetric region $p\in\partial T$ and the tangential tangent cone of $\partial T$ at $p$ is a hyperplane then $p$ is a regular point. Almgren's Theorem is thus optimal. Therefore, additional conditions are required to get more regularity in higher dimension.

\bigskip

The aim of this article is to show in Theorem \ref{T4} that an isoperimetric region, sufficiently close in the flat norm to a domain $B$ with smooth boundary $\partial B$, is also smooth and very close to $B$ in the $C^{2
,\alpha}$ topology 
and every $\alpha\in ]0,1[$. For further applications of this theorem, we also allow that the Riemannian metric $g$ of $M^n$ to be variable. We refer the reader to the last Section \ref{Section:VariableMetrics} for the precise meaning and definitions required to state our main Theorem \ref{T4}, especially Definitions \ref{Def:Petersen} and \ref{Def:WhitneyTopology} of topologies.
\begin{Res}
\label{T4}
Let $(M^n,p, g_{\infty})$ be a pointed Riemannian manifold of class $C^{\infty}$ with bounded geometry, $g_j$ a sequence of Riemannian metrics of class $C^{\infty}$ converging to $g_{\infty}$ in the fine $C^4$-topology or such that  $(M^n,p, g_j)$ converges to $g_{\infty}$ in the usual pointed $C^4$-topology. Let $B$ be an open relatively compact domain of  $M$ with smooth boundary $\partial B$. Consider $T_j$ a sequence of isoperimetric regions of $(M^n, g_j )$ such that 
\begin{equation}\label{Eq:Thm1Statement}
\mathbf{M}_{g_{\infty}} (B-T_j)\rightarrow 0,
\end{equation} 
where $\mathbf{M}_{g_{\infty}}$ denotes the mass of a current in the metric $g_{\infty}$. Then $\partial T_j$ is the graph  in normal exponential coordinates of a function  $u_j$ on $\partial B$. Furthermore, for all $\alpha\in ]0,1[$, $u_j \in C^{2,\alpha}(\partial B)$ and $||u_j||_{g_{\infty},C^{2,\alpha}(\partial B)}\rightarrow 0$. If the convergence of the metric is in the fine $C^{m,\alpha}$-topology then we also have $||u_j||_{g_{\infty},C^{m+1,\alpha}(\partial B)}\rightarrow 0$.
\end{Res}
\begin{Rem}
In the assumptions of Theorem \ref{T4}, rather than \eqref{Eq:Thm1Statement} we can use the following equivalent condition $V_{g_j}(B-T_j)\to 0$.
\end{Rem}
\begin{Rem} Notice that $||u_j||_{g_{\infty},C^{2,\alpha}(\partial B)}\rightarrow 0$ is equivalent to $||u_j||_{g_j,C^{2,\alpha}(\partial B)}\rightarrow 0$.
\end{Rem}
\begin{Rem} Observe that, since $\partial B$ is compact the spaces $C^{2,\alpha}(\partial B)$ are independent of the metric $g$, although the norm $||u_j||_{g,C^{2,\alpha}(\partial B)}$ depends on $g$. 
\end{Rem}
\subsection{Previous results}
We can find a particular case of Theorem \ref{T4} in the article \cite{MJ} of David L. Johnson and F. Morgan, from which we took a lot of inspiration. Indeed, these authors show that  in a compact manifold isoperimetric regions for small volumes are close to small geodesic balls. One can also consult Theorem $5$ of \cite{Ros} for analogous results in the case of a compact Riemannian manifold, but with a different proof. Here we follow the same ideas of \cite{MJ} adapted to our more general situation, paying attention to give intrinsic metric proofs from which very accurate estimates of the $C^{2,\alpha}$ norms and their dependence on the geometric bounds of $M$ and on the geometric bounds of the isometric imbeddings $\partial B\hookrightarrow M\hookrightarrow \R^N$ are given. 

\subsection{Some applications of Theorem \ref{T4}}

From Theorem \ref{T4}, we can argue that if, for a $\bar{v} >0$, all the isoperimetric regions in volume $\bar{v}$ are smooth, then the isoperimetric regions  for volumes $v$ close to $\bar{v}$ are smooth too. Under this condition, we could be able to reduce the isoperimetric problem for volumes close to $\bar{v}$ to a variational problem in finite dimension, as developed in \cite{NarAnn} and \cite{NarCalcVar} under small volumes assumptions. An analogous program is carried out in a separate paper see \cite{GNP} Lemma $6$. In 
\cite{NarAnn} and \cite{NarCalcVar} we used Theorem \ref{T4}, for showing that for small volumes the isoperimetric regions are pseudo-bubbles. But the range of application of this theorem and its straightforward generalizations is much wider.
\subsection{Sketch of the proof of Theorem \ref{T4}}

 First, assume that the metric is fixed, i.e., $g_j =g$, for every $j$. We make essential use of Allard's regularity theorem, see \cite{All}, Theorem $8.1$, which states that, if a varifold $V=\partial T\ni a$, has in a ball $B(a,r)$, a weight $||V||(B(a,r))$ sufficiently close to $\omega_{n-1}R^{n-1}$ (where $\omega_{n-1}$ is the volume of the unit ball of $\mathbb{R}^{n-1}$), and controlled first variation (i.e., mean curvature) in suitable $L^p$-norm, then $V$ is, locally, the graph of a function $u\in C^{1,\alpha}$. A regularity theorem for elliptic partial differential equations and a bootstrap argument imply that $u\in C^{\infty}$, and also give upper bounds for $||u||_{C^{2,\alpha}}$, via Schauder estimates.

	In order to show that $\partial T_j$ satisfies the conditions of Allard's regularity theorem, we  compare $\partial T_j$ to suitably chosen deformations with fixed enclosed volume. This is the matter in which is involved subsection \ref{Compvol}. 

Unfortunately for our purposes, Allard's theorem is stated in Euclidean spaces, hence we have to give a Riemannian version via isometric embedding of Riemannian manifolds in Euclidean spaces. Furthermore we need to control the second fundamental form of the isometric embeddings relative to different metrics on $M$. To make this possible we use a fine analysis of the proof of the Nash's isometric embedding theorem that M. Gromov did in \cite{Gr2}, this highlights the fact that free isometric embeddings can be chosen to depend continuously on the metric. It is worth noting here that we follow the scheme of the proof of Theorem 2.2 \cite{MJ}, with the difference that our context is more general because we consider arbitrary volumes instead that only small volumes, noncompact manifolds instead of compact ones and the proofs are made intrinsic as we can, when in the paper \cite{MJ} all is done using an isometric embedding of a compact manifold. Moreover we make an extra effort to find effective bounds for the constants involved and some improved arguments and some details that in \cite{MJ} are not mentioned, especially the way in which constants are calculated. In section \ref{calculs} explicit calculations of the mean curvature operator of a normal graph over a smooth $(n-1)$-dimensional submanifold are done. 
Lemma \ref{bornecm} on how to bound uniformly the curvatures is like in Theorem 2.2 of \cite{MJ} with the suitable modifications required to fit the case of a noncompact manifold. The compensation Lemma \ref{lr2} from Section \ref{Compvol} is new in the literature for the intrinsic metric geometry context and because it is done in a small ball centered on the boundary of the fixed domain $B$, rather than on the boundary of $T$. This features allows us to avoid the classical dependence of the constant measuring the distortion of area, on $T$. We want  that the constants involved depends just on $B$, uniformly on $T$, for values of $Vol_g(B\Delta T)$  smaller than a constant that depends just on $B$, and the bounds of the geometry of $M$ that in turn depend only on $g,\partial g, \partial^2g$. As a consequence of Lemma \ref{lr2} we have Lemma \ref{lr5} that permits to check one of the hypothesis of the Allard's regularity  Theorem, the intrinsic feature of these arguments is also new. The confinement Lemma \ref{tflat1} via the monotonicity formula I did not find in the literature but it is classical in the Euclidean case and perhaps it already exists for Riemannian manifolds at least for minimal submanifolds using Nash's embedding Theorem. An alternative proof of Lemma \ref{tflat1} is that of the confinement Lemma \ref{Lemme:UniformlyBoundedDiameter}, which is inspired by arguments on boundedness of isoperimetric regions in Euclidean space that can be found in \cite{Morgan94} combined with the technical compensation Lemma \ref{lr2}, and it is new in this form. The main difficulty encountered to adapt the Euclidean argument of \cite{Morgan94} to the Riemannian bounded geometry case consists in using an Euclidean type isoperimetric inequality for small volumes, but just with this we only produced Theorem $3$ of \cite{NarAsian}. To obtain the full generality of Lemma \ref{Lemme:UniformlyBoundedDiameter} we need the technical Lemma \ref{lr2} that permits us to control the variation of the perimeter of a deformation of volume $\Delta v$ of $T$, $|\Delta A|$ by the quantity $c(n, k, v_0, B, g, \partial g,\partial g^2)\Delta v$. How to apply the Allard's regularity theorem keeping track of the constants and the way in which they depend on the geometry of the noncompact manifold and on $B$, I did not find in the literature. The Schauder's estimates are classical and Hopf comparison theorems are classical. I put all the details required in our context for completeness and again to keep track of the constants and the quantities on which they depend in view of the subsequent application to the case of variable metrics. All these features make the arguments available for the case of variable metrics, arbitrary volumes, intrinsically without using an embedding of $M^n$ into some higher dimensional Euclidean $\R^N$. Unfortunately to achieve this program I need a version of the Allard regularity results stated in the Riemannian case, keeping track of how the constants depend on geometric quantities. This is a task that I did not write up because the details are cumbersome. On the other hand, using the remark of Gromov about the Nash embedding Theorem and keeping track of the constants involved I can overcome this difficulty paying the price of loosing the optimal $C^2$ convergence of the metrics in favor of the stronger $C^4$ convergence. The reader is invited to compare the last page for more details about this point. 	
\subsection{Plan of the article}
\begin{enumerate}
           \item In subsection \ref{calculs} a useful purely differential geometric formula for the mean curvature operator of a normal graph is given. Section \ref{chap3} provides Riemannian versions of three classical results of geometric measure theory: Allard's regularity theorem, the link between first variation and mean curvature in the case of currents and varifolds, the monotonicity formula. 
           \item Section \ref{Sec:NormalGrphaTheorem} is the core of the paper and gives the proof of Theorem \ref{T4} in case of a fixed metric, namely Theorem \ref{tr2}. It starts by a detailed  sketch of the proof. This part has the aim of elucidating the basic ideas involved in the proof of  Theorem \ref{tr2}.
           \item Section \ref{Section:VariableMetrics} deals with the general case of variable metrics and the final part of the proof of Theorem \ref{T4}.
\end{enumerate}

\subsection{Acknowledgements}  This paper is an extended and improved version of part of author's Ph.D. thesis, written under the supervision of Pierre Pansu at the University of Paris-Sud (XI-Orsay). The author gratefully acknowledges Renata Grimaldi of University of Palermo for the fruitful discussions that he had during his Ph.D. studies. The author wishes to express his gratitude to Pierre Pansu, for the many helpful suggestions during the preparation of the paper. It is worth to thank Frank Morgan for useful comments, making possible to improve the presentation of the manuscript. Finally I wish to thank Michael Deutsch for improving the english of the final text and my student Luis Eduardo Osorio Acevedo who helped me with the figures. The author wishes to thank also the Istituto Nazionale di Alta Matematica "Francesco Severi" of Rome (Italy) for financial support via the grant "Borse per l'estero".
\section{Regularity Theory}\label{chap3}
The aim of this section is to adapt several classical results of geometric measure theory stated in Euclidean spaces to arbitrary Riemannian manifolds.

\subsection{Notations}

In this section we are concerned with a Riemannian manifold $(M^n,g)$ of class at least $C^2$, with bounded second fundamental form, and we keep fixed an isometric embedding $i:(M^n,g)\hookrightarrow(\mathbb{R}^N,\delta)$, where $\delta$ is the Euclidean metric. We denote 
$$\beta_i =||II_{i\hookrightarrow M}||_{\infty ,g}<+\infty,$$ where  
$II_{i\hookrightarrow M}$ is the second fundamental form of the embedding $i$ and $||.||_{\infty ,g}$ is the supremum norm taken on $(M^n,g)$. In fact this is not a big restriction for the proof of our Theorem \ref{T4} because by Lemma \ref{Lemme:UniformlyBoundedDiameter} all happen in a bounded neighborhood of $B$ and the proof of Lemma \ref{Lemme:UniformlyBoundedDiameter} is completely intrinsic and independent from any isometric immersion into Euclidean spaces. Since in this neighborhood the second fundamental form is always bounded we loose no generality in making this assumption.
We observe, incidentally, that the second fundamental form depends on first and second derivatives of the embedding $i$ by continuous functions. Hence, if we have $2$ embeddings $i_1$, $i_2$ that are  $\varepsilon$ close in the $C^2$ topology, then $\beta_{i_1}$, $\beta_{i_2}$ will be $const.\varepsilon$ close and the constant is independent of embeddings $i_1$, $i_2$. Indeed the constant depends only on $M$ and the intrinsic metric but unfortunately this dependence is $C^4$.
\subsection{Mean curvature vector based on an hypersurface}\label{calculs}
For further applications, we now give a formula for mean curvature of an hypersurface $N^{n-1}\hookrightarrow M^n$ which is a graph over $\partial B$ in normal exponential coordinates inside a tubular neighborhood. In this subsection we consider a purely differential geometric context without needing any technical assumptions or any geometric measure theory, or any isometric embedding into Euclidean spaces. For every $y\in N$ let us define
        \begin{equation}\label{h-intrhyp} 
                        H_\nu ^{N}(y)=\underset{1}{\overset{n-1}\sum }_{i}II^N_y (e_{i,N} ,e_{i, N})=
                                                -\underset{1}{\overset{n-1}\sum }_{i} <\nabla_{e_{i, N}}\nu_N , e_{i, N} >_g (y),
        \end{equation}
        where $ (e_{1,N} ,\ldots ,e_{n-1, N} )$ is an orthonormal basis of $T_y N$,  $II^N_y(v, w)$ is the second fundamental form of $N$ at the point $y$ and evaluated on the tangent vectors $v,w\in T_y N$, $\nu_N$ is a unit normal vector field of $N$, where $\nu_N$ could be interpreted as a section of the normal bundle of $\nu N$ embedded in $TM$, and $\nabla$ is the Levi-Civita connexion of $M$. In what follow, write $\nu_N=a_N+b_N\theta$, with $a_N\perp\theta$, extend $\nu_N$ to a vector field $\nu=a+b\theta$ over an entire neighborhood $\nu:\mathcal{U}_{r}(\partial B)\rightarrow TM$ such that 
\begin{equation}
[\theta, a]=\theta(b)=0,
\end{equation} 
where $[\cdot ,\cdot]$ indicates the Lie bracket of two vector fields, and $\theta:=\nabla_g \tilde{d}_g(.,\partial B)$, the gradient of the function $\tilde{d}_g(.,\partial B)$, the signed distance function to $\partial B$ having positive values outside $B$. In particular $[\theta,\nu]=0$.\\

Let us introduce a chart $\phi$ of $M$. First, choose a chart $\Theta$ in a neighborhood of $\partial B$, and set
\Fonct{\phi}{]-r,r[\times \mathcal{U}}{\mathcal{U}_r\subseteq M}{(t,x)}{exp_{\Theta(x)}\left( t\theta(\Theta(x))\right),}
where $\mathcal{U}\subseteq\mathbb{R}^n$ is the domain of $\Theta$.
By choosing $r$ less than the normal injectivity radius of $\partial B$, we have that $\phi$ is a diffeomorphism. The functions $(t,x)$ are called \textbf{Fermi coordinates based at $\partial B$}. By Gauss's lemma, the metric $g$ of $M$ restricted to $\mathcal{U}_r$ is expressed in these coordinates as $dt^2+g_t$, where $g_t=i_t^* (g)$ and $i_t:\partial B\rightarrow M$ is the function defined as $i_t: x\mapsto\exp_{\Theta(x)}(t\theta(\Theta(x)))\in M$. In the local chart $\Theta^{-1}$, we can write $g_t=(g_t)_{ij}(x)dx^idx^j=g(t,x)_{ij}dx^idx^j$. It is useful for subsequent geometrical constructions and generalization to note that the family of embeddings $i_t$ can be interpreted as the images of $\partial B$ under the flow $\Phi_t$ of the vector field $\theta$ on $M$. Now we proceed to the explicit calculations of the mean curvature of a hypersurface $N\hookrightarrow\mathcal{U}_r$ isometrically embedded in $(\mathcal{U}_r, g_{|\mathcal{U}_r})$. Set $\nu =a+b\theta$, and note that on $N$,
\begin{equation}
\left(|a|^2+b^2\right)_{|N}=1,
\end{equation} 
but at $p\in\mathcal{U}_r\setminus N$, in general one could have 
\begin{equation}\left(|a|^2+b^2\right)_{|\mathcal{U}_r\setminus N}(p)\neq 1.
\end{equation}
As a last remark, one can see that $$H_{\nu} (t,x):=-\sum_{i=1}^{n-1}\left<\nabla_{e_i}\nu, e_i\right>,$$ is a function actually defined on all of $\mathcal{U}_{r_0}$, provided that we extend the vector fields $e_i^N$ to vector fields $e_i$ defined on $\mathcal{U}_r$, in such a way that $[e_i, \theta]=0$. Furthermore we observe that $H_{\nu} (t,x)$ is equal to $H_{\nu}^N$ when restricted to the subset $N$.
Since the trace of a linear operator is independent of the basis employed to compute it,  we will use two different basis adapted to our problem, namely $(e_1,\dots,e_{n-1}, \nu)$ and $(\partial_1,\dots,\partial_{n-1},\theta=\frac{\partial}{\partial t})$ and we obtain
\begin{eqnarray}
tr_M\left(\nabla_{(\cdot)}\nu\right) & = & \left(\sum_{i=1}^{n-1}\left\langle\nabla_{e_i}\nu,e_i\right\rangle\right)+\left\langle\nabla_{\nu}\nu,\nu\right\rangle\\
& = & \left\langle\nabla_{\partial_i}\nu,\partial_j\right\rangle g^{ij}+\left\langle\nabla_{\theta}\nu,\theta\right\rangle\\
& = & \left\langle\nabla_{\partial_i}\nu,\partial_j\right\rangle g^{ij},
\end{eqnarray}
where the Einstein summation convention is used with indexes $i$ and $j$ running on $\{1,...,n-1\}$, and by Gauss lemma $\left\langle\theta,\partial_i\right\rangle=g_{ni}=0$, which implies that $g^{ni}=0$ too. We continue the computation remarking
\begin{eqnarray}
\left\langle\nabla_{\partial_i}\nu,\partial_j\right\rangle g^{ij} & = & \left\langle\nabla_{\partial_i}a,\partial_j\right\rangle g^{ij}+\left\langle\nabla_{\partial_i}(b\theta),\partial_j\right\rangle g^{ij}, 
\end{eqnarray}
\begin{eqnarray}
\left\langle\nabla_{\partial_i}(b\theta),\partial_j\right\rangle g^{ij} & = & \left\langle\nabla_{\partial_i}(b)\theta,\partial_j\right\rangle g^{ij}+bg^{ij}\left\langle\nabla_{\partial_i}\theta,\partial_j\right\rangle,
\end{eqnarray}
but 
\begin{equation}
\left\langle\nabla_{\partial_i}(b)\theta,\partial_j\right\rangle=0,
\end{equation}
so 
\begin{equation}
\left\langle\nabla_{\partial_i}(b\theta),\partial_j\right\rangle g^{ij}=bg^{ij}\left\langle\nabla_{\partial_i}\theta,\partial_j\right\rangle.
\end{equation}
Recalling that 
\begin{equation}
g^{ij}\left\langle\nabla_{\partial_i}\theta,\partial_j\right\rangle=-H_{\theta}^{{\partial B}^t}(y),
\end{equation}
and
\begin{equation}
 \left\langle\nabla_{\partial_i}a,\partial_j\right\rangle g^{ij}=div_{\partial B^t }(a),
\end{equation}
we get
\begin{equation}
H_{\nu}^N(y)=-div_{\partial B^t }(a)+bH_{\theta}^{{\partial B}^t}+\left\langle\nabla_{\nu}\nu,\nu\right\rangle.
\end{equation}
Now it remains to examine the term
\begin{eqnarray}
\left\langle\nabla_{\nu}\nu,\nu\right\rangle & = & \left\langle\nabla_{a+b\theta}\nu,a+b\theta\right\rangle\\
 & = &  \left\langle\nabla_{a}\nu,a\right\rangle+b\left\langle\nabla_{\theta}\nu,a\right\rangle+b\left\langle\nabla_{a}\nu,\theta\right\rangle+b\left\langle\nabla_{\theta}\nu,\theta\right\rangle,
\end{eqnarray}
but
\begin{equation}
\left\langle\nabla_{\theta}\nu,\theta\right\rangle=\left\langle\nabla_{\nu}\theta,\theta\right\rangle=\frac{1}{2}\nabla_{\nu}\left\langle \theta,\theta\right\rangle=0,
\end{equation}
because $\nabla_{\theta}\nu=\nabla_{\nu}\theta$, since $[\theta, \nu]=0$, the Levi-Civita connection $\nabla$ is torsion free, and $\left\langle \theta,\theta\right\rangle=||\theta||^2\equiv 1$ on an entire neighborhood of $\partial B$. So we get
\begin{eqnarray}
\left\langle\nabla_{\nu}\nu,\nu\right\rangle & = & \left\langle\nabla_{a}a,a\right\rangle+\left\langle\nabla_{a}(b\theta),a\right\rangle+b\left\langle\nabla_{\theta}a,a\right\rangle\\
& + & b\left\langle\nabla_{\theta}(b\theta),a\right\rangle+b\left\langle\nabla_{a}a,\theta\right\rangle+b\left\langle\nabla_{a}(b\theta),\theta\right\rangle,
\end{eqnarray}
remark that
\begin{equation}
\left\langle\nabla_{\theta}(b\theta),a\right\rangle=\left\langle\nabla_{\theta}(b)\theta,a\right\rangle+b\left\langle\nabla_{\theta}(\theta),a\right\rangle=0,
\end{equation}
\begin{equation}
\left\langle\nabla_{a}(b\theta),a\right\rangle=\left\langle\nabla_{a}(b)\theta,a\right\rangle+b\left\langle\nabla_{a}(\theta),a\right\rangle=-bII_{\theta}^{{\partial B}^t} (a,a),
\end{equation}
because $\theta\perp a$, and that
\begin{equation}
b\left\langle\nabla_{\theta}a,a\right\rangle=\frac{1}{2}\nabla_{\theta}|a|^2,
\end{equation}
\begin{equation}
b\left\langle\nabla_{a}a,\theta\right\rangle=-\frac{1}{2}\nabla_{\theta}|a|^2,
\end{equation}
\begin{eqnarray}
\left\langle\nabla_{a}(b\theta),\theta\right\rangle & = & \nabla_{a}(b)\left\langle\theta,\theta\right\rangle+b\left\langle\nabla_{a}\theta,\theta\right\rangle\\
 & = & \nabla_{a}(b)+\frac{1}{2}\nabla_{a}(|\theta|^2)=\nabla_{a}(b),
\end{eqnarray}
since $|\theta|^2=1$ on $\mathcal{U}_r$.
From all this follows 
\begin{equation}
\left\langle\nabla_{\nu}\nu,\nu\right\rangle=\left\langle\nabla_{a}a,a\right\rangle-bII_{\theta}^{{\partial B}^t} (a,a)+b\nabla_{a}(b),
\end{equation}
but 
\begin{equation}
\left\langle\nabla_{a}a,a\right\rangle=\frac{1}{2}\nabla_{a}(|a|^2),
\end{equation}
\begin{equation}
b\nabla_{a}b=\frac{1}{2}\nabla_{a}(b^2),
\end{equation}
so
\begin{equation}
\left\langle\nabla_{\nu}\nu,\nu\right\rangle=-bII_{\theta}^{{\partial B}^t} (a,a)+\frac{1}{2}\nabla_{a}(|\nu|^2).
\end{equation}
Finally, for every $y\in N$ we obtain
\begin{equation}
H_\nu ^{N}(y)=-div_{\partial B^t }(a)+bH_{\theta}^{{\partial B}^t}-bII_{\theta}^{{\partial B}^t} (a,a)+\left\{\frac{1}{2}\nabla_{a}(|\nu|^2)\right\}_{|\mathcal{N}}(y).
\end{equation}
We give another way to compute of $H_\nu ^{\mathcal{N}}(y)$ in the particular case that $a\neq 0$, so that one possible choice for $e_{n-1}$, is $e_{n-1}=\frac{-b}{|a|}a+|a|\theta$, and one for $(e_1 ,\ldots ,e_{n-2})$, is $(e_1 ,\ldots ,e_{n-2})=T_y \mathcal{N}\cap T_y \partial B^t$, where $B^t$ is the domain whose boundary is the equidistant hypersurface at distance $t$ to $\partial B$.
We set $(\tilde{e}_1 =e_1 ,\ldots ,\tilde{e}_{n-1}=e_{n-2},\tilde{e}_{n-1}=\frac{a}{|a|})$.
The calculation that follows will be independent of the extensions of $\nu$, $e_i$, and thus $e_i$'s, $\nu$ can be chosen in such a way that $[e_i, \theta]=0$ for all $i\in\{1,..., n-1\}$, and $[a,\theta]=0$.\\  Note the following useful relations:
\begin{equation}
<\theta, \theta>\equiv 1, in\; \mathcal{U}_{r},
\end{equation}
\begin{equation}
\nabla_{\theta}\theta=0, in\; \mathcal{U}_{r},
\end{equation}
\begin{equation} 
<(\nabla_{e_i}b)\theta, e_i>=0, \forall i\in\{1,...,n-2\},
\end{equation}

Straightforward computations yield the equations
\begin{eqnarray}
H & = & -\sum_{i=1}^{n-1}<\nabla_{e_i }\nu , e_i>\\
    & = & -(\sum_{i=1}^{n-2}<\nabla_{e_i }\nu , e_i>)-<\nabla_{e_{n-1}}\nu ,e_{n-1}>,\label{Eq:Normalgraph}
\end{eqnarray}
expanding each term of the right hand side of the (\ref{Eq:Normalgraph}) we get
\begin{equation}
\sum_{i=1}^{n-2}<\nabla_{e_i }\nu , e_i>=\sum_{i=1}^{n-2}(<\nabla_{e_i }a , e_i>+<\nabla_{e_i }(b\theta), e_i>),
\end{equation}
\begin{equation}
<\nabla_{e_i }(b\theta), e_i>=<(\nabla_{e_i}b)\theta+b\nabla_{e_i}\theta, e_i>=b<\nabla_{e_i}\theta, e_i>, 
\end{equation}
\begin{equation}
\sum_{i=1}^{n-2}<\nabla_{e_i }\nu , e_i>=\sum_{i=1}^{n-2}<\nabla_{e_i }a , e_i>+b<\nabla_{e_i}\theta, e_i>,
\end{equation}
\begin{equation}
<\nabla_{e_{n-1}}\nu , e_{n-1}>=<\nabla_{e_{n-1}}a , e_{n-1}>+<\nabla_{e_{n-1}}(b\theta) , e_{n-1}>,
\end{equation}
\begin{eqnarray}
<\nabla_{e_{n-1}}(b\theta) , e_{n-1}> & = & <(\nabla_{e_{n-1}}b)\theta , e_{n-1}>\\
 & + & b<\nabla_{e_{n-1}}\theta, e_{n-1}>,
\end{eqnarray}
\begin{eqnarray} 
<(\nabla_{e_{n-1}}b)\theta, e_{n-1}> & = & (\nabla_{e_{n-1}}b)<\theta, |a|\theta>\\
& = & -\frac{b}{|a|}|a|\nabla_ab=-b\nabla_ab,
\end{eqnarray}
the final formula, coming from the preceding equalities, is the following
\begin{equation}
\nabla_{e_{n-1}}b=\nabla_{-\frac{b}{|a|}a+|a|\theta}b=-\frac{b}{|a|}\nabla_ab+|a|\nabla_{\theta}b=-\frac{b}{|a|}\nabla_ab.
\end{equation}
\begin{eqnarray}
b<\nabla_{e_{n-1}}\theta, e_{n-1}> & = &b<-\frac{b}{|a|}\nabla_a\theta+|a|\nabla_{\theta}\theta, -\frac{b}{|a|}a+|a|\theta>\\
 & = & b\left(\frac{b^2}{|a|^2}\right)<\nabla_a\theta, a>-b^2<\nabla_a\theta, \theta>\\
 & = & b\left(\frac{1}{|a|^2}-1\right)<\nabla_a\theta, a>\\
 & = & b<\nabla_{\tilde{e}_{n-1}}\theta, \tilde{e}_{n-1}>-b<\nabla_a\theta, a>.
\end{eqnarray}
Thus
\begin{eqnarray}
<\nabla_{e_{n-1}}(b\theta) , e_{n-1}> & = & b<\nabla_{\tilde{e}_{n-1}}\theta, \tilde{e}_{n-1}>-b<\nabla_a\theta, a>\\
 & - & b\nabla_ab, \nonumber
\end{eqnarray}
\begin{eqnarray}
<\nabla_{e_{n-1}}a , e_{n-1}> & = & \frac{b^2}{|a|^2}<\nabla_a a , a>\\
& = & \left(\frac{1}{|a|^2}-1\right)<\nabla_a a , a>\\
& = & <\nabla_{\tilde{e}_{n-1} }a , \tilde{e}_{n-1}>-<\nabla_a a , a>,
\end{eqnarray}
\begin{eqnarray}
    <\nabla_{e_{n-1}}\nu ,e_{n-1}> & = & <\nabla_{\tilde{e}_{n-1} }a , \tilde{e}_{n-1}>-<\nabla_a a , a>\\ \nonumber
     & + & b<\nabla_{\tilde{e}_{n-1}}\theta, \tilde{e}_{n-1}>-b<\nabla_a\theta, a>\\ \nonumber
     & - & b\nabla_ab.
\end{eqnarray}
Hence
\begin{eqnarray}\label{cmpphyp-}
       -\sum_{i=1}^{n-1}<\nabla_{e_i }\nu , e_i> & = & -\sum_{i=1}^{n-2}<\nabla_{e_i }a,e_i >-\sum_{i=1}^{n-2}b<\nabla_{e_i }\theta,e_i>\\ \nonumber
       & - & <\nabla_{e_{n-1}}\nu ,e_{n-1}>\\ \nonumber
       & = & -\sum_{i=1}^{n-2}<\nabla_{e_i }a,e_i >-\sum_{i=1}^{n-2}b<\nabla_{e_i }\theta,e_i>\\ \nonumber
       & - & <\nabla_{\tilde{e}_{n-1} }a , \tilde{e}_{n-1}>-b<\nabla_{\tilde{e}_{n-1}}\theta, \tilde{e}_{n-1}>\\ \nonumber
       & + & b<\nabla_a\theta, a> +\; b\nabla_ab \;+ <\nabla_aa,a>.\nonumber
\end{eqnarray}
Before to give the final formula we observe  
\begin{eqnarray}
d\frac{1}{2}|\nu|^2(a) & = & \nabla_{a}\frac{1}{2}|\nu|^2\\ \nonumber
 & = & \frac{1}{2}<\overrightarrow{\nabla}_g|\nu|^2, a>\\ \nonumber
 & = & \nabla_{a}\left[\frac{1}{2}(|a|^2(t,x)+b^2(t,x))\right](t,x)\\ \nonumber 
 & = & b_2\nabla_{a} b+<{\nabla}_{{a}}a, a>.
\end{eqnarray}
Finally we have again  
\begin{equation}\label{cmpphyp}
      H_{\nu} (t,x)=-div_{\partial B^t }(a)+bH_{\theta}^{{\partial B}^t}-bII_{\theta}^{{\partial B}^t} (a,a)+\frac{1}{2}\nabla_{a}|\nu|^2,
\end{equation} 
where $II_{\theta}^{{\partial B}^t}$ and $H_{\theta}^{{\partial B}^t} $ are respectively the second fundamental form and the mean curvature in the direction of  $\theta$ of the equidistant hypersurface at distance $t$ from $\partial B$ computed at the point $exp^{\partial B}(t\theta(\Theta(x)))\in N$. Equation (\ref{cmpphyp}) comes from a geometrical interpretation of the terms in (\ref{cmpphyp-}), by observing that
\begin{enumerate}[(i):]
          \item $div_{\partial B_t }(a)=\sum_{i=1}^{n-2}<\nabla_{e_i }a,e_i >+<\nabla_{\tilde{e}_{n-1} }a, \tilde{e}_{n-1}>$,
          \item $bH_{\theta}^{{\partial B}^t}=-\left[\sum_{i=1}^{n-2}b<\nabla_{e_i }\theta,e_i>+b<\nabla_{\tilde{e}_{n-1}}\theta, \tilde{e}_{n-1}>\right]$,
          \item $bII_{\theta}^{{\partial B}^t} (a,a)= -b<\nabla_a\theta, a> $.
\end{enumerate}
As a last remark, one can see that $H_{\nu} (t,x)$ is a function actually defined on all of $\mathcal{U}_{r_0}$, and it is equal to $H^N$ when restricted to the subset $N$.
Another interesting and simpler formula is obtained by choosing extensions $\nu_1$, $a_1$, $b_1$, respectively of $\nu_N$, $a_N$, and $b_N$, in such a way that 
\begin{equation}
\nabla_{\theta}\nu_1=0,
\end{equation} 
is satisfied on the entire $\mathcal{U}_{r_0}$. The same kind of computation leading to (\ref{cmpphyp}) leads to the following formula for the mean curvature:
\begin{equation}\label{Eq:parallelmc}
H_{\nu_1}(t,x)=-div_{\partial B^t }(a_1)+b_1H_{\theta}^{{\partial B}^t}.
\end{equation}
This latter formula is the analog of formula (\ref{cmpphyp}). It is easy to check that ${H_{\nu_1}}_{|N}={H_{\nu}}_{|N}$, but in general ${H_{\nu_{1}}}_{|(\mathcal{U}_{r_0}\setminus N)}\neq{H_{\nu_2}}_{|(\mathcal{U}_{r_0}\setminus N)}$ for every fixed $N\subset\mathcal{U}_{r_0}$. Another very simple way to prove (\ref{Eq:parallelmc}) is to observe that $H$ is the trace of an appropriate endomorphism and computing with respect to two different choices of orthonormal basis. First observe that, by construction,
\begin{equation}
<\nu_1,\nu_1>=1,
\end{equation} 
is a constant function on all $\mathcal{U}_{r_0}$, so for an arbitrary $p\in M$ we have
\begin{equation}
\nabla_{X_p}(<\nu_1,\nu_1>)=0, \forall X_p\in T_p\mathcal{U}_{r_0}.
\end{equation}
In particular
\begin{eqnarray}
\nabla_{\nu_p}(<\nu_1,\nu_1>)=0, \forall p\in\mathcal{N},\\
\nabla_{\theta_p}(<\nu_1,\nu_1>)=0, \forall p\in\mathcal{U}_{r_0}.
\end{eqnarray}
Now we are ready to prove (\ref{Eq:parallelmc}) as follows. We have that
\begin{equation}
H_{|\mathcal{N}}=-div_{\mathcal{N}}\nu=-div_{\mathcal{N}}(\nu_1)_{|\mathcal{N}},
\end{equation}  
because of the independence from extension of $\nu$ defined only on $\mathcal{N}$, to any $\nu_1$ defined on all $\mathcal{U}_{r_0}$.
The divergence of the vector field $\nu_1$ calculated in the orthonormal frame $(\tilde{e}_1,...,\tilde{e}_{n-1},\theta)$ is
\begin{eqnarray}
div_{M}(\nu_1)_{|\mathcal{N}} & = & \sum_{i=1}^{n-1}<\nabla_{\tilde{e}_i}\nu_1,\tilde{e}_i>+<\nabla_{\theta}\nu_1,\theta>\\ \nonumber
             & = & div_{\mathcal{N}}\nu_1\\ \nonumber
             & = & div_{\mathcal{N}}\nu.
\end{eqnarray}
But
\begin{equation}
          div_{M}(\nu_1)_{|\mathcal{N}}= div_{\partial B^t }(a_1)-b_1H_{\theta}^{{\partial B}^t}.
\end{equation}
Combining the last three equations, it is easy to check the validity of (\ref{Eq:parallelmc}).
Assume now that $N=\{p\in\mathcal{U}_{r_0}|\exists x\in\mathcal{U}, s.t.\; p=exp_p(u(x)\theta(\Theta(x)))\}$ is the normal graph (i.e., in normal exponential coordinates) of a function $u\in C^{2,\alpha}(\partial B)$. 
 Let
$$W_u :=\sqrt{1+\|\overrightarrow{\nabla}_{i_u ^* (g)}u \| _{i_u ^* (g)}^2 }.$$
Here we consider $W_u$ as both a function on $M$ independent of $t$ and a function on $\partial B$, and so we will make no distinction between $W_u$ and $\pi_{\partial B}^*(W_u)$, where $\pi_{\partial B}$ denotes the projection $\pi_{\partial B}:]\varepsilon, \varepsilon[\times\partial B\rightarrow \partial B$,  $\pi_{\partial B}:(t,x)\mapsto x$.
Then 
$$b= \left\{ 
         \begin{array}{llll}
                      \frac{1}{W_u }, & <\nu , \theta >\geq 0, & \nu & outward\\
                     -\frac{1}{W_u }, & <\nu , \theta >\leq 0, & \nu & inward\\
         \end{array}\right. 
         $$
Let $b=\frac{1}{W_u }$. 
In Fermi coordinates, the preceding equation (\ref{cmpphyp}) can be written as
\begin{equation}
a(t, x)=-\frac{1}{W_u }g(u(x),x)^{ij}\frac{\partial u}{\partial x^j}(x)\frac{\partial}{\partial x^i}(t, x).
\end{equation}  

\noindent This leads to
\begin{eqnarray}
div_{\partial B^t}(a) & = & -div_{(\partial B,g_t)}(\frac{\overrightarrow{\nabla}_{(\partial B,g_u)} u}{W_{u}}) \nonumber \\
& = & -\frac{1}{\sqrt{det(g_t)}}\frac{\partial}{\partial x^i}\left(\sqrt{det(g_t)}\frac{1}{W_u }g_u^{ij}\frac{\partial u}{\partial x^j}(x)\right).
\end{eqnarray}
We observe here that in general $g_t=g_t(x)$ depends on $x$, although it is independent in some important cases, including warped product normal bundles.
The mean curvature of the graph of $u$ is thus
\begin{eqnarray}\label{h-norm1hyp}  
                   \H_g[u](x) & = & \left(div_{(\partial B^t)}(\frac{\overrightarrow{\nabla}_{g_u} u}{W_{u}})\right)_{|t=u(x)}\\  \nonumber                             
                           & - &\frac{1 }{W_{u} ^3}II_\theta ^{u} (\overrightarrow{\nabla}_{g_u} u,
                                     \overrightarrow{\nabla}_{g_u} u)\\ \nonumber
                            & + & \frac{1}{W_{u}}H_\theta ^{u} (u)\\ \nonumber
                            & - & \frac{1}{W_u^2}g_u^{ij}\frac{\partial u}{\partial x^j}(x)g_u^{lm}\frac{\partial u}{\partial x^m}(x)<\nabla_{\frac{\partial}{\partial x^i}(t, x)}\left(\frac{\overrightarrow{\nabla}_{g_u} u}{W_{u}}\right),  \frac{\partial}{\partial x^m}(t, x)>_{|t=u(x)}\\ \nonumber
                            & - & \frac{1}{W_u^2 }g_u^{ij}\frac{\partial u}{\partial x^j}(x)\left(\nabla_{\frac{\partial}{\partial x^i}(t, x)}\left(\frac{1}{W_u }\right)\right)_{|t=u(x)}.                                
\end{eqnarray}
with $<\nu_{ext},\theta>_g\ge 0$ and $\overrightarrow{\nabla}_{g_u} u=g_{u(x)}^{ij}\frac{\partial u}{\partial x^j}(x)\frac{\partial}{\partial x^i}(t, x)$. 
Regarding $\H$ as an opertor $\H_g:C^{2,\alpha}(\partial B)\rightarrow C^{0,\alpha}(\partial B)$, we easily see that it is semilinear elliptic, which is essentially the only property of $\H_g$ we will use in this paper. But the exact expression (\ref{h-norm1hyp}) for $\H$ demonstrates that the coefficients of the constant mean curvature equation 
\begin{equation}\label{cmcequation}
\H_g[u]=const,
\end{equation}
are bounded in various Sobolev and H\"older spaces. As a result, one can apply standard bootstrap arguments of elliptic regularity theory to show higher order regularity of solutions $u$ of the constant mean curvature equation (\ref{cmcequation}). To obtain the estimates needed to apply elliptic regularity theory, one need not appeal to (\ref{h-norm1hyp}). In fact, this is an immediate consequence of the definition of the mean curvature function as the partial divergence with respect to $TN$ of the smooth vector field $\nu$, i.e., is a general operator in divergence form to which classical result applies. The interest in a formula like (\ref{h-norm1hyp}) is more geometric and lies in the possibility of applying (\ref{h-norm1hyp}) to solve (\ref{cmcequation}) and to help to give a qualitative description of solutions knowing the geometry of the equidistant foliation of $\mathcal{U}_r$, in the ambient manifold. One instance of this philosophy can be found in \cite{NarAnn}.  
\subsection{Allard's Regularity Theorem}

The proof of the Theorem \ref{tr2} is mainly based on a regularity theorem for almost minimizing varifolds. In the article \cite{All}, it is stated in an Euclidean ambient context. Using isometric embeddings we can deduce a Riemannian version of it.

\bigskip

We restate, here, for completeness sake, the regularity theorem of Chapter $8$ p. 466 of \cite{All} that will be of frequent use in the sequel. For this statement we use the notations of the original article  \cite{All}. 
\begin{Def}
For any $0\le m\le n$, we say that $V$ is a \textbf{$m$-dimensional varifold in $M$}, if $V$ is a nonnegative real extended valued $($compare section $2.6$ of $\cite{All}$$)$ Radon measure on $G_m(M)$ the Grassmannian manifold whose underlying set is the union of the sets of $m$-dimensional subspaces of $T_xM$ as $x$ varies on $M$. For every $m\in\{0,...,n\}$, we define $\mathbf{V}_m(M)$ to be the space of all $m$-dimensional varifolds on $M$ endowed with the weak topology induced by $C^0_c(G_m(M))$ say the space of continuous compactly supported functions on $G_m(M)$ endowed with the compact open topology.  
\end{Def}
\begin{Def}
Let $V\in\mathbf{V}_m(M)$, $g$ is a Riemannian metric on $M$, we say that the nonnegative Radon measure on $M$, $||V||$ is the \textbf{weight} of $V$, if $||V||=\pi_{\#}(V)$, here $\pi$ indicates the natural fiber bundle projection $\pi:G_m(M)\rightarrow M$, $\pi:(x, S)\mapsto x$, for every $(x,S)\in G_m(M)$, $x\in M$, $S\in G_m(T_xM)$,
 $$||V||(A):=V(\pi^{-1}(A)).$$ 
\end{Def} 
Notice that the notion of a varifold is independent of the choice of any Riemannian metric $g$ on $M$. This reflects the phenomenon that on a differentiable manifold one can have a fixed submanifold but whose metric datas like volume, curvature, second fundamental form, etc. depends on the metric that we put on it. If we consider a varifold $V\in\mathbf{V}_m(M)$ we can construct without the help of a metric the support of $||V||$ that is a set contained in $M$, however starting from a set $E\subseteq M$ even a good one like a $m$-dimensional smooth submanifold of $M$, there is no canonical way to come back to a uniquely determined varifold $V\in\mathbf{V}_m(M)$, such that $Supp||V||=E$. One way to proceed is to chose a metric $g$ and to associate to a $\mathcal{H}^m_g$-countably $m$-rectifiable set $E$, the varifold $V_g(E)\in\mathbf{V}_m(M)$, where
\begin{equation}
V(E,g)(A):=\mathcal{H}_g^m(\{x\in E:(x, T_xE)\in A\}),\:\forall A\in G_{m}(M^n),
\end{equation} 
in this way the manifold associated is unique and canonical in the sense that depends only on the choice of the metric $g$. When $(M^n,g)$ is $(\R^n, \xi)$ we find again the classical theory of varifolds as developed by Almgren, Allard et al. The way in which classically one proceed to study the theory of varifolds in Riemannian manifolds is well explained in \cite{All} and consists in embedding isometrically $(M^n,g)$ in some higher dimensional Euclidean space via Nash's Theorem, and then using the existing theory on $\R^n$ of \cite{All}. The point of view that we will adopt here is an intrinsic one, without having to choose an isometric embedding. This is needed because in the Euclidean monotonicity formula will appear an upper bound of the second fundamental form of the particular isometric embedding chosen and it is not clear to us how to bound the second  fundamental form of the isometric embedding just starting with intrinsic bounded geometry assumptions on the manifold $(M,g)$. The intrinsic approach avoid this technical difficulty and permits to have a monotonicity formula which  depends only on an upper bound of the sectional curvature. This means that locally the geometric measure theory of $\R^n$ is mutatis mutandis the same as the corresponding theory developed on a Riemannian manifold, with just the constants involved depending on the bound of the sectional curvature. This is what one could expects since locally a Riemannian manifold is bi-Lipschitz equivalent to an Euclidean ball via the exponential map. The importance of making rigorous the details and the proofs appears clear when we deal with problems in ambient manifolds with variable metric. Let us introduce at this point some standard notions that will be useful in our further developments.       
\begin{Def}\label{Def:Density} Let $\mu$ be a Borel regular measure on a locally compact Hausdorff topological space $X$. Define
$$\Theta_*^m(\mu, a):=\liminf_{r\rightarrow 0^+}\frac{\mu(B(a,r))}{\omega_mr^m},$$ the \textbf{$m$-lower density of $\mu$ at $a\in M$},
$$\Theta^{*m}(\mu, a):=\limsup_{r\rightarrow 0^+}\frac{\mu(B(a,r))}{\omega_mr^m},$$ the \textbf{$m$-upper density of $\mu$ at $a\in M$}, and if
$$\Theta_*^m(\mu, a)=\Theta^{*m}(\mu, a),$$ then we set $$\Theta^m(\mu, a):=\Theta^{*m}(\mu, a)=\Theta_*^m(\mu, a)=\lim_{r\rightarrow 0^+}\frac{\mu(B(a,r))}{\omega_mr^m}.$$ We call $\Theta^m(\mu, a)$ the \textbf{$m$-density of $\mu$ at $a\in X$}.  
\end{Def}
According to \cite{All} we give the following definition for the first variation of a varifold.
\begin{Def}
Let $V\in\mathbf{V}_m(M)$. Let $\mathfrak{X}_c(M)$ denotes the set of smooth vector fields on $M$ with compact support, we denote by the linear function $\delta_g V(X):\mathfrak{X}_c(M)\rightarrow\mathbb{R}$, the \textbf{first variation} of the varifold $V$ in the direction of the vector field $X\in\mathfrak{X}_c(M)$, defined as follows
\begin{eqnarray}\label{Eq:DefofFirstVariation}\nonumber
\delta_gV(X) & := & \int_{\xi\in G_m(M)}\left<(\nabla^g X(\pi(\xi))\circ \pi_S), \pi_S\right>_gdV(\xi)\\ 
& := & \int_{\xi\in G_m(M)}\sum_{i=1}^n\left<\nabla^g_{\pi_S(e_i)}X, \pi_S(e_i)\right>_gdV(\xi)\\
& := & \int_{\xi\in G_m(M)}div_{S,g}XdV(\xi),
\end{eqnarray}
for every $X\in\mathfrak{X}_c(M)$, where $S\leq T_xM$ is such that $\xi=(x,S)\in G_m(M)$, i.e., a $m$-dimensional subspace of $T_xM$, $\pi_S$ is the orthogonal projection $\pi_S:T_xM\to S$ with respect to the metric $g$, $(e_1,..., e_n)$ is an orthonormal basis of $(T_{\pi(\xi)}M,g_{\pi(\xi)})$, and $div_{S,g}X=\sum_{i=1}^m \langle \nabla_{\tilde{e}_i}X, \tilde{e}_i\rangle_g$, with $\{\tilde{e}_1,\dots,\tilde{e}_m\}$ being an orthonormal basis over $S$ with respect to $g$.
\end{Def}
\begin{Rem} The first variation is a metric concept and depends on $g$.
\end{Rem}
\begin{Rem} In the rest of this paper we adopt the convention to denote real variables with letters without subscripts and real constants by letters with subscripts.
\end{Rem}
\begin{Thm}[Allard's Regularity Theorem $8.1$ \cite{All}, Euclidean version ]\label{alleucl}
Let $p>1$ be a real number. Let $q$ be the conjugate exponent of $p$, i.e., $q$ satisfies $\frac{1}{p}+\frac{1}{q}=1$. Let $k$ be a integer number, $1\leq k\leq n$. We assume that $k<p<+\infty$, if $k>1$, and that $p\geq 2$, if $k=1$.

For all $\varepsilon\in ]0,1[$ there exists $\eta_1=\eta_1(\varepsilon) >0$, (that depends on $\varepsilon$) such that for all reals $R>0$, for all integer $d\geq 1$, for all varifolds $V\in V_k (\mathbb{R}^n )$ and for all points $a\in spt||V||$, if 
      \begin{enumerate}
            \item $\Theta^k (||V||,x)\geq d$ for $||V||$ almost all $x\in B_{\mathbb{R}^n}(a,R)$; 
            \item $||V||(U(a,R))\leq (1+\eta_1 )d\omega_k R^k$;
            \item $\delta_g V(X)\leq\eta_1 d^{\frac{1}{p}}R^{\frac{k}{p}-1}\left(\int_{\mathbb{R}^n}|X|^q ||V||(dx)\right)^{\frac{1}{q}}$,
                  with $X\in\mathfrak{X} (\mathbb{R}^n)$ and $supp(X)\subset U(a,R):=\left\{x\in\R^n\:|\: 0\leq|x-a|<R\right\}$. 
      \end{enumerate}
      Then there exists a map $F_1:\mathbb{R}^k\rightarrow\mathbb{R}^n$ such that
      \begin{enumerate}
            \item $F_1\in C^1 (\mathbb{R}^k,\mathbb{R}^n)$ and $T\circ F_1=Id_{\R^k}$, where $T:\mathbb{R}^n \to\mathbb{R}^k$ is an orthogonal projection,
            \item $U(a,(1-\varepsilon)R)\cap spt||V||=U(a,(1-\varepsilon )R)\cap F_1(\mathbb{R}^k )$,
            \item $\forall y,z\in \mathbb{R}^k$, $||dF_1(y)-dF_1(z)||\leq\varepsilon \left(\frac{|y-z|}{R}\right)^{1-\frac{k}{p}}$.
            \item Moreover $\eta_1(\varepsilon)\to0^+$, when $\varepsilon\to 0^+$.
      \end{enumerate}
\end{Thm}

\begin{Thm}[Allard's Regularity Theorem, Riemannian version]\label{allriem}
      Let $(M^n,g)$ be a compact Riemannian manifold, $i_g:M\hookrightarrow\mathbb{R}^N$ be an isometric embedding. Let $p>1$ be a real number. Let $q$ be the conjugate exponent, $\frac{1}{p}+\frac{1}{q}=1$. Let $k$ be an integer number, $1\leq k\leq n$. We assume that $k<p<+\infty$ if $k>1$, and that $p\geq 2$, if $k=1$.

For all $\varepsilon\in ]0,1[$ there exists $\tilde{\eta}_1=\tilde{\eta}_1(\varepsilon, i_g)>0$, such that there exists a $\tilde{R}_1=\tilde{R}_1(\varepsilon, i_g)=\tilde{R}_1(\varepsilon, g,\partial g, \partial^2 g, \partial^3 g, \partial^4 g)>0$ satisfying the property that for all $0<\tilde{R}\leq \tilde{R}_1$, for all integer number 
$0<\tilde{d}<+\infty$, for all varifolds $V\in V_k (M^n )$, and for all point $a\in spt||V||$, if 
      \begin{enumerate}
            \item $\Theta^k (||V||,x)\geq \tilde{d}$ for $||V||$ almost every $x\in B_{M}(a,\tilde{R})$; 
            \item $||V||(B(a,\tilde{R}))\leq (1+\tilde{\eta}_1)\tilde{d}\omega_k \tilde{R}^k$,
            \item $\delta_g V(X)\leq\tilde{\eta}_1 \tilde{d}^{\frac{1}{p}}\tilde{R}^{\frac{k}{p}-1}\left(\int_{M}|X|_g^q ||V||(dx)\right)^{\frac{1}{q}}$, with $X\in\mathfrak{X} (M)$ and $supp(X)\subset B(a,\tilde{R})$, 
      \end{enumerate}
      then there exists a function $\tilde{F}_1 :\mathbb{R}^k\rightarrow M$, $R_0=R_0(i_g,\tilde R,\varepsilon)<\tilde{R}$ $($$\tilde{F}_1$ and $R_0$ are mutually independent$)$ such that      
      \begin{enumerate}
            \item $\tilde{F}_1\in C^1 (\mathbb{R}^k,M)$, $d\tilde{F}_1(0)$ is an isometry,\\
            \item $B(a,(1-\varepsilon) R_0)\cap spt||V||=B(a,(1-\varepsilon)R_0)\cap \tilde{F}_1(\mathbb{R}^k )$,\\ 
        \item $||d\tilde{F}_1(y)-d\tilde{F}_1(z)||\leq\varepsilon \left(\frac{|y-z|}{R_0}\right)^{1-\frac{k}{p}}$ for all $y,z\in \mathbb{R}^k$.\\
        \item $\tilde{\eta}_1\to0$, when $\varepsilon\to0^+$.
      \end{enumerate}
\end{Thm}
\begin{Rem} $R_0$ is independent of $V$.
\end{Rem}
\textbf{Remark:}  In the statement of the theorem the constant $\tilde{\eta}_1$ depends on the embedding $i$ and on $\eta_1$ produced by Theorem \ref{alleucl}.

\textbf{Idea of the proof:}
      At this point we try to apply Theorem \ref{alleucl}  to the varifold $i_{\#}(V)$.
Actually, if $V$ satisfies the assumptions $1$, $2$ and $3$ of Theorem \ref{allriem}, then $i_{\#}(V)$ satisfies  the hypothesis of  Allard's Regularity Theorem, Euclidean version (see Theorem \ref{alleucl}) but, with different constants.\\
\bigskip

To this aim, we need to compare the intrinsic distance of a submanifold and the distance of the ambient manifold restricted to the submanifold.

\begin{Lemme}\label{dgeod}
            Let $M$ be an embedded manifold into $\mathbb{R}^N$ of arbitrary codimension.
            $i:M\hookrightarrow\mathbb{R}^N$ an isometric embedding and $\beta_i$ its second fundamental form. Fix a point $a\in\mathbb{R}^N$, $a\in M$ and consider a second point $y\neq a$ different from $a$ on $M$, now take the geodesic $\sigma$ of $M$ of length $\tilde{R}$ that joins $a$ and $y$ on $M$ and the Euclidean segment $[a,y]$ of length $R$. Then there exists $R_0>0$ and a constant $\delta_i >0$ depending only on $\beta_i $ and $R_0$ such that for all $R<R_0$, results 
$\tilde{R}\leq R(1+\delta_i R^2)$. 
\end{Lemme}
\begin{Dem}
We take as origin of coordinates the point $a$ and parametrize $\sigma$ by its arc length $s$. Consider the function $f(s)=|\sigma (s)|^2$. Then
           $f(\tilde{R})=R^2$, $f'(s)=2<\sigma ',\sigma >(s)$, 
$$ \begin{array}{ccl}
           f''(s) & = & 2(<\sigma '',\sigma >(s)+<\sigma ',\sigma '>(s))\\
                   & = & 2+2<\sigma '',\sigma >(s)\\
                   & = & 2+2<\sigma '',\sigma -s\sigma '>(s).
      \end{array}$$
Since $(\sigma -s\sigma ')'=\sigma ' -\sigma ' -s\sigma ''$, $||(\sigma -s\sigma ')'||\leq s||\sigma ''||\leq s\beta_i$, we get $||\sigma -s\sigma '||\leq \frac{s^2}{2}\beta_i $. It follows that
$f''(s)\geq 2-s^2\beta_i ^2$, $f'(s)-f'(0)=f'(s)\geq 2s-\frac{s^3}{3}\beta_i ^2$, $f(s)\geq s^2-\frac{s^4}{12}\beta_i ^2 $, which implies  
\begin{equation}
       f(\tilde{R})=R^2\geq \tilde{R}^2-\frac{\tilde{R}^4}{12}\beta_i ^2 .
\end{equation} 
Finally $\tilde{R}\leq R(1+\frac{R^2 const.}{24}\beta_i ^2)=R(1+R^2\delta_i)$ where $\delta_i$ is a constant that depends only on $\beta_i$ and $R_0$. More explicitely could be taken as $\delta_i=\sum_ja_n(\beta_i){R_0}^j$, for some positive general terms $a_n$ that depends only on $\beta_i$. 
\end{Dem}

\begin{Dem}[Proof of Allard's Regularity Theorem, Riemannian version]
      In this context, variables and constants respect the previous convention and furthermore constants and variables relative to intrinsic objects of $M$ are denoted with a tilda. 
      From the following formula [4.4 (1) in \cite{All}]:
      \begin{equation}\label{allcm}
            \delta(i_{g,\#}V)(X)=\delta_gV(X^{\top})-\int_{G_k (M)}X^{\bot}(x)\cdot h(M,(x,S))dV(x,S), 
      \end{equation} 
      with $X\in\mathfrak{X}_c(U(a,R_0))$, $X(x)=X^{\top}(x)+X^{\bot}(x)$, $X^{\top}(x)\in T_x M$, $X^{\bot}(x)\in T^{\bot}_x M$,
      we can deduce that assumption $3$ of Theorem \ref{alleucl} is satisfied with some suitably chosen constant.  To see this, it is sufficient to control the Euclidean mean curvature of $i_{g,\#}V$.\\
      Now, we assume that $R_0$, $\tilde{\eta}_1$, $\tilde{R}$ verify the following conditions:
                     \begin{small}
                     \begin{equation}
            0<R_0<\min\left\lbrace \inf_{x\in i_g(\partial B(a,\tilde{R}))}\left\{|x-a|_{\mathbb{R}^N}\right\},\sqrt{\frac{(1+\eta_1(\varepsilon))^{\frac{1}{k}}-1}{\delta_i }}\right\rbrace,
            \end{equation} 
            \end{small}
       $\tilde{d}=d$,
            \begin{equation}
            0<\tilde{\eta}_1(\varepsilon)\leq\min \left\lbrace \frac{\eta_1}{2},\frac{1+\eta_1}{(1+\tilde{\delta}_i R_0 ^2 )^k}-1\right\rbrace,
            \end{equation}
            \begin{equation}
            0<\tilde{R}\leq\frac{\tilde{\eta}_1(\varepsilon)}{\beta_i(1+\tilde{\eta}_1(\varepsilon))^{\frac{1}{p}}\omega_k ^ {\frac{1}{p}}}=:\tilde{R}_1(\varepsilon). \end{equation}
      \begin{Rem} First we choose $R_0 >0$, then $\tilde{\eta}_1$ and after that, $\tilde{R}_1$ with dependences in this order. \end{Rem}
      \begin{Rem} $\tilde{\eta}_1(\varepsilon)\to0$, $\tilde{R}_1(\varepsilon)\to0$, when $\varepsilon\to\varepsilon$.
      \end{Rem}
                                       The condition $0<R_0<\sqrt{\frac{(1+\eta_1)^{\frac{1}{k}}-1}{\delta_i }}$ serves to assert that $\frac{1+\eta_1}{(1+\delta_i R_0^2 )^k}-1 >0$ and there exists  $\tilde{\eta}_1$ such that $(1+\tilde{\eta}_1)\omega_k \tilde{R}^k\leq (1+\eta_1 )\omega_k R^k$.\\  
The condition $0<R_0<\inf_{x\in i_g(B(a,\tilde{R}))}\left\{|x-a|_{\mathbb{R}^N}\right\}$ serves to assert that\\ $spt||i_{g,\#}V||\cap i_g(B(a,R_0))\subseteq i_g(spt||V||\cap B(a,\tilde{R}))$.\\
From what is said, it follows  
\begin{equation}
||i_{g,\#}V||(B_{\mathbb{R}^N}(a, R_0))\leq ||V||(B_{M}(a,\tilde{R}))\leq d(1+\tilde{\eta}_1)\omega_k \tilde{R}^k\leq d(1+\eta_1 )\omega_k R_0 ^k .
\end{equation}

The first term on the right hand side of equation (\ref{allcm}) is estimated thanks to assumption 3,
$$
|\delta_g V(X^\top)|\leq \tilde{\eta}_1 d^{\frac{1}{p}}\tilde{R}^{\frac{k}{p}-1}\left(\int_{M}|X^{\top}|^q ||V||(dx)\right)^{\frac{1}{q}}\leq\tilde{\eta}_1 d^{\frac{1}{p}}\tilde{R}^{\frac{k}{p}-1}||X||_{L^q (||V||)}.
$$
To the second term, we apply H\"older's inequality,
$$
|\int_{G_k (M)}X^{\bot}(x)\cdot h(M,(x,S))dV(x,S)|\leq \beta_i \left\lbrace \int_{Supp(X)}d||V||\right\rbrace ^{\frac{1}{p}}||X||_{L^q (||V||)}.
$$
Choosing vector fields $X$ supported in the $R_0$-ball makes 
$$
\left\lbrace \int_{Supp(X)}d||V||\right\rbrace ^{\frac{1}{p}}\leq \left\lbrace ||i_{g,\#}V||(B(a,R_0 ))\right\rbrace ^{\frac{1}{p}}\leq d^{\frac{1}{p}}(1+\eta_1 )\omega_k R_0 ^k .
$$
It follows that
\begin{equation}
          \begin{array}{lll}
\delta(i_{g,\#}V)(X) & \leq & \eta d^{\frac{1}{p}}R_0 ^{\frac{k}{p}-1}\left(\int_{\mathbb{R}^n}|X|^q ||V||(dx)\right)^{\frac{1}{q}}.                      
          \end{array} 
\end{equation}
Now we can apply Theorem \ref{alleucl} (Allard's Euclidean) to $i_{g,\#}V$ at point $a$ with $R=R_0$ as described previously to obtain (with a little abuse of notation for $i_g^{-1}$), $\tilde{F}_1=i_g^{-1}\circ F_1$ where $F_1$ is given by Theorem \ref{alleucl} (Allard Euclidean).
It can be easily seen that $dF_1(0)=Id$ and that  $i$ is an isometric embedding. 
This implies that $d\tilde{F}_1(0)$ is an isometry.
\end{Dem}
\subsection{First Variation of isoperimetric regions}

In this subsection, we check that varifold isoperimetric regions have constant mean curvature. This will be used later, in Lemma \ref{bornecm}, where Levy-Gromov's inequality will be used to verify the third assumption in Allard's theorem.

 \begin{Lemme}\label{vpremiere} Let $(M^n,g)$ be a smooth Riemannian manifold. Let $V$ be the varifold associated to a current $\partial D$ of dimension $n-1$, that is the boundary of an isoperimetric region $D$. 
             Then there exists a constant $H_g$ so that for every vector field $X\in\mathfrak{X}(M)$ we have 
             $$\delta_g\partial D(X)=-H_g\int_{Spt ||\partial D||}<X,\nu>_g||\partial D||(x),$$ 
             where $\nu$ is the outward normal to the boundary of $D$ defined $||\partial D||$-a.e.  
      \end{Lemme}
\begin{Rem} We observe that it is the first time that we encounter in this paper a concept used to study a varifold depends on the metric, namely the mean curvature.
\end{Rem}
\begin{Rem} The vector $\nu$ is the exterior normal to $\partial D$ that by regularity theory exists $\mathcal{H}_g^{n-1}$ a.e. on $\partial D$. (The reader can consult \cite{Morg1}) 
\end{Rem}
      \begin{Dem}
            As $\mathfrak{X}(M)$ is the space of sections of the tangent bundle $TM\rightarrow M$,
            it has a natural structure of vector space (possibly of infinite dimension).
            Consider the following linear functionals on this vector space:
 \Fonct{Flux_g}{\mathfrak{X}(M)}{\mathbb{R}}{X}{\int_{\partial D}<X,\nu>_gdA_{g,\partial D}(x)}
                        \Fonct{\delta_g\partial D}{\mathfrak{X}(M)}{\mathbb{R}}{X}{\delta_g\partial D(X)}
            \begin{Lemme}\label{flux}
                  If $Flux_g(X)=0$, then there exists a variation $h(t,x)$ such that $\mathbf{M}_g((h_t )_{\#}D)=\mathbf{M}_g(D)$ and 
                  $\left[\frac{\partial h}{\partial t}\right]_{t=0}=X$.    
            \end{Lemme}
            \begin{Dem}
\textbf{Construction of $h$.}
We start with the flow $\tilde{h}(t,x)$ of $X$ (i.e: $X(x):=\frac{\partial}{\partial t}\tilde{h}(t,x)_{|t=0}$) and we make a correction by a flow $H_s$ of a vector field $Y$ that has $Flux_g(Y)\neq 0$.
Now, we consider the function 
\Fonct{f}{I^2}{M}{(s,t)}{\mathbf{M}_g((H_s \circ h_t )(D))-Vol_g(D)}
where $I$ is an interval of the real line. It is smooth by classical theorems of differentiation of an integral, since we make an integration on rectifiable currents. We apply the implicit function theorem at point $(0,0)$ to function $f$ in order to find an $s(t)$ that satisfies $$\mathbf{M}_g((H_{s(t)} \circ \tilde{h}_t )(D))-Vol_g(D)=0.$$
Such an application of implicit function theorem is possible since
$$\frac{\partial}{\partial s}f(0,0)=Flux_g(Y)\neq 0.$$
We have also $s'(0)=0$. Indeed  
$$\frac{d}{dt}f(s(t),t)=s'(t)\int_{h_t(D)}div_g(Y) +\int_D div_g({H_{s(t)}}_* X)$$
and an evaluation at $t=0$ gives 
$$s'(0)Flux_g (Y)+Flux_g(X)=0$$
hence $s'(0)=0$ since $Flux_g(Y)\neq 0$ and $Flux_g(X)=0$.\\ 
Now if we apply the previous argument to $h(t,x)=H_{s(t)}\circ \tilde{h}(t,x)$ we can see that 
$$\frac{\partial}{\partial t}h(0,x)=s'(0)Y(h_t (x))+{H_{s(0)}}_* X =X,$$
by the fact $s'(0)=0$.
            \end{Dem}

\bigskip

\textbf{End of the proof of Lemma \ref{vpremiere}.}

Let $X$ be a vector field with $Flux_g(X)=0$. Applying Lemma \ref{flux}, there exists a variation
            $h(t,x)$ satisfing the following two properties
             \begin{enumerate}
                        \item $\mathbf{M}_g((h_t )_{\#}D)=\mathbf{M}_g(D)$ 
                        \item  $\frac{\partial h}{\partial t}_{t=0}=X$,
             \end{enumerate}
 provided $Flux_g(X)=0$ and
            $$\frac{d}{dt}\left[\mathbf{M}_g((h_t )_{\#}\partial D)\right]_{t=0}=\delta_g\partial D(X)=0.$$
            In other words, $Ker(Flux_g)\subseteq Ker(\delta_g\partial D)$.
            Hence there exists $\lambda\in\mathbb{R}$ for which it is true that $\delta_g\partial D=\lambda Flux$. We set $H_g=-\lambda$. This notation is justified by the fact that on the smooth part of $\partial D$, $H_g$ is equal to the genuine mean curvature.
      \end{Dem}

\subsection{Riemannian Monotonicity Formula using isometric embedding}

\begin{Thm}[Riemannian Monotonicity Formula]\label{Monriem}
           Let $T\in\mathbb{R}V_n(M)$ be a varifold solution of the isoperimetric problem, consider $x\in Spt||\partial T||$, and $R>0$. Then
           \begin{equation}
             \Theta (||i_{g\#}\partial T||,x)\omega_{n-1}R^{n-1}e^{-(|H_g|_g+\beta_{i,g})R}\leq ||i_{g\#}\partial T||B_{\R^n}(x,R),
           \end{equation}  
           where $H_{g}$ is the generalized mean curvature of the varifold $\partial T$ viewed as a varifold on $M$, $\beta_{i_g}$ is an upper bound on the norm of the second fundamental form of the isometric embedding $i_g:M\hookrightarrow\R^N$.
\end{Thm}
\begin{Dem} When $M$ is Euclidean space, this result is due to W. K. Allard, Theorem $5.1$ of \cite{All}. In order to adapt it to the situation considered here, we make use of an isometric embedding $i_g$ of $M$ (whose existence is guaranteed by Nash's theorem) and then we look at the current ${i_g}_{\#}T$ in order to apply the Euclidean statement. In this case we see that the term to consider, instead of simply taking into account the mean curvature of $T$ in $M$, involves the mean curvature of $i_{g\#}T$ into $\mathbb{R}^N$. This is not really a problem because of our control on the norm of the second fundamental form of the embedding of  $M$ in $\mathbb{R}^N$ by the upper bound $\beta_{i_g}$. Therefore 
$$
\Theta (||\partial T||,x)\omega_{n-1}R^{n-1}e^{-(|H_g|_g+\beta_i )R}\leq ||\partial T||B(x,R).
$$
\end{Dem}
\newpage
\section{The Normal Graph Theorem}\label{Sec:NormalGrphaTheorem}
\begin{Def}\label{Def:BoundedGeometry}
A complete Riemannian manifold $(M, g)$, is said to be of \textbf{bounded geometry}, if there exists a constant $k\in\mathbb{R}$, such that $Ric_M\geq k(n-1)$ (i.e., $Ric_M\geq (n-1)kg$ in the sense of quadratic forms) and $V(B_{(M,g)}(p,1))\geq v_0$ for some positive constant $v_0$, where $B_{(M,g)}(p,r)$ is the geodesic ball (or equivalently the metric ball) of $M$ centered at $p$ and of radius $r>0$.
\end{Def}
\begin{Thm}\label{tr2}
Let $(M^n,g)$ be a smooth Riemannian manifold endowed with a Riemannian metric $g$ of class $C^\infty$ with bounded geometry. Let $i_g:M\hookrightarrow\mathbb{R}^N$ be an isometric embedding. Let $B$ be an open relatively compact domain whose boundary $\partial B$ is smooth, 
$\alpha\in ]0,1[$, $\varepsilon>0$, given real numbers. 
           Then there exist\\ $\varepsilon_0=
           \varepsilon_0(n,B, \xi, g, \partial g, \partial^2g, \partial^3g, \partial^4g, \alpha, \varepsilon)>0$ and $C^*(1,\varepsilon, \varepsilon_0)>0$,\\ such that  for every current $T$ solution of the isoperimetric problem that satisfies the following condition 
\begin{equation}\label{Eq:Thm31Statement}
Vol_g(B\Delta T)\leq\varepsilon_0,
\end{equation}  
$\partial T$ is the normal graph of a function $u_T$ on $\partial B$, $u_T\in C^{1,\alpha}(\partial B)$, and $||u_T||_{C_g^{1,\alpha}(\partial B)}\le C^*(1, B, \partial^4g, \varepsilon, \varepsilon_0)$. Moreover  $C^*(1,\varepsilon, \varepsilon_0, \partial^4g)$ tends to $0$ as $\varepsilon,\varepsilon_0\to 0^+$ and the constant $\varepsilon_0$ is continuous with respect to its arguments and so in particular with respect to convergence of metrics in the $C^4$ topology. In particular, if $T_j$ is a sequence of isoperimetric regions such that $Vol_g(B\Delta T_j)\to 0$, then $||u_{T_j}||_{C_g^{1,\alpha}(\partial B)}\to 0$, $H^{\partial T_j}-H^{\partial B}\to 0$, $\partial B$ is a constant mean curvature hypersurface, and actually $B$ is an isoperimetric region. This convergence is uniform with respect to $g$. Furthermore for any positive integer $m\ge1$ and real number $\alpha\in]0,1[$ there exists a positive constant $C^*_{m+1}:=C^*(m,\varepsilon,\varepsilon_0, ||g||_{m, \alpha}, \partial^4g)>0$ such that 
$||u_T||_{C_g^{m+1,\alpha}(\partial B)}\le C^*_{m+1}$ where $C^*_{m+1}\to0$ as $\varepsilon\to0$, where $||g||_{m, \alpha}$ is the $C^{m,\alpha}$ norm of the metric tensor over a suitable compact neighborhood of $B$.
\end{Thm}
\begin{Rem} All the constants that bound the geometry of the ambient space are calculated on a tubular neighborhood of $\partial B$ contained in a compact $\mathcal{V}$ where the normal exponential map of $\partial B$ is a diffeomorphism, except for the confinement Lemmas \ref{Lemme:UniformlyBoundedDiameter}, \ref{tflat1}.
\end{Rem}
The proof of Theorem \ref{tr2} occupies paragraphs \ref{sketch} to \ref{preuve}.\\
We give at first an informal sketch of this proof and then a series of lemmas that are used in the rigorous proof.\\ 

\subsection{Sketch of the Proof of Theorem \ref{tr2}}
\label{sketch}
\begin{enumerate}
          \item Lemma \ref{Lemme:UniformlyBoundedDiameter} allows us to locate the entire picture of Theorem \ref{tr2} inside a compact tubular neighborhood of $B$. So all the quantities needed in the proof are bounded above and are bounded below away from $0$, hence the proof go in the same way as in the compact case. Furthermore we notice that Lemma \ref{Lemme:UniformlyBoundedDiameter} does not make any use of an isometric embedding of the Riemannian manifold $(M,g)$ into some Euclidean space. 
          \item We continue as in the compact case and we make use of an \textit{a priori} estimate of the mean curvature for isoperimetric regions, this is L\'evy-Gromov's lemma, stated in \ref{bornecm}. From the discussions contained in the proof of Lemma \ref{bornecm} we have that if the length of the mean curvature vector of $\partial T$ is strictly bigger than $\sqrt{k}$ then $T$ is always mean convex.
         \item Secondly, we apply Allard's regularity theorem (Riemannian, but still non intrinsic, version) to prove that $\partial T$ is a $C^{1,\alpha}$ submanifold and to prove $C^{1,\alpha}$ convergence at this point we make a crucial use of Nash's isometric embedding theorem.\\
To this aim we proceed as in the following steps: 
                  \begin{enumerate}
                           \item We stand on a sufficently small scale $R$ in order to estimate the first variation like required by Theorem \ref{allriem}. 
                           \item   We estimate  the volume of the intersection of $\partial T$ with a ball  $B_{M}(x,R)$ and we proceed as follows: we cut $\partial T$ with $B_{M}(x,R)$ and replace $T$ by $T'$ of equal volume thanks to the construction (Lemma \ref{lr2}) of a one parameter family of diffeomorphisms that perturbes $T$ preserving the volumes of perturbed domains. This leads to the estimates of Lemmas \ref{lr2}, \ref{lr5}.
                            \item We apply Allard's theorem and we conclude that $\partial T$ is of class $C^{1,\alpha}$. The tangent cone is hence a vector space. As showed by Frank Morgan in \cite{Morg1}, it follows that $\partial T$ is as smooth as the metric. We shall give a direct proof of this.    
                  \end{enumerate}
                   \item We confine $\partial T$ in a tubular neighborhood of $\partial B$, of sufficiently small thickness, in Lemma \ref{tflat1}. For this, \ref{lr2} is combined with the Riemannian monotonicity formula \ref{Monriem}. 
                   \item We calculate a bound on $r$ (the tubular neighborhood thickness) so that the projection  $\pi$, of the tubular neighborhood $\mathcal{U}_{r_0 }(\partial B)$ of thickness $r$ on $\partial B$, restricted to $\partial T$ is a local diffeomorphism and, after, via a topological argument we argue that $\pi_{|\partial T}$ is a global diffeomorphism. This shows that $\partial T$ is the global normal graph on $\partial B$ of a function $u$.   
By an application of the implicit function theorem, $u$ is then of class $C^{1,\alpha}$. Notice that $r=r(Vol_g(B\Delta T))\rightarrow 0$ when $Vol_g(B\Delta T)\rightarrow 0$.
                   \item The estimates presented in the conclusions of the Allard's regularity theorem shows that $||u||_{C^{1,\alpha}}\rightarrow 0$ when $Vol_g(B\Delta T)\rightarrow 0$. A geometric argument also shows that the $C^1$ norm of $u$ goes to zero if $r\rightarrow 0$, i.e., if $Vol_g(B\Delta T)\rightarrow 0$. Alternatively an appeal to Ascoli-Arzel\`a's theorem could be used to show that $||u||_{C^{1, \alpha}}\rightarrow 0$ when $r\rightarrow 0$.  
                \item Now we are ready to use elliptic regularity theory, Schauder's estimates, in order to find upper bounds on $||u||_{C^{2, \alpha}}$ and with the same technique of Ascoli-Arzel\`a of point $5$, we show $||u||_{C^{2, \alpha}}\rightarrow 0$ when $Vol_g(T\Delta B)\rightarrow 0$. In particular $H_{\partial T}\rightarrow H_{\partial B}$.
                \item Finally, when $B$ is the limit in flat norm of isoperimetric regions then by the continuity of the isoperimetric profile in bounded geometry and by lower semicontinuity of the perimeter we get that $B$ is isoperimetric, so with constant mean curvature. 
\end{enumerate}
\subsection{A priori estimates on mean curvature}
Set $$k:=Min\left\{-1,\inf_{\overline{\mathcal{U}_{r_0 }(\partial B)}}\mathcal{K}_1^{M}(x)\right\},$$ $$\delta:=Max\left\{\sup_{\overline{\mathcal{U}_{r_0 }(\partial B)}}\mathcal{K}_2^{M}(x),1\right\},$$
where $\mathcal{K}_1^{M}(x)$ is a lower bound on the sectional curvatures of $M$ at $x$, and $\mathcal{K}_2^{M}(x)$ is an upper bound on the sectional curvatures of $M$ at $x$.  
Denote by $H_{g,\partial T}$ the mean curvature vector of $\partial T$.  It is constant for isoperimetric domains. This means that the mean curvature vector have a constant scalar product with the fixed global defined inward pointing unit normal vector defined $\mathcal{H}^{n-1}$-a.e. on the support of the measure $||\partial T||$. The following Lemma is inspired by Theorem $2.2$ of \cite{MJ} in which only the case of $(M,g)$ compact is treated.
\begin{Lemme}\label{bornecm}
          Let $M^n$ be a complete not necessarily compact Riemannian manifold satisfying $Ricci\geq (n-1)k$, $k\in\R$. 
          Let $B$ an open bounded domain whose boundary $\partial B$ is smooth. Then there exists $\varepsilon_1 >0$ and $H_1 >0$ such that for every current $T$ solution of the isoperimetric problem that satisfies the condition $$Vol_g( T\Delta B)\leq\varepsilon_1 ,$$
          we have
\begin{equation}
          |H_g^{\partial T}|\leq H_{1,g},
\end{equation}
where $H_1=H_1 (n, k, Vol_g(B), Vol_g(M))=H_1(B, g)$, if $M$ is compact and $H_1=H_1(n,k, Vol_g(B))$, if  $M$ is noncompact. 
\end{Lemme}
\begin{Dem} We can assume in this proof without loss of generality that $\partial T$ is smooth. As we know from regularity theory, compare \cite{Morg1} or Theorem $2.28$ of \cite{Maggi}, if $T$ is an isoperimetric region every point of $Supp(\partial T)$ that have the tangent cone being an half space is a regular point, hence for any point $p\in M$ a minimizing geodesic issued from $p$ hits $\partial T$ orthogonally in a regular point, because we can put a tangent ball to $\partial T$ at half the distance between $p$ and $\partial T$. The tangency condition implies that the tangent cone is a half space. For this reason the arguments of this proof are not affected at all by the possible presence of singularities in $\partial T$. Compare on this issue $(4)$ of page $297$ of \cite{BessonGallotBerard} or the original paper \cite{Gr1}. Settled this first technical point we proceed with our proof. Set 
 \Fonct{c_{k}}{\mathbb{R}}{\mathbb{R}}{t}{\left\{\begin{array}{ll}
                                                                       cos(\sqrt{k}t), & if\: k>0,\\
                                                                       1, & if\: k=0,\\
                                                                       cosh(\sqrt{k}t), & if\: k<0,   
                                                                 \end{array}\right. }
      \Fonct{s_{k}}{\mathbb{R}}{\mathbb{R}}{t}{\left\{\begin{array}{ll}
                                                                       \frac{1}{\sqrt{k}}sin(\sqrt{k}t), & if\: k>0,\\
                                                                       t, & if\: k=0,\\
                                                                       \frac{1}{\sqrt{-k}}sinh(\sqrt{-k}t), & if\: k<0.    
                                                                 \end{array}\right. }
Let $h:=|H_{\partial T}|\in[0,+\infty[$ denote the length of the mean curvature vector of the regular part of the boundary $\partial T$ of an isoperimetric region $T$. It is well known that $h$ is a constant. The mean curvature vector of the regular part of $\partial T$ could point toward the interior or the exterior of the support of $T$. 
Now, fix $x\in\partial T$, denote by $\xi:=\xi(x)\in T_xM$ a unit vector normal to $\partial_rT$ at $x$, where $\partial_rT$ is the regular part of $\partial T$. Let us define $r_{x,\xi(x)}:=\sup\{t\in[0,+\infty[: d_{(M,g)}(\gamma_{\xi}(t), \partial T)=t\}$, where $\gamma_{\xi}$ is a geodesic parametrized by the arc length such that $\gamma(0)=x=\pi(\xi)$ and $\dot{\gamma_{\xi}}(0)=\xi$. Using Theorem $2.1$ of \cite{HeintzeKarcher} (see also \cite{BuragoZalgaller} Corollary $34.4.1$, or \cite{Chavel} Theorem $IX.3.2$) is not too hard to verify that
\begin{equation}\label{Eq:NormalCutLocus}
r_\xi\leq\tau,
\end{equation}
where $\tau$ is the first positive zero of $t\mapsto f_{n, k, h}(t):=c_k(t)-\frac h{n-1}s_k(t)$. Notice that when $k\le 0$ and $h>(n-1)\sqrt{-k}$, there exists a first positive zero $\tau\in]0,+\infty[$, otherwise when $h\le(n-1)\sqrt{-k}$ there is no first positive zero of $f_{n, k, h}$ and we set $\tau:=+\infty$. If $k>0$, then $\tau\in]0, \frac \pi{\sqrt{k}}[$. Assume for the moment that $h>(n-1)\sqrt{|k|}$, again by Theorem $2.1$ of \cite{HeintzeKarcher}
\begin{eqnarray}\label{Eq:MeanCurvatureBound}
v\le v+\tilde{v}_T\le v_{n,k,h,\partial T}\le A_g(\partial T)f(\tau),
\end{eqnarray}
where $\tilde{v}_T$ is the volume of a tubular neighborhood of $\partial T$ outside $T$,  
\begin{eqnarray}
v_{n,k,h,\partial T} & := & \int_{\partial T}\int_0^{r_{\xi}}\left(c_k(t)-\frac h{n-1}s_k(t)\right)^{n-1}\chi_{[0,\tau]}dt\\
& + &  \int_{\partial T}\int_0^{r_{-\xi}}\left(c_k(t)+\frac h{n-1}s_k(t)\right)^{n-1}dt,
\end{eqnarray} 
and for every $s\ge0$ we set
\begin{small}
$$f(s):=\int_0^s\left[\left(c_k(t)-\frac h{n-1}s_k(t)\right)^{n-1}\chi_{[0,\tau]}(t)+\left(c_k(t)+\frac h{n-1}s_k(t)\right)^{n-1}\right]dt.$$
\end{small}
As it is easy to check $f$ is a strictly increasing function, moreover we have that $f(s)\ge s$ for every $s$, hence by \eqref{Eq:MeanCurvatureBound} we get
\begin{equation}\label{Eq:MeanCurvatureBound0}
\frac v{A_g(\partial T)}\le f^{-1}\left(\frac v{A_g(\partial T)}\right)\le\tau.
\end{equation} 
From the last inequality it is easy to see that for every constant $0<c<1$, (say $c=\frac 12$) there exists $\varepsilon_1$ such that if $Vol_g(T\Delta B)\leq\varepsilon_1$, then by \eqref{Eq:MeanCurvatureBound0}
\begin{equation}\label{Eq:bornecm}
\tilde H_1:=(n-1)\cot_k\left(c\frac{Vol_g(B)}{I_{\mathbb{M}^n_k}(Vol_g(B))}\right)\geq h,
\end{equation}
since $(n-1)\cot_k(\tau)=h$, $\cot_k$ is a strictly decreasing function, $I_M\le I_{\mathbb{M}^n_k}$, $I_{\mathbb{M}^n_k}$ is a continuous function, and the perimeter is lower semicontinuous with respect to the convergence in flat norm. Thus we proved that $h\le\max\left\{\tilde H_1, (n-1)\sqrt{|k|}\right\}=H_1$ and the lemma follows. Now we can compare this proof with that of Theorem $2.2$ of \cite{MJ} in which the case when $M$ is compact is treated and a little better estimates are provided in that case. 
\end{Dem}

\subsection{Volume of the Intersection of a smooth hypersurface with a ball of the ambient Riemannian manifold} 
         Let $\tau_{\delta,\beta}>0$ be the first positive zero of the function $c_{\delta}-\beta s_{\delta}$.\\
         Set $\lambda(\beta, \delta)(t)=\frac{1}{c_{\delta}(t)-\beta s_{\delta}(t)}$ for $t\in [0,\tau_{\delta,\beta}[$.\\       
\begin{Lemme}\label{lr1}
          Let $M$ be a Riemannian manifold, $V\subset M$ be a smooth hypersurface.
There exists $R_2=R_2(V, g,\partial g,\partial^2 g)>0$ and $C_2(V, g,\partial g,\partial^2 g)>0$ such that for every $R<R_2$ and for every $x\in M$ at distance $d<R_2$ from $V$, if $R'=d+R$, then
\begin{eqnarray*}
Vol_g(V\cap B(x,R))\leq (1+C_2 R')\omega_{n-1}R'^{n-1}.
\end{eqnarray*}
$R_2$ depends only on $\beta$, $r_0$, $inj_{(M,g)}$ (bound on the second fundamental form of $V$, normal injectivity radius of $V$, injectivity radius of $M$ ), $\delta_0$ (geometry of the ambient Riemannian manifold) and $C_2$ depends on the same quantities plus a lower bound on Ricci curvature of $V$. 
      \end{Lemme}
\begin{Rem} 
In the proof of Theorem \ref{tr2} we apply Lemma \ref{lr1} with $V=\partial B$, $d\leq R^3$, but,  $d\leq R^2$ is enough too.
\end{Rem}
\begin{Rem} $\beta=\beta(V,g,\partial g, \partial^2g)$, $r_0=r_0(V,g,\partial g, \partial^2g)$, $inj_M=inj_{(M,g)}(V,g,\partial g, \partial^2g)$. 
\end{Rem}
\bigskip

\textbf{Idea of the proof}. Using comparison theorems  for distortion of the normal exponential map based on a submanifold, we can compare the intrinsic and extrinsic distance functions on $V\hookrightarrow M$. This allows us to reduce the problem to the estimation of the volume of an intrinsic ball of $V$, i.e., to Bishop-Gromov's inequality.\\

      \begin{Dem}
            Whenever $y\in V$ such that $d_{(M,g)}(x,V)=d_{(M,g)}(x,y)=d$ there exists $R''>0$ for which 
$$V\cap B(x,R)\subseteq B_{V}(y,R'').$$
            We can take for example $R''\geq\sup_{z\in V\cap B(x,R)} \{ d_{V}(y,z)\}$.\\ 
            Set $$k_2:=Min\{\frac{\inf\{Ric_{V}\}}{n-2},-1\}.$$ 
            then 
            \begin{eqnarray*}
   Vol_g(V \cap B(x,R))&\leq& Vol(B_{V}(y,R''))\\
&\leq& Vol_{\mathbb{M}^{n-1}_k }(B(o,R''))\\
&=&\alpha_{n-2}\int_0 ^{R''}s_{k_2} (t)^{n-2}dt,
            \end{eqnarray*}
            where the second inequality follows from Bishop-Gromov's Theorem. We have then 
            $$Vol_g(V \cap B_{M}(x,R))\leq (1+C'(k_2)(R'')^2 )\omega_{n-1}R''^{n-1}$$
            after expanding the term 
            $$\frac{V_{\mathbb{M}^{n-1}_k }(B(o,R''))-Vol_{\mathbb{R}^{n-1}}(B(o,R''))}{\omega_{n-1}R''^{n-1}}$$
            by a Taylor-Lagrange type formula.
            Let $\pi$ be the projection of $\mathcal{U}_{r_0}$ on $V$. Following a comparison result of 3.2.1 Main inequality and Corollary 3.3.1 of \cite{HeintzeKarcher} we get 
            \begin{equation}\label{Eq:HKCorollary}
                  (c_{\delta}(t)-\beta s_{\delta}(t))^2 g_0 \leq g_t \leq (c_{k}(t)+\beta s_{k}(t))^2 g_0  ,   
            \end{equation}
 where $g_t$ is the induced metric on the equidistant hypersurface $V_t:=\{x\in M: d_M(x, V)=t\}$ and the preceding expression is understood in the sense of quadratic forms.  
            Let $z\in V$ so that $d_{M}(x,z)=R$, $d_{V}(y,z)=R''$ and 
            $d_{M}(x,z)=b$. If we consider the minimizing geodesic $\gamma$ of $M$ that joins $y$ to $z$ parameterized by arc length $s$ and let us denote  
$\tilde{\Delta}=Sup_{s\in [0,b]}\{ d_{M}(\gamma (s),\partial B)\}$, there are points 
$p\in\partial B$, $q\in\gamma$, $p,q\in B_{M}(y,b)$ for which $\tilde{\Delta}=d_{M}(p,q)$ and conclude $\tilde{\Delta}\leq 2R$. 
If we take $R_2$ such that $0<R_2:=Min\{ \frac{\tau_{\delta,\beta}}{4}, inj_{M}\}$ this provides that $c_{\delta}-\beta s_{\delta}$ is decreasing and positive on $[0,R_2]$, we then infer 
\begin{eqnarray*}
R'' & \leq & l(\pi\circ\gamma)_{g_0}\\
& = & \int_0^b |d\pi(\gamma')|_{g_0}(s)ds\\
& \leq & \int_0^b \lambda(\beta, \delta)(s)|\gamma'|_{g_{d_M(\gamma(s), V)}}(s)ds\\
& = & \int_0^b \lambda(\beta, \delta)(s)ds\\
& \leq & \int_0^b \lambda(\beta, \delta)(2R)ds ,
\end{eqnarray*}
this last inequality leads certainly to
$$R''\leq\lambda(\beta, \delta)(2R)b.$$
But, $b\leq d+R$, by triangle inequality, hence
\begin{equation}
       R''(R)\leq\lambda(\beta, \delta)(b)b\leq (1+C(\beta,\delta)b)b.
\end{equation}
Incidentally we observe that the preceding equation gives us an analogue result to Lemma \ref{dgeod} in case of an arbitrary Riemannian ambient manifold, but still in codimension $1$. 
If we look at the Taylor expansion of $\lambda(\beta, \delta)(t)=1+\beta t+\mathcal{O}(t^2)$, we notice at a qualitative level that $$R''(R)\leq(1+\beta 2R+\mathcal{O}(R^2))(d+R)=(1+\mathcal{O}(R))(d+R)=(1+C R)(d+R),$$ where the constant $C=Sup_{R\in [0,R_2]}\{\frac{\lambda(\beta, \delta)(2R)}{R} \}$.
So we get
$$Vol(V \cap B_{M}(x,R))\leq (1+C'(k_2)((1+C R)(d+R))^2 )\omega_{n-1}((1+C R)(d+R))^{n-1} $$ 
and finally 
$$Vol_g(V \cap B_{M}(x,R))\leq (1+C_2 R')\omega_{n-1}R'^{n-1}$$
for $C_2$ depending on a lower bound on Ricci curvature tensor of  $V$, on an upper bound on the second fundamental form of $V$ and un upper bound on curvature tensor of ambient manifold.  
      \end{Dem}
\subsection{Compensation of Volume Process}\label{Compvol}
\textbf{Remark:} In this subsection we make no assumption on the distance of an arbitrary point $x$ of $\partial T$ to the boundary $\partial B$. 
Let $$R_3 :=Min\{inj_{(M,g)}, r_{0,g} , \frac{diam_g(B)}{4}\}=R_3(B, g,\partial g,\partial^2 g).$$
\begin{Lemme}[Deformation Lemma first version]\label{lr2}
            There exists $C_3=C_3(B, g,\partial g,\partial^2 g)>0$ such that whenever $R<R_3$, $a<\frac{R}{2}$, there is $\varepsilon_3>0$ so that, for every $x\in \partial T$, there exists a vector field $\xi_x$ with the following properties 
\begin{enumerate}
  \item the support of $\xi_x$ is disjoint from $B(x,R)$ ;
  \item the flow $\phi_t$ is defined for $t\in[-R,R]$, and for $t\in[-\frac{R}{2},\frac{R}{2}]$, $\xi_x$ restricted to a sufficiently small ball centered at a point $y'\in\partial B$, coincides with the gradient of the signed distance function to $\partial B$;
  \item the norm of the covariant derivative $|\nabla_g\xi_x |_g<C_3$.
\end{enumerate}
Furthermore, for every solution $T$ of the isoperimetric problem whose boundary contains $x$, and  $Vol_g(T\Delta B)<\varepsilon_3$, there exists $t\in[-a,a]$ such that $T'=(B\cap B(x,R))\cup (\phi_t (T)\setminus B(x,R))$ has volume equal to the volume of $T$.
In particular,
\begin{equation}\label{compaires}
\begin{array}{lll}
A_g(\partial T\cap B(x,R))& \leq & A_g(\partial B\cap B(x,R))+A_g((T\Delta B)\cap\partial B(x,R))\\
& + & A_g(\phi_{t\#}(\partial T))-A_g(\partial T).
\end{array}
\end{equation}
Constants $C_3$ and $\varepsilon_3$ depend only on the geometry of the problem, of the a priori choice of a vector field fixed once and for all on $\mathcal{U}_{\partial B}(r_0)$ and on a bump function $\psi$ defined once at all also.  
\end{Lemme}
\textbf{Remarks:}\begin{enumerate}
                                     \item In the proof of Theorem \ref{tr2} we use Lemma \ref{lr2} with $\varepsilon_0\leq\varepsilon_3$, among other contraints that will be clear in the sequel.
                                           \item Furthermore if  $\delta v:=Vol_g(B\cap B(x, R))-Vol_g(T\cap B(x,R))\leq 0$ then 
                                                  $t\geq 0$ and if $\delta v >0$ then $t<0$ (balancing of volume).
                                         \item The parameter $a$ serves to control  that $t$ be small, as
                                               this $t$ will control the term 
                                            $|Vol_g(T'\cap Supp(\varphi ))-Vol_g(T\cap Supp(\varphi ))|$   
                                \end{enumerate}
\bigskip

\textbf{Idea of proof}. The vector field $\xi_x$ is obtained with the classical technique of multiplication by a bump function the metric vector gradient of the signed distance function $\partial B$. This bump function
has support in a neighborhood of a point that belongs to $\partial B$ and that is far away from $x$. We provide also that the flow of this vector field significantly increases the volume of $B$. This is sufficient to suitably change the volume of $T$. We can then operate a balancing of a given volume variation.\\ \\
      \begin{Dem}
             First, we make the following geometric construction of a vector field $\nu$. 
             Fix a point $y'\in\partial B$ with $B(x,R)\cap B(y,R)=\emptyset$ (it suffices to take $y'$ so that 
              $d(x,y')\geq R+\frac{1}{2}diam(B)$, for example).\\
             Let $\mathcal{U}_{\partial B}(r_0):=\left\{ x\in M| d(x,\partial B)<r_0\right\}$. By the choice of $r_0$, the normal exponential map  
      \Fonct{exp^{\partial B}}{\partial B\times ]-r_0 ,r_0[}{\mathcal{U}_{\partial B}(r_0)}{(q,t)}{exp_q(t\nu(q))} 
      is a diffeomorphism.\\
       Let $\nu$ be the extension by parallel transport on normal (to $\partial B$) geodesics of the exterior normal issuing from $\partial B$ (equivalently, $\nu$ is the gradient of the signed distance function to $\partial B$), in a vector field defined on $\mathcal{U}_{r_0 }(\partial B)$.\\
       Let
           \Fonct{\psi}{\mathbb{R}}{[0,1]}{s}{\chi_{[0,1/2]}(|s|)+e^{4/3}e^{\frac{1}{s^2 -1}}\chi_{]1/2,1[}(|s|).}
Now, we modulate $\nu$ with the smooth function $\psi$ and we set $$\xi_x:=\psi(\frac{d(y',.)}{R})\nu=\psi_1 \nu .$$ 
It can be seen that $||\nabla_X \xi_x ||\leq ||\psi'||_{\infty,[-1,1]}||X||+||\nabla_X \nu||\leq C_3 ||X||$,
 establishing that $C_3$ depends on geometric quantities and on the choice of $\psi$. Let $\left\{\varphi_t\right\}$ be the flow (one parameter group of diffeomorphisms of $M$) of the vector field $\xi_x$. It's true that $Supp(\varphi)\subset B_{M}(y',R)$. Now consider, whenever $a\in ]0,\frac{R}{2}[$ the functions $f,f_1,h$ defined as follows:
                  \Fonct{f_1}{[-a,a]}{\mathbb{R}}{t}{Vol_{g,n}(\varphi_t (B))}
                  \Fonct{f}{[-a,a]}{\mathbb{R}}{t}{Vol_{g,n}(\varphi_t (\tilde{T})),} 
                  \Fonct{h}{[-a,a]}{\mathbb{R}}{t}{Vol_{g,n}(\varphi_t (T)),}
 where $\tilde{T}:=(T-B(x,R))\cup (B\cap B(x,R))$.\\
For the aims of the proof, we need to show that $Vol_g(T)\in f([-a,a])$ with an argument independent of $x$ as  $f$ depends on $x$.\\
 By construction 
              \begin{equation}
                  \begin{array}{ccc}
                  \frac{d}{dt}\left[Vol_g(\varphi_t (B))\right] & = & \int_{\varphi_t (\partial B)}\psi_1<\nu,\nu>dV_{\varphi_t (\partial B)} \\ 
& \geq & \psi(t)A_g(\partial B_t \cap Supp(\psi_1))\\
& = & A_g(\partial B_t \cap Supp(\psi_1)),  
                  \end{array}
              \end{equation}
hence letting $R':=\frac{R}{2(c_k +\beta s_k )(\frac{R}{2})}$ and \\
$(c_{\delta}-\beta s_{\delta})(\frac{R}{2})(Inf_{y'\in\partial B}V(\partial B\cap B(y', R')):=C'_3 ,$
\begin{equation}
\begin{array}{ccc}
f_1'(t) & \geq & Vol(\partial B_t \cap Supp(\psi_1))\\
& \geq & C'_3  ,
\end{array}
\end{equation} 
whenever $t<\frac{R}{2}$.\\
Hence $f_1$ is strictly increasing and $f_1(a)-f_1(-a)\geq 2aC'_3=:\Delta_3$.\\
                  Let 
                  $$J:=\left|det\left(\frac{\partial\varphi_t (y)}{\partial y}\right)\right|_{\infty,[-a, a]\times \overline{\mathcal{U}_{r_0 }(\partial B)}}\leq e^{nC_3 a},$$ 
                  by similar arguments to those of the proof of Lemma \ref{lr3}.\\ 
                  From
                  $$
                    \begin{array}{ccl}
                      |f(t)-h(t)| & = & |Vol_n (B\cap B(x,R))-Vol_n(T\cap B(x,R))|\\ 
                      & \leq & Vol((T\Delta B)\cap B(x,R))\\ 
                      & \leq & \varepsilon_3 ,\\
                      |h(t)-f_1(t)| & \leq & |Vol_g(\varphi_t (T\Delta B))|\\
                      & \leq & JVol(T\Delta B)\\
                      & \leq & e^{nC_3 a} \varepsilon_3 ,
                    \end{array}
                  $$
                  it follows that 
                  $$|f(t)-f_1(t)|\leq\varepsilon_3+J\varepsilon_3\leq (1+ e^{nC_3 a})\varepsilon_3=:\sigma ,$$
                  $\sigma$ is independent on $x$.\\
                  If we take 
                  \begin{equation}
                        0<\varepsilon_3\leq\frac{1}{2(1+ e^{nC_3 a})}aC'_3 ,
                  \end{equation}  
                  then  
                  \begin{equation}
\sigma\leq\frac{1}{2}\min\{f_1(0)-f_1(-a),f_1(a)-f_1(0)\},
\end{equation} 
                  therefore $$[f_1(-a)+\sigma, f_1(a)-\sigma]\subseteq f([-a,a]).$$
                  With this choice for $\varepsilon_3$ we obtain 
                  $$Vol_g(T)\in[f_1(-a)+\sigma, f_1(a)-\sigma],$$
                  so, there exists $t\in [-a,a]$ depending on $x$ such that 
 $f(t)=Vol_g(T)=Vol_g(\varphi_t (\tilde{T}))$ and we conclude the proof by taking $T':=\varphi_t (\tilde{T})$. \\
\newline
Finally
\begin{equation*}
          \begin{array}{ccc}
            A_g(\partial T) & = & I_{(M,g)}(Vol_g(T))\\
             & \leq & A_g(\partial T'), 
           \end{array}
      \end{equation*}
whence
           \begin{eqnarray}\label{Eq:CompensationFinal}\nonumber
            A_g(\partial T') & \le & A_g(\partial B\cap B(x,R))+Vol_g((T\Delta B)\cap\partial B(x,R))\\
           & + &A_g({\varphi_t }_{\#} (\partial T))-A_g(\partial T\cap B(x,R)),
           \end{eqnarray}
      which implies (\ref{compaires}).   
      \end{Dem}
\subsection{Comparison of the area of the boundary of an isoperimetric domain with the area of a perturbation with constant volume} 
\begin{Lemme}\label{lr3}
Let $M$ be a Riemannian manifold. For every $C>0$, for every vector field $\xi$ on $M$ such that $|\nabla_g\xi |_g<C$, whose flow is denoted by $\phi_t$, and whenever $V$ is a hypersurface embedded in $M$, it holds
\begin{eqnarray*}
Vol_g(\phi_{t\#}V)\leq e^{(n-1)C|t|}Vol_g(V).
\end{eqnarray*}
\end{Lemme}
      \begin{Dem}
It suffices to majorate the norm of the differential of diffeomorphism $\phi_t$. 
$$|d_x \phi_t (v)|=\left( g(\phi_t (x))(d_{x}\phi_t (v))\right) ^{\frac{1}{2}}=(\phi_t ^* (g_{M})(x)(v))^{\frac{1}{2}}$$
$$(\phi_t ^* (g_{M})(x)(v))^{\frac{1}{2}}\leq e^{C|t|}g(x)(v)=e^{C|t|}|v|.$$
The last inequality comes from the following lemma.
\begin{Lemme}\label{bourguignonl}
$(\phi_t ^* (g_{M})(x)(v))\leq e^{2C|t|}g(x)(v)$.
\end{Lemme}
\begin{Dem} From well known properties of Lie derivative we know that 
\begin{equation}
\frac{\partial}{\partial t}\left( \phi_t ^* (g_{M})\right) = \phi_t ^* \mathcal{L}_{\xi }g_{M}. 
\end{equation} 
We assume for the moment that we can show the following inequality 
\begin{equation}\label{bourguignon}
\mathcal{L}_{\xi }g_{M}= 2\times\textrm{symmetric part of } \nabla\xi .
\end{equation}
We use this fact to establish
 $$\mathcal{L}_{\xi }g_{M}\leq 2|\nabla\xi | g_{M}\leq 2Cg_{M},$$
hence $\phi_t^* \mathcal{L}_{\xi }g_{M}\leq 2C\phi_t ^* (g_{M})$. Set $t\mapsto  \phi_t ^* (g_{M})=q_t$, on $T_x M$, then it is not too hard to see that $q_t$ satisfies $\frac{\partial}{\partial t}q_t\leq 2Cq_t$ with $q_0=g_{M}$. 
It follows that whenever $x\in M$ and $v\in T_x M$, $q_t (v)\leq e^{2C|t|}q_0 (v)$ we have $(\phi_t ^* (g_{M})(x)(v))\leq e^{2C|t|}g(x)(v)$. It remains to show  that $\mathcal{L}_{\xi }g_{M}= 2 \times\textrm{symmetric part of }  \nabla\xi$. Let $A_{\xi}:=\mathcal{L}_{\xi }-\nabla_{\xi}$. We look at this operator on 2 covariant tensor fields and evaluate it on the metric $g_{M}$. We obtain $\mathcal{L}_{\xi }g_{M}=A_{\xi}g_{M}$, since $\nabla_\xi g=0$ and then  
$$0=A_{\xi}(g(w_1 , w_2 ))=(A_{\xi}g)(w_1 ,w_2)+g(-\nabla_{w_1}\xi , w_2 )+g(-w_1 ,\nabla_{w_2}\xi).$$
The first equality comes from the fact that the Lie derivative and the covariant derivative coincide when acting on functions, the others are straightforward consequences of the definition of $A_\xi$. So we conclude that $|\mathcal{L}_{\xi }g_{M}|\leq2|\nabla\xi|$.
\end{Dem}

\bigskip

\textbf{End of the proof of Lemma \ref{lr3}.}

We apply the inequality of Lemma \ref{bourguignonl} to the members of an orthonormal basis $(v_1 ,\ldots,v_{n-1})$ of the tangent space $T_x V$, we find
\begin{eqnarray*}
|\phi_{t\#}(v_1 \wedge\cdots\wedge v_{n-1})|_g\leq e^{(n-1)C|t|}.
\end{eqnarray*}
By an integration on $V$, one gets  
\begin{eqnarray*}
A_g(\phi_{t\#}V)\leq e^{(n-1)C|t|}A_g(V).
\end{eqnarray*}
      \end{Dem} 
\begin{Lemme}\label{lr6}
             Whenever $R>0$, $x\in Spt||\partial T||$ there exists $R_4$, $\frac{R}{2}<R_4< R$, such that 
             $$A_g((T\Delta B)\cap\partial B(x,R_4 ))\leq\frac{2}{R}Vol_g(T\Delta B).$$
\end{Lemme}
\begin{Dem} By a straightforward application of the coarea formula and the mean value theorem for integrals.
\end{Dem}  

\begin{Rem} At this point of the article we cannot put restrictions on the distance of $x\in\partial T$ to $\partial B$.
\end{Rem}
This lemma is used in the confinement lemma to majorate the volume of $\partial T$ in a geodesic ball.
In Lemma \ref{lr5}, we need to control the $(n-1)$-dimensional volume of the intersection of $\partial T$ with a geodesic ball of radius $R$ centered in $x$.
To make it possible we need to have the quantity $\frac{d_g(x,\partial B)}{R}$ very small.
\begin{Lemme}\label{lr5}
             Whenever $\eta>0$, there is $R_5$ such that whenever $R<R_5=R_5(B,\xi,g,\partial g, \partial^2g)$, (depending on the geometry of the problem) 
             there are $R_6$, $\varepsilon_6>0$ (depending only on $R$ and on the geometry of the problem, i.e., $B$, $\xi$, $g,\partial g, \partial^2g$) such that $0<\frac{R}{2}<R_6<R$ and if $T$ is a current solution of the isoperimetric problem with the property $Vol_g(B\Delta T)\leq\varepsilon_6$, then whenever $x\in Spt||\partial T||$ with $d_g(x,\partial B)\leq\left( \frac{R}{2}\right)^3$ we have
                   \begin{equation}
                         A_g(\partial T\cap B_{M}(x,R_6))\leq(1+\eta )\omega_{n-1}R_6^{n-1}.
                   \end{equation}
      \end{Lemme}
\begin{Rem}
                          In this context there are $2$ distance scales. The scale of $R_6$ the radius of the cutting geodesic ball of the ambient Riemannian manifold, that is the same as the scale of $R$ and that of $r_6$ that is the distance between an arbitrary point of $\partial T$ and a point of $\partial B$.  This is an important point in the estimates required by Allard's theorem, as the proof of Lemma \ref{lr1} shows.
Without this control on the scales involved we cannot have good control on the volume of the intersection of the hypersurface $\partial B$ with an ambient geodesic ball.
\end{Rem}
\begin{Rem}
The presence of interval $]\frac{R}{2},R[$ is just a technical complication due to the mean value theorem for integrals in the estimates of the $(n-1)$-dimensional volume of the part of $\partial T \cap B(x,R)$ that is $T\Delta B$. 
\end{Rem}      
      \begin{Dem} \\
                 Let $A:=C_2 s \left( 1+\frac{s^2}{2^3}\right)$, $B:=\left( 1+\frac{s^2}{2^3}\right) ^{n-1}-1$.\\
                 Let $R_5$ be the greatest positive real number $s$ such that  
\begin{enumerate}
           \item $s\leq Min\{inj_{M}, r_0 , \frac{diam(B)}{4}, R_3\}$,
           \item \begin{equation}\label{cmp2}
                              AB+B+A \leq \frac{1}{3}\eta.
\end{equation}
\end{enumerate}  
                 We fix $r_6 >0$ with $r_6\leq \left( \frac{R}{2}\right) ^3$.\\
                 Let $x\in Spt||\partial T||$.
Let $a$ be the greatest positive real number $s< \frac{R}{2}$ with 
 \begin{equation}\label{cmp1}
            (e^{(n-1)C_3s}-1)M\leq\frac{1}{3}\eta\omega_{n-1}\left( \frac{R}{2}\right) ^{n-1},  
 \end{equation} 
      where $M$ is the maximum of the isoperimetric profile on the interval $$[vol(B)/2,2vol(B)],$$ 
      i.e. 
      $$a\leq Min\{\frac{1}{(n-1)C_3}log\left[ 1+\frac{\eta\omega_{n-1}\left( \frac{R}{2}\right)^{n-1}}{3M}\right] ,\frac{R}{2}\}.$$
Set $\varepsilon_6 :=Min\{\varepsilon_3 , \frac{Vol(B)}{2}, \frac{1}{3}\eta\omega_{n-1}\left( \frac{R}{2}\right)^n \}$. 
Let $T$ be a solution of the isoperimetric problem such that $Vol(T\Delta B)<\varepsilon_6$.
By (\ref{compaires}) we find $t(x)\in[-a,a]$ and $\varepsilon_3$ (given by Lemma \ref{lr2}) 
satisfying 
 \begin{eqnarray}
A_g(\partial T\cap B(x,R) ) & \leq & A_g(\partial B\cap B(x,R) )\\ \nonumber
    & + & A_g((T\Delta B)\cap\partial B(x,R) )\\ \nonumber
& + & A_g({\varphi_t }_{\#} (\partial T))-A_g(\partial T). \nonumber
      \end{eqnarray}
From Lemmas \ref{lr2} and \ref{lr3} we have 
\begin{eqnarray}
            A_g(\partial T\cap B(x,R) ) & \leq & A_g(\partial B\cap B(x,R) )\\ \nonumber
             & + & A_g((T\Delta B)\cap\partial B(x,R) )\\ \nonumber
             & + & (e^{(n-1)C_3 t}-1)A_g(\partial T).\nonumber
      \end{eqnarray}  
By Lemma \ref{lr6} we get $R_4$ satisfying 
$$\begin{array}{ccc}
              A_g((T\Delta B)\cap\partial B(x,R_4) ) & \leq & \frac{2}{R}Vol(T\Delta B)\\
               & \leq & \frac{2}{R}\varepsilon_6.
     \end{array}$$ 
Let $R_6:=R_4$. Lemmas \ref{lr3}, \ref{lr6} and \ref{lr1} combined give
                \begin{eqnarray}\label{ep}
                           Vol_g(\partial T\cap B(x,R_6) )& \leq & (1+\mathcal{O}(R_6))\omega_{n-1}R_6 ^{n-1}\\
                           & + & \frac{2}{R}Vol_g(T\Delta B) +(e^{(n-1)C_3a}-1)M, \nonumber           
                \end{eqnarray} 
      as, by Lemma \ref{lr1},
$$Vol_g(\partial B\cap B(x,R) )\leq (1+\mathcal{O}(R))\omega_{n-1}R^{n-1},$$ and by Lemma \ref{lr6}, $0<\frac{R}{2}<R_6<R$.\\
       By  (\ref{cmp1}), (\ref{cmp2}), and the choice of $\varepsilon_6$, equation (\ref{ep}) becomes       
      \begin{eqnarray}
            Vol_g(\partial T\cap B(x,R_6) )& \leq & (1+\frac{1}{3}\eta)\omega_{n-1}R_6^{n-1}+\frac{1}{3}\eta\omega_{n-1}R_6 ^{n-1}\\
            & + & \frac{1}{3}\eta\omega_{n-1}R_6^{n-1}.\nonumber
      \end{eqnarray}
Finally
\begin{equation}
Vol_g(\partial T\cap B(x,R_6) )\leq (1+\eta)\omega_{n-1}R_6^{n-1}.  
\end{equation} 
      \end{Dem}
\subsection{Confinement of  an Isoperimetric Domain by Monotonicity Formula}
\begin{Lemme}\label{tflat1} Let $M^n$ be a Riemannian manifold. Let $B$ a compact domain whose boundary $\partial B$ is smooth. For every $s\in]0,R_3[$, there exists $\varepsilon_7(s)>0$ with the property that if $T$ is a current solution of the isoperimetric problem with
           $$Vol_g(B\Delta T)<\varepsilon_7 ,$$
           then $\partial T$ is contained in a tubular neighborhood of thickness $s$ of $\partial B$.
\end{Lemme}

\bigskip

           \textbf{Idea of the proof:}
By contradiction, we assume that there is a current $T$ and a point $x\in\partial T$ at distance $>s$ of $\partial B$. We choose $R\in ]s/2,s[$ so that the intersection $T\Delta B$ with the sphere $\partial B(x,R)$ has small area. The mechanism of balancing gives an estimation of the area of $\partial T\cap B(x,R)$, as $\partial B\cap B(x,R)=\emptyset$. This estimates from above contradicts the estimates from below given by monotonicity formula   (Lemma \ref{Monriem}), if $Vol_g(T\Delta B)$ is sufficiently small.\\                   
\begin{Dem} 
Set $s>0$. Let $H_1$ be the constant produced by Lemma \ref{bornecm}.
Let $C_3$ be the constant given by Lemma \ref{lr2}. Let $M_0$ be the maximum of the isoperimetric profile  on the interval $[Vol(B)/2,2Vol(B)]$. Let $\beta_i$ be a bound on the second fundamental form of an isometric immersion of $M$ in $\mathbb{R}^N$ the Euclidean space. 
We can choose $a$ so that 
 \begin{equation}
              (e^{(n-1)C_3a}-1)M_0<\frac{1}{2}\omega_{n-1}\left( \frac{s}{2}\right) ^{n-1}e^{-(H_1+\beta_i )s}.
      \end{equation}
Let $\varepsilon_3$ be the second constant given by Lemma \ref{lr2}, when, in this lemma, we take $R=s/2$. 
Let $\varepsilon_7 <\varepsilon_3$, $\varepsilon_7 <vol(B)/2$ and 
\begin{eqnarray*}
\frac{2\varepsilon_7}{s}<\frac{1}{2}\omega_{n-1}\left( \frac{s}{2}\right) ^{n-1} e^{-(H_1+\beta_i )s}.
\end{eqnarray*}
Let $T$ be a current solution of the isoperimetric problem satisfying
\begin{eqnarray*}
Vol_g(T\Delta B)<\varepsilon_7 .
\end{eqnarray*}
We argue by contradiction. Assume there is a point $x\in \partial T$ placed at distance $>s$ from $\partial B$. \\
The balancing of volume (Lemma \ref{lr2}) gives for all $R\leq Min\{s, R_3\}$
\begin{eqnarray*}
A_g(\partial T\cap B(x,R))\leq A_g((T\Delta B)\cap\partial B(x,R))+A_g(\phi_{t\#}(\partial T))-A_g(\partial T),
\end{eqnarray*}
as $B(x,R)\cap B=\emptyset$.
We apply Lemma \ref{lr3} with $C=C_3$ and we set $R_7 \in ]s/2, s[$ defining $R_7:=R_4$ obtained by applying Lemma \ref{lr6} with $R=s$ such that 
\begin{eqnarray*}
A_g((T\Delta B)\cap\partial B(x,R_7))\leq \frac{2}{s}Vol_g(T\Delta B).
\end{eqnarray*}
 It follows
\begin{eqnarray*}
A_g(\partial T\cap B(x,R_7))\leq \frac{2\varepsilon_7}{s}+(e^{(n-1)C_3a}-1)A_g(\partial T) 
\end{eqnarray*}
\begin{eqnarray*}
A_g(\partial T\cap B(x,R))\leq \frac{2\varepsilon_7}{s}+(e^{(n-1)C_3a}-1)M_0. 
\end{eqnarray*}
 Invoking Lemma \ref{bornecm} (L\'evy-Gromov), the mean curvature of $\partial T$ satisfies
\begin{eqnarray*}
|H|\leq H_1 .
\end{eqnarray*}
Monotonicity inequality (Lemma \ref{Monriem} ) gives us 
\begin{eqnarray*}
A_g(\partial T\cap B(x,R_7))\geq \omega_{n-1}R_7 ^{n-1}e^{-(|H|+\beta_i )R_7},
\end{eqnarray*}
which, by our choice of $\varepsilon_7$, contradicts the preceding inequality.
We conclude that $\partial T$ is contained in a tubular neighborhood of thickness $s$ of $\partial B$.
\end{Dem}

Loosely speaking the next theorem asserts that $\partial T$ is contained in a tubular neighborhood of thickness at most $C'_7Vol_g(T\Delta B)$ of $\partial B$, where $C'_7$ is a constant that depends only on $n$.
\begin{Lemme}\label{Lemma:ConfinementEstimates} Let $M^n$ be a Riemannian manifold. Let $B$ a compact domain whose boundary $\partial B$ is smooth. Then there exists a constant $C'_7=C'_7(n)>0$ such that if $s:=\sup\{x\in\partial T, d_g(x, \partial B)\}<\min\left\{R_3, \frac{\ln(2)}{H_1+\beta_{i_g}}\right\}$, then $s\le C'_7Vol_g(T\Delta B)^{\frac{1}{n}}$. 
\end{Lemme}
\begin{Dem} 
Let $H_1$ be the constant produced by Lemma \ref{bornecm}.
Let $C_3>0$ be the constant given by Lemma \ref{lr2}. Let $M_0$ be the maximum of the isoperimetric profile  on the interval $[Vol(B)/2,2Vol(B)]$. Let $\beta_i$ be a bound on the second fundamental form of an isometric immersion of $M$ in $\mathbb{R}^N$ the Euclidean space. 
We can choose $a$ so that 
 \begin{equation}
              (e^{(n-1)C_3a}-1)M_0<\frac{1}{2}\omega_{n-1}\left( \frac{s}{2}\right) ^{n-1}e^{-(H_1+\beta_i)s}.
      \end{equation}
Assume there is a point $x\in \partial T$ placed at distance $s$ from $\partial B$. \\
The balancing of volume (Lemma \ref{lr2}) gives for all $R<s$
\begin{eqnarray*}
A_g(\partial T\cap B(x,R))\leq A_g((T\Delta B)\cap\partial B(x,R))+A_g(\phi_{t\#}(\partial T))-A_g(\partial T),
\end{eqnarray*}
as $B(x,R)\cap B=\emptyset$.
We apply Lemma \ref{lr3} with $C=C_3$ and we set $R_7 \in ]s/2,s[$ defining $R_7:=R_4$ obtained by applying Lemma \ref{lr6} with $R=s$ such that 
\begin{eqnarray*}
A_g((T\Delta B)\cap\partial B(x,R_7))\leq \frac{2}{s}Vol_g(T\Delta B).
\end{eqnarray*}
 It follows
\begin{eqnarray*}
A_g(\partial T\cap B(x,R_7))\leq \frac{2Vol_g(T\Delta B)}{s}+(e^{(n-1)C_3a}-1)A_g(\partial T), 
\end{eqnarray*}
hence
\begin{eqnarray*}
A_g(\partial T\cap B(x,R))\leq \frac{2Vol_g(T\Delta B)}{s}+(e^{(n-1)C_3a}-1)M_0. 
\end{eqnarray*}
 Invoking Lemma \ref{bornecm} (L\'evy-Gromov), the mean curvature of $\partial T$ satisfies
\begin{eqnarray*}
|H|\leq H_1 .
\end{eqnarray*}
Monotonicity inequality (Lemma \ref{Monriem}) gives us 
\begin{eqnarray*}
A_g(\partial T\cap B(x,R_7))\geq \omega_{n-1}R_7 ^{n-1}e^{-(|H|+\beta_i )R_7},
\end{eqnarray*}
thus 
\begin{eqnarray*}
\omega_{n-1}R_7 ^{n-1}e^{-(|H|+\beta_i )R_7}\le\frac{2Vol_g(T\Delta B)}{s}+(e^{(n-1)C_3a}-1)M_0, 
\end{eqnarray*}
which in turn gives
\begin{eqnarray}
\omega_{n-1}\left(\frac{s^n}{2^{n+1}}\right) & \le & s\omega_{n-1}R_7 ^{n-1}e^{-(H_1+\beta_{i_g})R_7}\\
& \le & Vol_g(T\Delta B). 
\end{eqnarray}
Setting $C'_7:=\frac{2^{\frac{n+1}{n}}}{\omega_{n-1}^{\frac{1}{n}}}$
We conclude that $\partial T$ is contained in a tubular neighborhood of thickness $s$ of $\partial B$.
\end{Dem}
\subsection{Alternative proof of confinement under weaker bounded geometry assumptions} 
We present here an alternative proof of the results contained in the preceding section under weaker assumptions on the way the geometry of $(M,g)$ is bounded. The main result of this section is Lemma \ref{Lemme:UniformlyBoundedDiameter}. Before stating and proving it, we need an important technical deformation lemma in the spirit of what is called today Almgren's Lemma. Instances of this kind of lemma are Lemma \ref{lr2}, Lemma $4.8$ of \cite{NardulliOsorio}, Lemma $17.21$ of \cite{Maggi} and Lemma $4.5$ of \cite{GalliRitore}, but in the literature there are plenty of ad-hoc versions of it . Roughly speaking we deform an isoperimetric region $\Omega$ by a small amount of volume $\Delta v$ controlling the amount of variation of area $\Delta A$ by a constant $C$ times $\Delta v$, i.e., $\Delta A\le C\Delta v$. In general the constant $C$ depends on $\Omega$, but in our specific situation we need to have an uniform constant $C>0$ independent of $\Omega$ if $\Omega$ is close enough in flat norm to $B$. To overcome this difficulty we prove the following uniform deformation lemma which needs the notion of normal injectivity radius of an arbitrary codimension submanifold, which in turn generalizes the notion of injectivity radius at a point.
\begin{Def} Let $(M^n,g)$ be a Riemannian manifold, $0\le m\le n$, and $N^m\subseteq M$ be a $m$-dimensional submanifold of $M$. Consider $T^1N:=\{(p,w)\in TM: w\in T_pN^\perp, ||w||_g=1\}$ the \textbf{unit tangent bundle of $N^k$}. For any $(p,w)\in T^1N$ let us define the nonnegative extended real numbers $r_{0,g, N}(p,w):=\sup\{t>0:d_g(exp_p(tw), N)=t\}\in ]0,+\infty]$ and $r_{0,g, N}:=\inf\{r_{0,g}(p,w):(p,w)\in T^1N\}\in ]0,+\infty]$. We call $r_{0,g, N}$ the \textbf{normal injectivity radius of $N$}. 
\end{Def}
\begin{Rem}\label{Rem:NormalInjectivityRadius} Notice that in the language of the Definition at the end of page $145$ of $\cite{Gra}$ we have $r_{0,g, N}=minfoc(\partial N)$.
\end{Rem}
\begin{Rem} By the choice of $r_{0,g, N}$ and standard comparison results for the shape operator, see for instance Equation $(7.23)$ and Lemma $8.51$ of \cite{Gra} we know that 
\begin{equation*}
r_{0,g, N}\ge\cot_{\Lambda_{g, N}}^{-1}(\max\{1,\beta_g\})=r_0=r_0(N, g, \partial g, \partial ^2 g)>0, 
\end{equation*}
where $\Lambda_{g, N}:=\sup\{K_{(M,g)}(x):x\in d_g(x, N)\le1\}$ with $K_{(M,g)}(x)$ being the maximum taken over all the sectional curvature of $2$-plane in $T_xM$ with respect to the Riemannian metric $g$ and $\beta_g$ is an upper bound on the second fundamental form of the isometric embedding of $(N,g_{|_{N}})$ into $(M,g)$.
\end{Rem}
To simplify the notation in what follows we set $r_{0,g}:=r_{0,g, \partial B}$. In first we make the following geometric construction of a vector field $\nu$. Fix a point $y'\in\partial B$. Let $\mathcal{U}_{\partial B}(r_{0,g}):=\left\{ x\in M| d_g(x,\partial B)<r_{0,g}\right\}$. It is well known that the normal exponential map  
      \Fonct{exp_g^{\partial B}}{\partial B\times ]-r_{0,g} ,r_{0,g}[}{\mathcal{U}_{\partial B}(r_{0,g})}{(q,t)}{exp_q(t\nu(q))} 
      is a diffeomorphism. Let $\nu$ be the extension by parallel transport on normal (to $\partial B$) geodesics of the exterior normal issuing from $\partial B$ (equivalently, $\nu$ is the gradient of the signed distance function to $\partial B$), in a vector field defined on $\mathcal{U}_{r_{0,g}}(\partial B)$.\\
       Let
           \Fonct{\psi}{\mathbb{R}}{[0,1]}{s}{\chi_{[0,1/2]}(|s|)+e^{4/3}e^{\frac{1}{s^2 -1}}\chi_{]1/2,1[}(|s|),} by a direct computation it is easy to check that $||\psi'||_{\infty,[-1,1]}\le 4$. Now, we modulate $\nu$ with the smooth function $\psi$ and we set 
           \begin{equation}\label{Eq:xiDef}
                     \xi:=\psi(\frac{d_g(y',.)}{r_{0,g}})\nu=\psi_1 \nu.
           \end{equation}          
\begin{Lemme}[Uniform Deformation Lemma second version]\label{Lemma:DeformationLemma2} Let $(M,g)$ be a complete Riemannian manifold $($without any further assumption on $g$$)$, $B\subseteq M$ be an open relatively compact set with smooth boundary, $y'\in\partial B$, $r_{0,g}$ the normal injectivity radius of $\partial B$, and $\xi$ the smooth vector field with $Supp(\xi)\subseteq B_g(y', r_{0,g})$ defined by \eqref{Eq:xiDef}. Then there exist $\varepsilon_8=\varepsilon_8(B, \partial B, \xi, g, \partial g,\partial ^2 g)>0$, $C_8=C_8(B,\partial B, \xi, g, \partial g,\partial ^2 g)>0$, and $\sigma_0=\sigma_0(B, \partial B, \xi, g, \partial g,\partial ^2 g)>0$ such that for every finite perimeter set $\Omega$ with $Vol_g(\Omega\Delta B)\le\varepsilon_8$ and $\sigma\in[-\sigma_0, \sigma_0]$ there exist $T'$ a finite perimeter set $($or $n$-rectifiable current$)$ such that $Vol_g(T')=Vol_g(\Omega)+\sigma$, $T'\Delta\Omega\subseteq B_g(y', r_0)$, and 
\begin{equation}
A_g(\partial T')\le A_g(\partial\Omega)+C_8A_g(\partial\Omega\cap Supp(\xi))|V_g(T')-V_g(\Omega)|.
\end{equation}
\end{Lemme}
\begin{Dem}  
It can be seen easily that $$||\nabla_X \xi||\leq\frac{1}{r_{0,g}}||\psi'||_{\infty,[-1,1]}||X||+||\nabla_X \nu||\leq C_3 ||X||,$$ where $$C_3=C_3(\tilde{\beta}_g, r_{0,g})=\frac{4}{r_{0,g}}+\sup\{||II_g||_{\infty,\partial B_t}:t\in[-r_{0,g}, r_{0,g}]\}=\frac{4}{r_{0,g}}+\tilde{\beta}_g,$$ being $\partial B_t:=\tilde d_g^{-1}(t)$ the level set of the signed distance function $\tilde d$ to $\partial B$. This last equation establishes readily that $C_3=C_3(B, \partial B, \xi, g, \partial g, \partial^2 g)>0$ depends on geometric quantities and on the choice of $\psi$. Let $\left\{\varphi_t\right\}$ be the flow (one parameter group of diffeomorphisms of $M$) of the vector field $\xi$. It is immediate to check that $Supp(\varphi)\subset B_{(M,g)}(y', r_{0,g})$. Now, consider, whenever $a\in ]0,\frac{r_{0,g}}{2}[$ for example $a:=\frac{r_{0,g}}{4}$ the functions $f_1,h$ defined as follows:
                  \Fonct{f_1}{[-a,a]}{[0,+\infty[}{t}{Vol_{g,n}(\varphi_t (B)),}
                  \Fonct{h}{[-a,a]}{[0,+\infty[}{t}{Vol_{g,n}(\varphi_t (\Omega)).}
For the aims of the proof, we need to show that $Vol_g(\Omega)+\sigma\in h([-a,a])$ for sufficiently small $\sigma$. First of all assume that 
\begin{equation}\label{Eq:DeformationLemmaVolumeConstraint}
\varepsilon_8<\min\{\frac{1}{2}Vol_g(B_g(y', r_{0,g})\setminus B), \frac{1}{2}Vol_g(B_g(y', r_{0,g})\cap B)\},
\end{equation} 
to have enough space to put enough volume inside $B_g(y', r_{0,g})$ and to have $A_g(\partial\Omega\cap B_g(y', r_0))>0$. By the first variation formula for volumes in maximal dimension $n$ we have  
              \begin{equation}
                  \begin{array}{ccc}
                  f_1'(t)=\frac{d}{dt}\left[Vol_g(\varphi_t (B))\right] & = & \int_{\varphi_t (\partial B)}\psi_1<\nu,\nu>dVol_{g,\varphi_t (\partial B)} \\ 
& \geq & A_g(\partial\varphi_t (B)\cap B_{(M,g)}(y',\frac{r_{0,g}}{2})), 
                  \end{array}
              \end{equation}
hence letting $R':=\frac{r_{0,g}}{2(c_k +\beta s_k )(\frac{r_{0,g}}{2})}$ and \\
$$(c_{\delta}-\beta s_{\delta})(\frac{r_{0,g}}{2})(Inf_{z\in\partial B}A_g(\partial B\cap B_g(z, R'))=:C'_3=C'_3(\partial B, g, \partial g, \partial ^2 g)>0,$$ thus
\begin{equation}\label{Eq:LemmaConfinementVolumeLowerboundB}
\begin{array}{ccc}
f_1'(t) & \geq & A_g(\partial\varphi_t (B)\cap B_{(M,g)}(y',\frac{r_{0,g}}{2}))\\
& \geq & C'_3>0,
\end{array}
\end{equation} 
whenever $t<\frac{r_{0,g}}{2}$. Hence $f_1$ is strictly increasing on $[-a, a]$ and $$f_1(a)-f_1(-a)\geq 2aC'_3=:\Delta_3>0.$$
                  Let us define $J$ as 
        $$J:=\left|det\left(\frac{\partial\varphi_t (y)}{\partial y}\right)\right|_{\infty,[-a, a]\times \overline{\mathcal{U}_{r_{0,g}}(\partial B)}},$$ 
                  by similar arguments to those of the proof of Lemma \ref{lr3} we obtain
                  $$J\leq e^{nC_3 a}.$$ 
                  From
                  $$
                    \begin{array}{ccl}
                      |h(t)-f_1(t)| & \leq & |Vol_g(\varphi_t (\Omega\Delta B))|\\
                      & \leq & JVol_g(\Omega\Delta B)\\
                      & \leq & e^{nC_3 a} \varepsilon_8,
                    \end{array}
                  $$
                  it follows that 
                 \begin{equation}
                 |h(t)-f_1(t)|\leq e^{nC_3 a}\varepsilon_8=:\delta.
                 \end{equation}               
                  Now we want to estimate $h'(t)$ from below. The idea behind this estimates is that when $\Omega$ is close in flat norm to $B$ the flux of $\xi$ through $\Omega$ is close to the flux of $\xi$ trough $B$. Formally we have  
\begin{tiny}
\begin{eqnarray}
|h'(t)-f'_1(t)| & = & \left|\int_{\partial\varphi_t(\Omega)}\langle\xi,\nu_{\varphi_t(\Omega)}\rangle_gd\mathcal{H}^{n-1}_g-\int_{\partial\varphi_t(B)}\langle\xi,\nu_{\varphi_t(B)}\rangle_gd\mathcal{H}^{n-1}_g\right|\\
 & = & \left|\int_{\varphi_t(\Omega)}div_g(\xi)dVol_g-\int_{\varphi_t(B)}div_g(\xi)dVol_g\right|\\
 & = & \left|\int_{\varphi_t(\Omega)\Delta\varphi_t(B)}div_g(\xi)dVol_g\right|\\
 & \le & ||div_g(\xi)||_{\infty, \partial B\times[r_{0,g}, r_{0,g}]}Vol_{g}(\varphi_t(\Omega)\Delta\varphi_t(B))\\
 & \le & Ce^{nC_3t}Vol_g(\Omega\Delta B)\\
 & \le & Ce^{nC_3r_{0,g}}Vol_g(\Omega\Delta B), 
\end{eqnarray}
\end{tiny}
where $C:=||div_g(\xi)||_{\infty, \partial B\times[r_{0,g}, r_{0,g}]}=C(\partial B, g, \partial g, \partial ^2 g)>0$.
Hence choosing 
\begin{equation}\label{Eq:DeformationLemmaVolumeConstraint0}
\varepsilon_8\le\frac{C'_3}{2Ce^{nC_3r_{0,g}}},
\end{equation} 
we get
\begin{equation}
h'(t)\ge f'_1(t)-Ce^{nC_3r_{0,g}}Vol_g(\Omega\Delta B)\ge\frac{C'_3}{2}>0,\;\; \forall t\in[0,r_{0,g}].
\end{equation}
Integrating over the interval $[0,t]$ this last inequality we easily conclude
\begin{equation}
\Delta v:=|Vol_g(\varphi_t(\Omega))-Vol_g(\Omega)|=h(t)-h(0)\ge \frac{C'_3}{2}t,\;\; \forall t\in[0,r_{0,g}].
\end{equation}
Combining this last equation with Lemma \ref{lr3} and putting $T':=\varphi_t(\Omega)$ leads to
\begin{eqnarray}
A_g(\partial T') & \le  & A_g(\partial\Omega\cap (M\setminus B_g(y', r_0)))\\
& + & e^{(n-1)C_3t}A_g(\partial\Omega\cap B_g(y', r_0))\\
& \le & A_g(\partial \Omega)\\
& + & \frac{e^{(n-1)C_3r_{0,g}}-1}{r_{0,g}}tA_g(\partial\Omega\cap B_g(y', r_0))\\ 
 & \le & e^{(n-1)C_3t}A_g(\partial\Omega)\\
& \le & \left(1+\frac{e^{(n-1)C_3r_{0,g}}-1}{r_{0,g}}t\right)A_g(\partial \Omega)\\
& \le & A_g(\partial \Omega)+\frac{2\Delta v(e^{(n-1)C_3r_{0,g}}-1)}{C'_3r_{0,g}}A_g(\partial \Omega)\\
& \le & A_g(\partial \Omega)+\frac{2\Delta v(e^{(n-1)C_3r_{0,g}}-1)}{C'_3r_{0,g}}A_g(\partial \Omega).
\end{eqnarray}

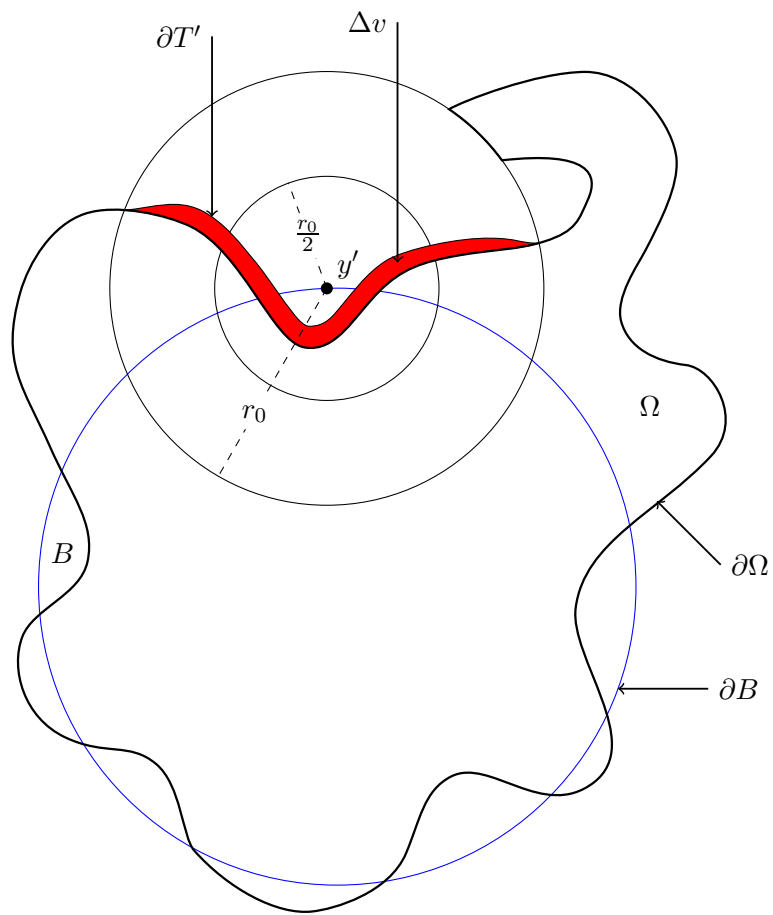
\begin{figure}
\label{fig:deformation}
\centering
\begin{tikzpicture}[scale=0.04]
\def\a{99.296};         
\def\b{72.141};   	
\def\c{37.279};   	
     
\coordinate (A) at (355.445,374.270);  
\coordinate (B) at ($(355.445,374.270)+(92:99.296)$); 
\coordinate (C) at (416.347,488.410);
\coordinate (D) at (415.347,490.200);
\coordinate (E) at (406.347,490.200);
\coordinate (F) at (396.347,490.500);

\coordinate (G) at (309.995,499.548);
\coordinate (H) at (301.894,503.334);
\coordinate (I) at (294.995,500.548);

\def\ballA{  (A) circle (\a)};
\def\ballB{  (B) circle (\b)};
\def\ballC{  (B) circle (\c)};

\def\curveA{
(272.388,382.608)
.. controls (275.034,395.275) and (267.069,404.488) .. (260.022,420.475)
.. controls (251.147,440.61) and (243.241,448.3) .. (250.026,472.092)
.. controls (252.882,482.106) and (260.474,494.532) .. (271.409,498.219)
.. controls (281.947,501.771) and (298.6,497.878) .. (306.605,494.026) 
.. controls (323.8,485.751) and (333.89,455.731) .. (344.367,453.825) 
.. controls (355.025,451.886) and (360.755,463.892) .. (369.57,472.647)
.. controls (376.634,479.663) and (383.165,482.244) .. (396.005,484.44) 
.. controls (408.424,486.564) and (415.036,486.175) .. (424.196,489.035)
.. controls (434.958,492.395) and (435.777,496.472) .. (438.698,502.956) 
.. controls (442.532,511.467) and (439.172,519.283) .. 
([shift={(36.2:\b)}]B) arc (36.2:55.7:\b)
.. controls (396.869,535.644) and (421.523,545.319) .. (437.689,545.443)
.. controls (454.131,545.569) and (467.535,527.299) .. (468.183,516.073)
.. controls (468.684,507.394) and (463.252,496.469) .. (457.677,485.801)
.. controls (455.375,481.397) and (447.606,470.316) .. (449.756,462.439)
.. controls (452.172,453.59) and (460.99,449.129) .. (471.261,447.927)
.. controls (480.49,446.847) and (490.791,429.719) .. (479.654,417.693)
.. controls (462.375,399.035) and (437.348,389.337) .. (434.622,367.001)
.. controls (432.611,350.523) and (457.816,322.315) .. (440.696,309.083)
.. controls (423.855,296.066) and (406.112,317.646) .. (393.221,311.042)
.. controls (384.065,306.352) and (382.926,296.564) .. (378.241,285.487) 
.. controls (372.664,272.299) and (354.545,266.937) .. (346.973,266.131)
.. controls (335.304,264.889) and (319.032,273.742) .. (307.809,286.436) 
.. controls (303.902,290.855) and (303.553,308.17) .. (295.182,315.2)
.. controls (286.195,322.746) and (278.607,317.637) .. (264.971,324.056)
.. controls (251.191,330.543) and (246.599,346.739) .. (250.699,357.528)
.. controls (254.318,367.051) and (270.034,371.341) .. (272.388,382.608)
--cycle;
}
\def\curveB{
([shift={(158.88:\b)}]B) .. controls ([shift={(158.88:\b)}]B) and 
(296.6,498.878) ..(306.605,494.026)
.. controls (323.8,485.751) and (333.89,455.731) .. (344.367,453.825) 
.. controls (355.025,451.886) and (360.755,463.892) .. (369.57,472.647)
.. controls (376.634,479.663) and (383.165,482.244) .. (396.005,484.44)
.. controls (408.424,486.564) and (415.036,486.175) .. 
([shift={(12.05:\b)}]B)
.. controls (C) and (D) .. (E)
.. controls (F) and (377.852,486.589) .. (371.44,482.537)
.. controls (364.946,478.433) and (359.825,471.219) .. (356.971,467.872) 
.. controls (354.162,464.577) and (351.681,460.625) .. (345.893,460.839) 
.. controls (343.133,460.941) and (338.896,467.572) .. (334.642,473.476)
.. controls (333.172,475.516) and (330.545,479.477) .. (327.484,483.093)
.. controls (323.894,487.334) and (317.995,495.548) .. (G) 
.. controls (H) and (I).. 
([shift={(158.88:\b)}]B) --cycle}
\draw[blue] \ballA;
\draw[] \ballB;
\draw[] \ballC;
\fill[red] \curveB;
\draw[] \curveB;
\draw[] [line width=0.3mm] \curveA;
\draw ([shift={(36.2:\b)}]B) arc (36.2:55.7:\b);
\fill (B) circle (2);
\node [above right] at  (B) {$y'$};
\node  at  ($(A)+(173:92)$) {$B$};
\node at  ($(A)+(30:120)$) {$\Omega$} ;
\draw[dashed] (B)--($(B)+(240:\b)$);
\fill [white] ($(B)+(240:{2*\b/3})$) circle (6);
\node at  ($(B)+(240:{2*\b/3})$) {$r_0$} ; 
\draw[dashed] (B)--($(B)+(110:\c)$);
\fill [white] ($(B)+(110:\c/2)$) circle (8);
\node at  ($(B)+(110:\c/2)$) {$\frac{r_0}2$} ; 
\draw[line width=0.7pt, <-]   ($(A)+(15:110)$)-- ($(A)+(15:110)+(-45:30)$)  node[right] {$\partial\Omega$} ;
\draw[line width=0.7pt, <-]   ($(A)+(-20:\a)$)-- ($(A)+(-20:\a)+(0:30)$)  node[right] {$\partial B$} ;
\draw[line width=0.7pt, <-]   ($(B)+(20:25)$)-- ($(B)+(20:25)+(90:80)$)  node[left] {$\Delta v$} ;
\draw[line width=0.7pt, <-]   ($(B)+(148:45)$)-- ($(B)+(148:45)+(90:60)$)  node[left] {$\partial T'$} ;
\end{tikzpicture}
\caption{Illustration of the Uniform Deformation Lemma.} 
\end{figure}
Hence we can choose $C_8:=\frac{2(e^{(n-1)C_3r_{0,g}}-1)}{C'_3r_{0,g}}=C_8(\partial B, B, \xi, g, \partial g,\partial ^2 g)>0$.                  
                                    Integrating \eqref{Eq:LemmaConfinementVolumeLowerboundB} over the interval $[0,a]$ we get $f_1(a)-f_1(0)>aC'_3>0$, and again integrating over the interval $[-a,0]$ we get $f_1(0)-f_1(-a)>aC'_3>0$ so if we choose 
                  \begin{equation}\label{Eq:DeformationLemmaVolumeConstraint1}
                        0<\varepsilon_8\leq\frac{1}{2e^{nC_3 a}}aC'_3 ,
                  \end{equation}  
                  then  
                  \begin{equation}
\delta\leq\frac{1}{2}\min\{f_1(0)-f_1(-a),f_1(a)-f_1(0)\},
\end{equation} 
                  therefore $$[f_1(-a)+\delta, f_1(a)-\delta]\subseteq h([-a,a]).$$
                  So for every $\sigma\in [-\tilde\sigma_0, \tilde\sigma_0]$ where $$\tilde\sigma_0:=\min\{|f_1(-a)+\delta-Vol_g(\Omega)|, |f_1(a)-\delta-Vol_g(\Omega)|\},$$ there exists $t\in [-a,a]$ such that $Vol_g(\varphi_t(\Omega))=Vol_g(T')=Vol_g(\Omega)+\sigma$. Taking $\varepsilon_8$ possibly smaller, i.e., 
                  \begin{equation}\label{Eq:DeformationLemmaVolumeConstraint2}
                  0<\varepsilon_8\le\frac{1}{2(1+e^{-nC_3 a})}\min\{-f_1(-a)+Vol_g(B),f_1(a)-Vol_g(B)\},
                  \end{equation} we have $$\tilde\sigma_0\ge\frac{1}{2}\min\{|f_1(-a)-Vol_g(B)|, |f_1(a)-Vol_g(B)|\}>0.$$ Thus we can choose $$\sigma_0:=\frac{1}{2}\min\{|f_1(-a)-Vol_g(B)|, |f_1(a)-Vol_g(B)|\}=\sigma_0(\partial B, B, \xi, g, \partial g,\partial ^2 g)>0.$$ Since   
                  $\varepsilon_8$ satisfies \eqref{Eq:DeformationLemmaVolumeConstraint}, \eqref{Eq:DeformationLemmaVolumeConstraint0}, \eqref{Eq:DeformationLemmaVolumeConstraint1}, and \eqref{Eq:DeformationLemmaVolumeConstraint2} we argue that we can choose $\varepsilon_8=\varepsilon_8(B, \partial B, \xi, g, \partial g,\partial ^2 g)>0$ and this finishes the proof. 
  \end{Dem}
 
We now state our desired confinement lemma in bounded geometry.
\begin{Lemme}[Confinement Lemma General Case]\label{Lemme:UniformlyBoundedDiameter} Let $(M^n,g)$ be a complete Riemannian manifold, with bounded geometry. Let $B$ be an open bounded domain with $V_g(B)>0$ and smooth boundary $\partial B$, $T$ an isoperimetric region, and $0\leq s_T:=Sup\left\{d_{(M,g)}(x, B):\:x\in Supp(||T||)\right\}$. Then there exist positive constants $\varepsilon^*_7=\varepsilon^*_7(B, \xi, g, \partial g, \partial^2 g)>0$ and $\tilde c=\tilde c(n,k,v_0)>0$ such that whenever $Vol_g(T\Delta B)\leq \varepsilon^*_7$, it holds
\begin{equation}\label{Eq:UniformlyBoundedDiameter}
s_{T,g}\leq\tilde cVol_g(Supp(||T||)\setminus B)^{\frac{1}{n}}.
\end{equation}
Furthermore for every $s\in]0, R_3[$, there exists $\varepsilon'_7(s, n, k_0, v_0, B)>0$ with the property that if $T$ is a current solution of the isoperimetric problem with
           $$Vol_g(B\Delta T)<\varepsilon'_7 ,$$
           then 
           $\partial T$ is supported in a tubular neighborhood of thickness $s$ of $\partial B$.
In other words, if $T_j$ is a sequence of isoperimetric regions such that $T_j\rightarrow B$ in flat norm, then $d_{G-H}(T_j,B)\to 0$. 
\end{Lemme}
\begin{Rem} As will appear evident from the proof below, $\frac{1}{\tilde c}=c=c(n,k, v_0)=\frac{C_{Heb}}{4n}>0$ where $C_{Heb}$ denotes the constant appearing in Lemma $3.2$ of $\cite{Hebey}$, that we restate here for completeness's sake. 
\end{Rem}
\begin{Lemme}[Lemma $3.2$ of \cite{Hebey}]\label{Lemma:Hebey3.2}
 Let $(M^n,g)$ be a smooth, complete Riemannian $n$-dimensional manifold with weak bounded geometry. There exist two positive constants $C_{Heb}=C_{Heb}(n,k,v_0)>0$ and $\bar{v}:=\bar{v}(n,k,v_0)>0$, depending only on $n,k$, and $v_0$, such that for any open subset $\Omega$ of $M$ with smooth boundary and compact
closure, if   $V_g(\Omega)\le\bar{v}$, then $C_{Heb}V_g(\Omega)^{\frac{n-1}n}<A_g(\partial\Omega)$.
\end{Lemme}
Now we are ready to prove Lemma \ref{Lemme:UniformlyBoundedDiameter}.

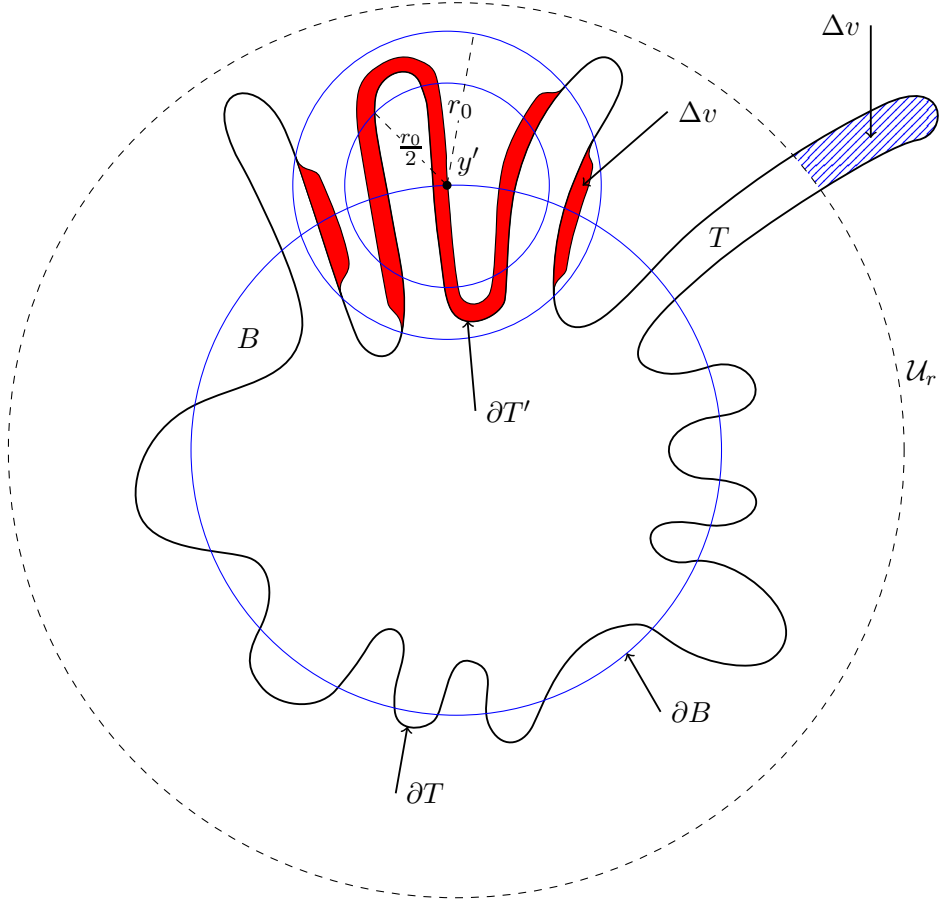
\begin{figure}
\label{fig:Confinement}
\centering
\begin{tikzpicture}[scale=0.3]
\def\curveA{(22.523,61.000)
.. controls (24.507,62.030) and (27.141,62.084) .. (27.397,64.456) 
.. controls (27.555,65.921) and (26.041,69.036) .. (25.048,71.067)
.. controls (24.719,71.740) and (24.429,72.040) .. (24.090,73.094)
.. controls (23.396,75.251) and (25.450,75.541) .. (26.501,73.529) 
.. controls (26.677,73.194) and (26.879,72.690) .. (27.100,72.082) 
.. controls (27.724,70.360) and (28.490,67.798) .. (29.168,65.828) 
.. controls (29.399,65.156) and (29.620,64.553) .. (29.822,64.075) 
.. controls (30.238,63.091) and (31.535,62.723) .. 
($(B)+(252.9:\c)$)
.. controls (32.128,66.021) and (31.706,67.463) .. (31.371,69.380) 
.. controls (30.822,72.521) and (30.473,73.618) .. (30.644,74.705)
.. controls (30.821,75.829) and (32.271,76.307) .. (32.650,75.286)
.. controls (33.031,74.259) and (32.930,74.486) .. (33.060,72.929)
.. controls (33.187,71.399) and (33.222,69.907) .. (33.465,68.213)
.. controls (33.659,66.864) and (33.673,65.893) .. (34.010,65.278)
.. controls (34.551,64.292) and (36.060,64.911) .. (36.345,65.622) 
.. controls (36.500,66.010) and (36.433,68.922) .. (36.810,70.612)
.. controls (37.275,72.693) and (38.197,74.265) .. 
($(B)+(40:\c)$)
.. controls (39.780,75.973) and (40.600,76.955) .. (41.296,76.231)  
.. controls (41.923,75.580) and (41.480,74.737) .. 
($(B)+(17.6:\c)$)
.. controls (40.149,72.601) and (39.897,72.081) .. (39.394,70.844) 
.. controls (38.564,68.804) and (38.506,66.793) .. 
($(B)+(313.9:\c)$)
.. controls (38.624,64.591) and (39.535,63.915) .. (41.036,65.199) 
.. controls (42.505,66.457) and (43.501,67.878) .. (45.135,69.204)
.. controls (47.742,71.319) and (50.289,72.997) .. (53.129,74.366) 
.. controls (56.024,75.762) and (56.078,73.064) .. (54.660,72.812) 
.. controls (53.109,72.537) and (47.410,68.798) .. (45.559,67.198)
.. controls (44.461,66.250) and (41.166,63.580) .. (42.665,62.730)
.. controls (43.456,62.282) and (44.323,62.407) .. (45.442,62.763)
.. controls (47.363,63.374) and (48.552,61.033) .. (46.158,60.615)
.. controls (45.001,60.412) and (45.451,60.513) .. (45.222,60.449)
.. controls (42.750,59.760) and (43.528,57.878) .. (45.260,57.796)
.. controls (48.969,57.622) and (47.525,55.313) .. (45.111,55.766) 
.. controls (41.569,56.432) and (42.782,53.739) .. (44.561,54.113)
.. controls (46.695,54.563) and (50.399,51.306) .. (48.142,49.761)
.. controls (47.418,49.265) and (45.478,49.495) .. (44.120,50.312)
.. controls (42.751,51.135) and (42.853,51.623) .. (41.090,51.138)
.. controls (38.880,50.367) and (38.150,47.905) .. (37.619,46.841)
.. controls (37.151,45.905) and (36.173,45.794) .. (35.690,46.841)
.. controls (35.176,47.959) and (36.325,49.378) .. (34.809,49.706)
.. controls (33.755,49.934) and (33.766,47.509) .. (33.209,47.033)
.. controls (32.720,46.615) and (31.590,46.573) .. (31.503,47.502)
.. controls (31.385,48.768) and (31.824,49.238) .. (31.896,50.117)
.. controls (31.988,51.236) and (31.028,51.557) .. (30.336,50.389)
.. controls (29.876,49.612) and (29.445,49.070) .. (28.859,48.384)
.. controls (27.986,47.454) and (26.274,47.705) .. (25.600,48.560)
.. controls (25.075,49.226) and (24.799,49.750) .. (25.498,51.083)
.. controls (26.242,52.502) and (26.034,53.922) .. (25.057,54.279)
.. controls (24.035,54.652) and (21.070,54.508) .. (20.264,56.042)
.. controls (19.574,57.355) and (20.335,59.864) .. (22.523,61.000)
--cycle;
}
\def\curveB{(28.655,67.340)
.. controls (28.655,67.340) and (28.132,68.975) .. (27.967,69.471)
.. controls (27.946,69.533) and (27.352,71.310) .. (27.102,72.083)
.. controls (27.361,71.296) and (27.968,71.646) .. (28.179,71.177)
.. controls (28.661,70.104) and (29.261,68.650) .. (29.458,67.525)
.. controls (29.601,66.711) and (28.852,66.711) .. (29.161,65.844)
.. controls (28.865,66.735) and (28.655,67.340) .. (28.655,67.340)
--cycle;
}
\def\curveF{(40.31275,72.834879883)
.. controls (39.807,71.829) and (40.397,72.076) .. (40.246,71.471)
.. controls (40.197,71.277) and (40.150,70.906) .. (39.982,70.458)
.. controls (39.731,69.783) and (39.711,69.582) .. (39.486,68.857)
.. controls (39.249,68.092) and (39.144,67.240) .. (39.112,67.092)
.. controls (39.025,66.685) and (38.483,66.695) .. (38.528,65.839)
.. controls (38.528,66.684) and (38.534,67.662) .. (38.647,68.409)
.. controls (38.804,69.441) and (39.074,70.204) .. (39.409,71.098)
.. controls (39.687,71.839) and (39.790,71.770) .. (40.313,72.835)
--cycle;
}
\def\curveC{
(39.028,75.143)
.. controls (38.252,74.168) and (37.077,71.461) .. (37.077,71.336)
.. controls (37.077,70.862) and (36.492,67.771) .. (36.625,67.362)
.. controls (36.717,67.076) and (36.569,66.202) .. (36.295,65.603)
.. controls (36.100,65.179) and (35.807,64.952) .. (35.649,64.911)
.. controls (34.778,64.683) and (34.327,64.852) .. (34.269,64.911)
.. controls (33.912,65.268) and (33.776,65.810) .. (33.650,66.315)
.. controls (33.631,66.392) and (33.397,67.893) .. (33.218,69.687)
.. controls (32.936,72.521) and (32.578,76.087) .. (31.794,75.929)
.. controls (31.498,75.869) and (30.568,75.084) .. (30.580,74.215)
.. controls (30.600,72.792) and (30.999,71.620) .. (31.270,70.265)
.. controls (31.342,69.908) and (31.752,67.380) .. (31.865,66.576)
.. controls (31.940,66.047) and (31.841,65.124) .. (31.841,64.363)
.. controls (31.839,65.149) and (31.298,65.161) .. (31.199,65.577)
.. controls (31.041,66.237) and (30.214,71.105) .. (29.938,72.763)
.. controls (29.898,73.002) and (29.700,74.221) .. (29.904,75.001)
.. controls (30.093,75.728) and (30.678,76.024) .. (30.818,76.119)
.. controls (31.230,76.398) and (31.836,76.629) .. (32.412,76.404)
.. controls (32.823,76.245) and (33.102,76.138) .. (33.269,75.809)
.. controls (33.573,75.211) and (33.803,72.072) .. (33.793,72.002)
.. controls (33.747,71.681) and (33.959,69.503) .. (33.959,69.503)
.. controls (33.983,68.976) and (34.218,66.291) .. (34.411,65.886)
.. controls (34.598,65.495) and (35.243,65.337) .. (35.665,65.932)
.. controls (35.937,66.315) and (35.983,68.633) .. (36.030,69.670)
.. controls (36.034,69.776) and (36.090,70.650) .. (36.264,71.538)
.. controls (36.409,72.276) and (36.754,72.975) .. (36.924,73.299)
.. controls (37.265,73.952) and (37.656,74.553) .. (37.846,74.752)
.. controls (38.208,75.133) and (38.607,74.645) .. (39.028,75.143)
--cycle;
}
\def\curveC{($(B)+(252.9:\c)$) 
.. controls (32.128,66.021) and (31.706,67.463) .. (31.371,69.380) 
.. controls (30.822,72.521) and (30.473,73.618) .. (30.644,74.705)
.. controls (30.821,75.829) and (32.271,76.307) .. (32.650,75.286)
.. controls (33.031,74.259) and (32.930,74.486) .. (33.060,72.929)
.. controls (33.187,71.399) and (33.222,69.907) .. (33.465,68.213)
.. controls (33.659,66.864) and (33.673,65.893) .. (34.010,65.278)
.. controls (34.551,64.292) and (36.060,64.911) .. (36.345,65.622) 
.. controls (36.500,66.010) and (36.433,68.922) .. (36.810,70.612)
.. controls (37.275,72.693) and (38.197,74.265) .. 
($(B)+(40:\c)$)
.. controls  (38.607,74.645) and (38.208,75.133) .. (37.846,74.752)
.. controls (37.656,74.553)  and (37.265,73.952)..  (36.924,73.299)
.. controls (36.754,72.975) and (36.409,72.276) .. (36.264,71.538)
.. controls (36.090,70.650) and (36.034,69.776) .. (36.030,69.670)
.. controls (35.983,68.633) and (35.937,66.315) .. (35.665,65.932)
.. controls (35.243,65.337) and (34.598,65.495) .. (34.411,65.886)
.. controls (34.218,66.291) and (33.983,68.976) .. (33.959,69.503)
.. controls (33.959,69.503) and (33.747,71.681) ..  (33.793,72.002)
.. controls (33.803,72.072) and (33.573,75.211) .. (33.269,75.809)
.. controls (33.102,76.138) and (32.823,76.245) .. (32.412,76.404) 
.. controls (31.836,76.629) and (31.230,76.398) .. (30.818,76.119)
.. controls (30.678,76.024) and (30.093,75.728) .. (29.904,75.001)
.. controls (29.700,74.221) and (29.898,73.002) .. (29.938,72.763) 
.. controls (30.214,71.105)and (31.041,66.237) .. (31.199,65.577)
.. controls (31.298,65.161)and (31.839,65.149) .. 
($(B)+(252.9:\c)$)
 ;
}
\def\curveD{($(B)+(17.6:\c)$)
.. controls (39.807,71.829) and (40.397,72.076) .. (40.246,71.471)
.. controls (40.197,71.277) and (40.150,70.906) .. (39.982,70.458)
.. controls (39.731,69.783) and (39.711,69.582) .. (39.486,68.857)
.. controls (39.249,68.092) and (39.144,67.240) .. (39.112,67.092)
.. controls (39.025,66.685) and (38.483,66.695) .. ($(B)+(313.9:\c)$)
.. controls (38.506,66.793) and (38.564,68.804)  .. (39.394,70.844) 
.. controls (39.897,72.081) and (40.149,72.601) .. ($(B)+(17.6:\c)$)
--cycle;
}

\def\a{19.847};  
\def\b{11.754}; 
\def\c{6.834};   
\def\d{4.533};     
\coordinate (A) at (34.227,59.055);  
\coordinate (B) at ($(A)+(92:\b)$); 
\def\ballA{  (A) circle (\a)};
\def\ballB{  (A) circle (\b)};
\def\ballC{  (B) circle (\c)};
\def\ballD{  (B) circle (\d)};
\begin{scope}
\path [pattern=north east lines, pattern color=blue] \curveA;
\fill[white]    \ballA;
\end{scope}
\begin{scope}
\clip \ballC;
\fill[red] \curveB;
\fill[white] \curveA;
\end{scope}
\draw [line width=0.25mm] \curveA;
\draw [line width=0.15mm, fill=red] \curveB;
\draw [line width=0.15mm,fill=red] \curveC;
\draw [line width=0.15mm, fill=red] \curveD;
\draw[dashed] \ballA;
\draw[blue] \ballB;
\draw[blue] \ballC;
\draw[blue] \ballD;
\fill (B) circle (0.2);
\node [above right] at  (B) {$y'$};
\node [below right ] at  ($(A)+(150:\b)$) {$B$} ;
\node at  ($(A)+(39:15)$) {$T$} ;
\node[right] at   ($(A)+(10:\a)$) {$\mathcal U_r$} ;
\draw[dashed] (B)--($(B)+(80:\c)$);
\fill [white] ($(B)+(80:\c/2)$) circle (0.5);
\node at  ($(B)+(80:\c/2)$) {$r_0$} ; 
\draw[dashed] (B)--($(B)+(135:\d)$);
\fill [white] ($(B)+(135:\d/2)$) circle (0.7);
\node at  ($(B)+(135:\d/2)$) {$\frac{r_0}2$} ; 
\draw[line width=0.7pt, <-]  ($(A)+(37:23)$)--($(A)+(37:23)+(90:5)$) node [left] {$\Delta v$};
\draw[line width=0.7pt, <-]  ($(A)+(64.5:13)$)--($(A)+(64.5:13)+(41:5)$) node [right] {$\Delta v$};
\draw[line width=0.7pt, <-]  ($(A)+(85:5.75)$)-- ($(A)+(85:5.75)+(-85:4)$)  node[right] {$\partial T'$} ;
\draw[line width=0.7pt, <-]  ($(A)+(-100:12.45)$)-- ($(A)+(-100:12.45)+(-100:3.0)$)  node[right] {$\partial T$};
\draw[line width=0.7pt, <-]  ($(A)+(-50:\b)$)-- ($(A)+(-50:\b)+(-60:3.0)$)  node[right] {$\partial B$} ;
\end{tikzpicture}

\caption{Illustration of the Confinement Lemma \ref{Lemme:UniformlyBoundedDiameter}.}
\end{figure}
\begin{Dem}[of Lemma \ref{Lemme:UniformlyBoundedDiameter}] 
Set $V_T(r):=Vol_g(Supp(||T||)\setminus\mathcal{U}_r)$ where $\mathcal{U}_r:=\{x\in M:\:d_g(x, B)\leq r\}$.
Looking at the proof of Theorem $3$ of \cite{NarAsian}, in which boundedness of isoperimetric regions in Riemannian manifolds with bounded geometry is proved (proof that was inspired by preceding works of Frank Morgan \cite{Morgan94} proving boundedness of isoperimetric regions in the Euclidean setting and Manuel Ritor\'e and Cesar Rosales in Euclidean cones \cite{RitoreRosales}), we have that if $Vol_g(T\Delta B)\leq const(n, k, v_0, B, \xi)$, then there exists a positive constant $$c=c(n,k,v_0, Vol_g(B), A_g(\partial B), ||H_{\partial B}||_{\infty,g}, inj_B)>0,$$ such that 
\begin{equation}\label{Eq:UniformlyBoundedDiameter1}
\left(V_T^{\frac{1}{n}}\right)'\leq -c, a.e.\: on\:[0, +\infty[.
\end{equation}
Integrating equation (\ref{Eq:UniformlyBoundedDiameter1}) on the support $[0, s_{T,g}]$ of $V_T$ we get 
\begin{equation}
V_T(s_{T,g})^{\frac{1}{n}}- V_T(0)^{\frac{1}{n}}\leq -cs_{T,g},
\end{equation}
but $V_T(s_{T,g})=0$ and $$V_T(0)=Vol_g(Supp(||T||)\setminus B),$$ hence $s_{T,g}\leq \frac{Vol_g(Supp(||T||)\setminus B)^{\frac{1}{n}}}{c}$, which proves (\ref{Eq:UniformlyBoundedDiameter}). Since, we have trivially that $$Vol_g(Supp(||T||)\setminus Supp(||B||))\leq Vol_g(T-B),$$ we easily finish the proof of the last assertion of the theorem.
To make rigorous the arguments that leads to (\ref{Eq:UniformlyBoundedDiameter1}) we rewrite here the modifications to the proof of Theorem $3$ of \cite{NarAsian} needed here. By Theorem $3$ of \cite{NarAsian}, $T$ have bounded support. Put $A_g(r):=A_g(\partial T\cap (M\setminus\bar{\mathcal{U}_r}))$ it is known that for any $r\in\R\setminus S$ where $S$ is a countable set we have $\mathcal{H}^{n-1}(\partial\mathcal{U}_r\cap\partial^*T)=0$. Fix $\varepsilon^*_7<\sigma_0$ and apply Lemma \ref{Lemma:DeformationLemma2} with $T\setminus\mathcal{U}_r$ in place of $\Omega$ and $\sigma=V_T(r)\ge 0$. In this way we obtain $T'$ such that $Vol_g(T')=Vol_g(T)$ and 
\begin{equation}
A_g(\partial T')\le A_g(\partial (T\setminus\mathcal{U}_r))+C_8A_g(\partial (T\setminus\mathcal{U}_r)\cap Supp(\xi))V_T(r),
\end{equation}
where $C_8=C_8(B, \partial B, \xi, g, \partial g, \partial^2 g)$.
We consider two cases. First $0\le r\le r_0$. Second $r>r_0$. If $0\le r\le r_0$, then $A_g(\partial (T\setminus\mathcal{U}_r)\cap Supp(\xi))=A_g(\partial T\cap Supp(\xi))+A_g(T^{(1)}\cap\partial\mathcal{U}_r\cap Supp(\xi))\le 2A_g(\partial B)+\sup\left\{A_g(\partial\mathcal{U}_r):0\le r\le r_0\right\}=:C_9(B, g, \partial g)$. If $r>r_0$, then $\partial\mathcal{U}_r\cap Supp(\xi)=\emptyset$ and thus $A_g(\partial (T\setminus\mathcal{U}_r)\cap Supp(\xi))\le 2A_g(\partial B)$. Hence in both cases
 \begin{equation}\label{Eq:LemmaConfinement3}
A_g(\partial T')\le A_g(\partial (T\setminus\mathcal{U}_r))+C_{10}V_T(r),
\end{equation} 
where $C_{10}:=C_9C_8=C_{10}(B, \partial B, \xi, g, \partial g, \partial^2 g)$.
We know that $T$ is an isoperimetric region, this implies that 
\begin{equation}\label{Eq:ConfinementFinal}
A_g(\partial T)\leq A_g(\partial T').
\end{equation}
 Hence by \eqref{Eq:LemmaConfinement3}, \eqref{Eq:ConfinementFinal}, and standard slicing theory for currents (or finite perimeter sets or varifolds depending on the taste of the reader) we get 
\begin{equation}\label{Eq:UniformlyBoundedDiameter2}
A_g(r)\le -V_T'(r)+KV_T(r),
\end{equation}
with $K:=C_{10}$.
Assuming $\varepsilon^*_7\le\bar{v}$ we are allowed to apply the isoperimetric inequality for small volumes as in Lemma \ref{Lemma:Hebey3.2} (see Lemma $3.2$ of \cite{Hebey}) to the domain $Supp(||T||)\setminus\mathcal{U}_r$, and again by standard slicing theory, readily follows
\begin{equation}\label{Eq:UniformlyBoundedDiameter3}
-V_T'(r)+A_g(r)\geq C_{Heb}(n,k, v_0)V_T(r)^{\frac{n-1}{n}}.
\end{equation}
Summing (\ref{Eq:UniformlyBoundedDiameter2}) and (\ref{Eq:UniformlyBoundedDiameter3}) we get
\begin{equation}
-\frac{C_{Heb}}{2n}+\frac{K}{2n}(V_T(r))^{\frac{1}{n}}\geq \left(V_T^{\frac{1}{n}}\right)'.
\end{equation}
Therefore, if we choose $\varepsilon^*_7<\left(\frac{C_{Heb}}{2K}\right)^n$ we obtain  $$Vol_g(T\Delta B)<\varepsilon^*_7<\left(\frac{C_{Heb}}{2K}\right)^n=const(B, \partial B, \xi, g, \partial g, \partial^2 g).$$ Remembering that $V_T(r)\leq Vol_g(T\Delta B)$, we obtain 
\begin{equation}
-\frac{C_{Heb}}{4n}=-c\geq \left(V_T^{\frac{1}{n}}\right)'.
\end{equation}
Thus putting 
\begin{eqnarray*}
\varepsilon^*_7 & := & \min\{\left(\frac{C_{Heb}}{2K}\right)^n, \bar{v}, \sigma_0\}\\
& = & \varepsilon^*_7(n,k,v_0, B,Vol_g(B), A_g(\partial B), ||II_{\partial B}||_g, r_{0,g})\\
& = & \varepsilon^*_7(B, \partial B, \xi, g, \partial g, \partial^2 g)>0,
\end{eqnarray*} the proof of the first part of the lemma, i.e., \eqref{Eq:UniformlyBoundedDiameter} is completed. Now to finish the proof we need just to note that what just shown until now permits to us to reduce to the compact case so Lemma \ref{tflat1} applies immediately to a suitable compact neighborhood of $B$ and we can conclude the proof of the lemma taking $\varepsilon'_7<\min\{\varepsilon_7(s),\varepsilon^*_7\}$.  
\end{Dem}
\subsection{Proof of Theorem \ref{tr2}}\label{preuve}
\paragraph{Application of Allard's Theorem} Before starting our proof we recall that the Allard regularity theorem is a regularity theorem with estimates on the the $C^{1,\alpha}$ norm.   
      We give now the proof of Theorem \ref{tr2}. We must show that solutions $T$ of the isoperimetric problem which are close to $B$ in flat norm are graphs of small functions in $C^{1,\alpha}$ norm. Therefore, we fix a real number $\varepsilon>0$ and will find $\varepsilon_0 (\varepsilon)>0$ such that $Vol_g(T\Delta B)< \varepsilon_0 (\varepsilon)$ implies that $\partial T$
is the graph of a function $u$ with $||u||_{\infty}<r(\varepsilon'_0)$, $||u_T||_{C^{1,\alpha}}(\partial B)\le C(\varepsilon'_0)+\varepsilon$. Later on, stronger norms of $u$ will be estimated in terms of $r$ and $\varepsilon$ by Schauder's estimates.\\ 
      
\begin{Dem} Set $\alpha\in ]0,1[$, $\varepsilon\in ]0,1[$, $d=1$ and $p=\frac{n-1}{1-\alpha}$  in the Riemannian Allard's theorem. Consider $R_3=\min\{inj_{(M,g)}, r_{0,g} , \frac{diam_g(B)}{4}\}=R_3(B, \partial B, g,\partial g, \partial^2 g)>0$ as defined in Section \ref{Compvol} and let\\ 
 \begin{eqnarray}\label{Eq:ThmFixedMetric0}
 R & = & \frac12\min\{ R_5,\tilde{R}_1(\varepsilon) ,R_3 , \frac{\tilde{\eta}_1(\varepsilon)}{H_1 [(1+\tilde{\eta}_1(\varepsilon))\omega_{n-1}]^\frac{1}{p}},1\}\\
  & = & R(B, \partial B, g,\partial g, \partial^2 g, \partial^3 g, \partial^4 g, \varepsilon)>0.
 \end{eqnarray}
 Without loss of generality we can assume that $\varepsilon$ is small enough to fill the following conditions
\begin{equation}\label{Eq:ThmFixedMetric}
\varepsilon<\min\left\{\frac{\alpha+1}3,1\right\},
\end{equation}
and
\begin{equation} 
R_0(\varepsilon)<\frac2{3C},
\end{equation}
\begin{equation} 
R_0(\varepsilon)>\frac{6R(\varepsilon)^3}{(1-\varepsilon)}.
\end{equation}
The preceding inequality is possible because by construction we have $R_0(\varepsilon)\sim Const(B, i_g) R(\varepsilon)$ as $\varepsilon\to0$.
 Theorem \ref{allriem} provides us with a constant $\tilde{\eta}_1$ and radius $\tilde{R}_1$ satisfying the conclusion of Theorem \ref{allriem}.  
Then from the comparison Lemma \ref{lr5} applied with $\eta=\tilde{\eta}_1(\varepsilon)$ and $R$ defined by \eqref{Eq:ThmFixedMetric0}, we obtain a  $R_6=R_6(\varepsilon)\in ]\frac{R}{2}, R[$ with the property 
\begin{equation}
||V||(B_g(x,R_6 ))\leq (1+\tilde{\eta}_1 )d\omega_k R_6 ^k .
\end{equation} 
We recall here that $R_6(\varepsilon)\to0$ when $\varepsilon\to0$.
From Lemmas \ref{vpremiere} and \ref{bornecm} we argue that whenever $X\in\mathfrak{X}_c(M)$ with\\ $Supp(X)\subset B_{M}(x,R_6)$,
      \begin{equation}
  \delta\partial T (X)\leq H_1 (A_g(\partial T\cap B(x,R_6 )))^{\frac{1}{p}}||X||_{L^q (\partial T)}.
      \end{equation}
      Hence, an application of comparison Lemma \ref{lr5} allows us to get 
      \begin{eqnarray}
            \delta\partial T (X) & \le & \left\{ H_1 [(1+\tilde{\eta}_1)\omega_{n-1}]^\frac{1}{p}R_6\right\}
R_6^{\frac{n-1}{p}-1}||X||_{L^q (\partial T)}\\ 
& \le & \tilde{\eta}_1 R_6^{\frac{n-1}{p}-1}||X||_{L^q (\partial T)},
      \end{eqnarray}
because  \begin{equation}
  \left\{H_1 [(1+\tilde{\eta}_1)\omega_{n-1}]^\frac{1}{p}R \right\}\leq\tilde{\eta}_1 ,
               \end{equation}  
by the choice of $R$ made in \eqref{Eq:ThmFixedMetric0}. The Riemannian version of Allard's theorem applies with $\tilde{R}=R_6(\eta_{\varepsilon})=R_6(\varepsilon)\to0$, when $\varepsilon\to0^+$. It provides us with a radius $R_0=R_0(\varepsilon)$ such that 
\begin{equation}\label{Eq:ThmFixedMetric2}
(1-\delta_{i_g}r^2)R_6\le R_0\le (1+\delta_{i_g}r^2)R_6,
\end{equation} 
and with a $C^1$ map $F^x:\mathbb{R}^{n-1}\rightarrow M$, for all  $x\in\partial T$, whose image of a neighborhood of the origin is exactly $||i_{g\#}(\partial T)||\cap B_{\mathbb{R}^N}(x,(1-\varepsilon)R_0)$, and whose differential satisfies
           $$||dF^x_z-dF^x_{z'} ||\leq\varepsilon \left(\frac{d_{\mathbb{R}^n}(z,z')}{R_0}\right)^{\alpha}, \forall z,z'\in \mathbb{R}^{n-1},\quad |z|,\,|z'|<R_0,$$
where $i_{g}:(M,g)\to(\R^N,can)$ is the Nash embedding of $M$ in $\R^N$. Now w.l.g. we can assume that $\varepsilon$ is small enough to get 
\begin{equation}\label{Eq:epsilonsmallness}
C(B,\xi,\partial^4 g,\varepsilon) R_0(\varepsilon)<\frac{1}{3}, \frac{c}{2}(1-\varepsilon)R_0(\varepsilon)<\frac13,
\end{equation}
where $C=C(B, g,\partial g,\partial^2 g, \partial^3 g, \partial^4g)$ and $c=c(B, g,\partial g,\partial^2 g, \partial^3 g, \partial^4g)>0$ are constants that will be defined in sequel in equation \eqref{Eq:Thm3.1Proof0}, and 
\begin{equation}\label{Eq:epsilonsmallness0}
\frac12<\frac{\eta}{\eta''}<\frac{3}{2},
\end{equation} 
where $\eta, \eta''$ are defined later in \eqref{Eq:Thm31ProofGlobalGraph} and \eqref{Eq:Thm31ProofGlobalGraph0}.
Pick a radius $$r=r(\varepsilon)\leq R_3(\varepsilon)^*:=\min\left\{\left(\frac{R(\varepsilon)}{2}\right)^3,\frac{(1-\varepsilon)R_0(\varepsilon)}{6}\right\},$$ and set
\begin{small}
$$\varepsilon'_0:=\min\left\{\left[\frac{1}{\tilde c}\left(\frac{R}{2}\right)^3\right]^n,\varepsilon_6,\varepsilon^*_7, \varepsilon'_7 (r)\right\}=\varepsilon'_0 (B, \partial B, \xi, g,\partial g, \partial^2 g, \partial^3 g, \partial^4 g,r(\varepsilon),\varepsilon)>0.$$
\end{small} 
Observe that $\varepsilon'_0\to$ when $\varepsilon\to0$.
Let $T$ be a solution of the isoperimetric problem satisfying 
$$Vol_g(T\Delta B)\leq\varepsilon'_0.$$
The confinement Lemma \ref{Lemme:UniformlyBoundedDiameter} allows us to state that the support of $\partial T$ is inside a tubular neighborhood of thickness $r$.    
      
\paragraph{$\pi|_{\partial T}$ is a local diffeomorphism.}

In what follows $r$ indicates again the thickness of a tubular neighborhood of $\partial B$ in which $\partial T$ is confined, $\pi$ is the projection of $\mathcal{U}_{r}(\partial B)$ on $\partial B$, $\theta$ is the gradient vector of the signed distance function to $\partial B$ and $g_0$ the induced metric by that of  $M$ on $\partial B$.
Let $\varepsilon_0 =min\{\varepsilon '_0 , Vol_g(\{x\in M|d(x,\partial  B)\leq r\})\}$. 
From now on, we assume that $Vol_g(T\Delta B)<\varepsilon_0$. Consider the functions  
\Fonct{f}{]-(1-\varepsilon)R_0(\varepsilon),(1-\varepsilon)R_0(\varepsilon)[}{\mathbb{R}}{t}{d_{(M,g)}(F^x(tv),\partial B),} 
where $ R_0$ is given by the Riemannian Allard's theorem goes to $0$ as $\varepsilon\to0^+$, $v$ is a unit vector in $T_x \partial T$.  Allard's theorem gives a $C^{1,\alpha}$ bound on $F$ therefore  for any $s\in]-(1-\varepsilon)R_0(\varepsilon),(1-\varepsilon)R_0(\varepsilon)[$ 
\begin{eqnarray}\nonumber
|f'(s)-f'(0)| & \le & |\langle dF^x_0(v),\theta_0\rangle_g-\langle dF^x_s(v), \theta_s\rangle_g|\\ \nonumber
 & = & |\langle dF^x_0(v)-dF^x_s(v),\theta_0\rangle_g-\langle dF^x_0(v),\theta_0-\theta_s\rangle_g|\\ \nonumber
  & \le  & \frac{\varepsilon s^\alpha}{((1-\varepsilon)R_0(\varepsilon))^\alpha}+|\langle dF^x_0(v),\theta_0-\theta_s\rangle_g|\\ \nonumber
  & \le & \frac{\varepsilon s^\alpha}{((1-\varepsilon)R_0(\varepsilon))^\alpha}+|dF^x_0(v)| |\theta_0-\theta_s|\\
  \label{Eq:Thm3.1Proof0-}
  & \le &  \frac{\varepsilon s^\alpha}{((1-\varepsilon)R_0(\varepsilon))^\alpha}+c(B,n,\partial^4g)s.
\end{eqnarray}
In particular when $s=(1-\varepsilon)R_0$ we have 
\begin{equation}\label{Eq:ThmFixedMetric3}
|f'(s)-f'(0)|\le\varepsilon+c(B,\partial^4g)(1-\varepsilon)R_0(\varepsilon).
\end{equation}
Thus we conclude that 
\begin{eqnarray}\label{Eq:Thm3.1Proof}
|f'(0)| & \le & \frac{2r(\varepsilon)}{(1-\varepsilon)R_0(\varepsilon)}+\frac{\varepsilon}{\alpha+1}+\frac{c}{2}(1-\varepsilon)R_0(\varepsilon)\\ \label{Eq:Thm3.1Proof0}
& \le & C(B,\xi,\partial^4 g,\varepsilon) R_0(\varepsilon)+\frac{\varepsilon}{\alpha+1}+\frac{c}{2}(1-\varepsilon)R_0(\varepsilon)\\
 & \le & C(B,\xi,\partial^4 g,\varepsilon)<1,
\end{eqnarray}
where $C(\partial^4 g,\varepsilon)\to 0^+$ when $\varepsilon\to0^+$. It is elementary to deduce from \eqref{Eq:Thm3.1Proof0-}, \eqref{Eq:Thm3.1Proof0}, that there exists $C^*_1=C^*(1, B, \partial^4g, \varepsilon, \varepsilon_0)>0$, such that $$||u||_{C^{1,\alpha}_g(\partial B)}\le C^*_1,$$
with $C^*(1, B, \partial^4g, \varepsilon, \varepsilon_0)\to0$, when $\varepsilon\to0$.  

Furthermore, as r gets smaller, the differential of $\pi|_{\partial T}$ gets closer and closer to an isometry.

\paragraph{$\pi|_{\partial T}$ is a global diffeomorphism} 
\begin{Lemme}
            Let $\mathcal{U}$ be a tubular neighborhood of $B$. There exists $\omega\in \Lambda^{n-1}(\mathcal{U})$ such that $d\omega=dVol_g$. 
      \end{Lemme}
      \begin{Dem} Working on each connected component of $B$ we can assume $\mathcal{U}$ being a connected non compact manifold of dimension $n$ implies $H^n (\mathcal{U},\mathbb{R})=0$, see [\cite{God}{Thm. $6.1$  p. 216]}.      
      \end{Dem} 
      
At this stage we just know that $\partial T$ is a locally $C^{1,\alpha}$ regular submanifold of $M$ lying in a tubular neighborhood and possibly composed of infinitely many layers parameterized by a family of functions $u_{T,i}\in C^{1,\alpha}(\partial B)$, $i\in\mathbb{N}$. We have to show that in fact $\partial T$ is a global defined normal graph over $\partial B$. With this aim in mind observe that the family $\{u_{T,i}\}$, as it is easily seen, is actually finite because the area of $\partial T$ is finite and the distortion of are from $\partial B$ to a the $i$-th leaf of $\partial T$ is uniformly bounded by the $C^1$ norm of $g$ and the $C^1$ norm $u_{T,i}$ being $||u_{T,i}||_{C^{1,\alpha}}$ bounded by a constant independent of $i$. By $C^{1,\alpha}$ regularity we can use classical Stoke's Theorem which combined with the preceding lemma gives
\begin{equation}\label{Eq:Thm31ProofGlobalGraph}
Vol_g (T)=\int_{T}d\omega\stackrel{Stokes}{=}\int_{\partial T}\omega=\eta Vol_g (B)=\eta\int_{\partial B}\omega,
\end{equation}
      with $(1-\frac{\varepsilon'_0}{Vol_g (B)})\le\eta\le\left(1+\frac{\varepsilon'_0}{Vol_g (B)}\right)$ close to $1$ for $\varepsilon$ close to $0$. On the other hand denoting by $l:=\#\{u_{T,i}\}$ yields
\begin{equation}\label{Eq:Thm31ProofGlobalGraph0}
\int_{\partial T}\omega=l\eta''\int_{\partial B}\omega\stackrel{Stokes}{=}l\eta''\int_{B}d\omega=l\eta''Vol_g(B),
\end{equation}
      with $\eta'=\eta'(B,\partial g,\varepsilon)=\frac1{\eta''(\varepsilon)}$ close to $1$, when $\varepsilon$ is close to $0$ as $\pi^* (\omega |_{\partial B})$ is close to $\omega |_{\partial T}$ as  $\partial T$ is $C^1$ close to $\partial B$ and 
      $$ \eta ' \int_{\partial T} \omega =\int_{\partial T}\pi^*(\omega |_{\partial B})=l\int_{\partial B}\omega .$$ To be convinced of the first equality of \eqref{Eq:Thm31ProofGlobalGraph0} and its dependence on the $C^1$ norm $u_T$ and on $g$, it is enough to write $\int_{\partial T}\omega$ in local Fermi coordinates based on an open coordinate set $\mathcal{U}\subseteq \partial B$ and then observing that 
\begin{equation}\label{Eq:leafs}   
\int_{\partial T\cap\mathcal{V}}\omega=\sum_{i=1}^l\int_{\U}F^*_{T,i}(\omega),
\end{equation} 
where $l<\infty$ is the number of leaves of $\partial T$ and $F_{T,i}:\mathcal{U}\subseteq\partial B\to\partial T$ is represented in local Fermi coordinates by $(x,u_{T,i}(x))$, i.e., $F_{T,i}(x):=exp_x(u_{T,i}(x)\nu(x))$ for every $x\in\U$, and $\mathcal{V}$ is a cylindrical neighborhood with base $\U$, i.e., $\mathcal{V}:=\U\times ]-r,r[$ of the normal bundle of $\partial B$. Expanding the terms that are in the sum of the right hand side of \eqref{Eq:leafs} we get 
\begin{equation}
\int_{\U}F^*_{T,i}(\omega_g)=\int_{\U} \omega_g(x,u_{T,i})(\partial_1+\frac{\partial u_{T,i}}{\partial x_1}\theta,...,\partial_{n-1}+\frac{\partial u_{T,i}}{\partial x_{n-1}}\theta).
\end{equation} 
Standard computations using basic multilinear algebra and basic elementary inequalities show that 
\begin{equation}\label{Eq:leafs0}
\left|\int_{\U}F^*_{T,i}(\omega_g)-\int_{\U} \omega_g\right|\le C(B,\partial g,||u_{T,i}||_{C^1(\partial B)})\le C(B,\partial g,\varepsilon),
\end{equation}
with $C(B,\partial g,\varepsilon)\to 0$ uniformly as $\varepsilon\to 0^+$, and $C(B,\partial g,\varepsilon)$ being continuous with respect to $g,\partial g$.
From \eqref{Eq:leafs0}, \eqref{Eq:Thm31ProofGlobalGraph}, and \eqref{Eq:Thm31ProofGlobalGraph0} we conclude $$\eta Vol_g(B)=l\eta''Vol_g(B).$$ Having already chosen $\varepsilon$ in \eqref{Eq:epsilonsmallness0} small enough to have $\frac12<\frac{\eta}{\eta''}<\frac{3}{2}$ we establish that $l=1$. In other words we have showed that $\pi|_{\partial T}$ is a global diffeomorphism allowing to set the following definition of the function $u_T$ belonging to $C^{1,\alpha}(\partial B)$ and representing the current (varifold) $\partial T$ as a normal global graph defined over the entire $\partial B$
\begin{equation}
          u_T:=d(\cdot ,\partial B)\circ F\circ (\pi\circ F)^{-1}.
\end{equation}   
\paragraph{$C^{2,\alpha}$ and Higher order Regularity.}

Let us first give a precise definition of the $C^{\ell,\alpha}$ norms.
\begin{Def}
Let $M$ be a compact Riemannian manifold, let $u$ be a function on $M$. We say that $u\in C^{\ell,\alpha}(M,\mathbb{R}^{m})$ if the representative of $u$ in every coordinates chart is of class $C^{\ell,\alpha}$.
\end{Def}
\begin{Def}
Let $u\in C^{\ell,\alpha}(M)$. We set 
$$||u||_{C^{\ell,\alpha}(M)}=\max_l \left\lbrace ||u_{|\Omega_l }||_{C^{\ell,\alpha}(\Omega_l )}\right\rbrace,$$ where 
$||u_{|\Omega_l }||_{C^{\ell,\alpha}(\Omega_l )}:=||u\circ \Theta^{-1}||_{C^{\ell,\alpha}(\mathcal{U}_l )}$ with
$\{\Omega_l \overbrace{\cong}^{\Theta} \mathcal{U}_l \subseteq\mathbb{R}^{n-1}\}$ be a fixed atlas of  $M$. 
\end{Def}

At this point we quote a standard regularity result. The $C^{2,\alpha }$ regularity follows by Schauder estimates, and higher regularity by bootstrap arguments. In order to show that $u$ is more regular we use the same argument used in \cite{Morg1} Proposition $3.3$ p. 5044 as indicated at the end of the proof of \cite{Morg1} Proposition $3.5$ p. 5047. For reader's convenience, we restate here the theorem.
\begin{Prop}[\cite{Morg1} Prop. 3.3]\label{Morg33}
         Let $f$ be a real $C^{1,\alpha}$ function defined on an open set $\Omega $ of $\mathbb{R}^{n-1}$ with the property  
$$\frac{d}{dt}\left[ \int_{\Omega}F(x,f(x)+tg(x),\nabla (f(x)+tg(x)))dx\right] _{t=0}=0$$ whenever $g$ is a $C^1$ function with $Supp(g)\subset\subset\Omega$. Assume $F$ and $\frac{\partial F}{\partial f_i }$ are $C^{\ell-1,\alpha}$ for some $l\ge2$, $\alpha\in]0,1[$, and $F$ is elliptic, i.e.  the matrix $\frac{\partial F}{\partial f_i \partial f_j}$ is positive definite.\\
Then\\
$f$ is $C^{\ell,\alpha}$.
\end{Prop} 
\begin{Dem}
The proof can be found in Proposition 3.3 of \cite{Morg1}.
\end{Dem} 

In local coordinates, we can see $\partial T$ locally like the graph of a function $f$ of class $C^{1,\alpha}$.\\
For smooth variations $g$ with compact support the area functional $\mathcal{A}(f):=\int A(x,f,\nabla f (x))dx$ and the volume functional $\mathcal{V}(f):=\int V(x,f(x))dx$ satisfy the relation:
\begin{equation}
          \frac{d}{dt}\left[ \mathcal{A}(f+tg)-\lambda\mathcal{V}(f+tg)\right]|_{t=0}=0 
\end{equation}
for some Lagrange multiplier $\lambda$ that is the mean curvature of $\partial T$. The functional $\mathcal{A}-\lambda\mathcal{V}$ then satisfies the regularity and ellipticity assumptions required by Proposition \ref{Morg33}, hence $\partial T$ is as regular as possible and at least of class $C^{2,\alpha}$, which implies  by an application of the implicit function theorem that $F$ given by Allard's theorem belongs to $C^{2,\alpha}$ and therefore that $u$ is also of class $C^{2,\alpha}$.\\
In other words, there exists $\tilde{F}$ of class $C^{2,\alpha}$ such that
$$u=d(\cdot ,\partial B)\circ \tilde{F}\circ (\pi\circ \tilde{F})^{-1},$$ and we conclude that $u$ is of class $C^{2,\alpha}$. By a standard bootstrap argument we conclude that $u$ is $C^{\infty}$, since $g$ is $C^{\infty}$.\\
\paragraph{$C^{2,\alpha}$ and higher order estimates.} 
Now we are in a position to exploit formula (\ref{h-norm1hyp}) for the mean curvature of a normal graph, represented as a function $u$ defined on $\partial B$. This allows to estimate the $C^{1,\alpha}$ norm and $C^{2,\alpha}$ norm of $u$. Straightforward computations will show that the $C^{2,\alpha}$ norm of $u$ goes to zero when $r\rightarrow 0$. 
We now give some details of these calculations. We consider a system of Fermi coordinates $(r,x)$ centerd at a point $p\in\partial B $, with $x$ normal coordinates on an open set of $\partial B$ centered in $p$. Let
$$ u_i :=\frac{\partial u}{\partial x^i},\quad u_{ij} :=  \frac{\partial^2 u}{\partial x^i x^j},$$
\begin{equation}
         g:=dt^2 +g_{ij}(t,x )dx ^i dx ^j ,
\end{equation}  
\begin{equation}
         ||\nabla_{g_u}u||^2 _{g_u}=g^{ij}(u,x )u_i u_j ,
\end{equation}
\begin{equation}
\begin{array}{ccl}
        \nabla_{g_u }W_u & = & -\frac{1}{2}\frac{1}{\sqrt{(1+||\nabla u||^2 )^3}}\left\lbrace \frac{\partial}{\partial r}g^{lj}(u,x )u_i u_j u_l  \right\rbrace\\
    & - & \frac{1}{2}\frac{1}{\sqrt{(1+||\nabla u||^2 )^3}}\left\lbrace \frac{\partial}{\partial x^i}g^{jl}(u,x)u_j u_l +2g^{lj}(u,x )u_i u_{ij} u_l \right\rbrace g^{im}\frac{\partial}{\partial x^m} 
\end{array}
\end{equation}
\begin{equation}
\frac{1}{W_u}\left[ div_{\partial B^r}\left( \nabla_{g_u} u\right)\right] _{|r=u}=\left[ \frac{1}{W_u}g^{ij}(u,x )+f^{ij}(x,u,\nabla u)\right] u_{ij}+f(x ,u,\nabla u) .
\end{equation}
Notice that $ f(x ,u,\nabla u)$ , $ f^{ij}(x,u,\nabla u) $ $ \rightarrow 0$ ,  $||u||_{C^1 }\rightarrow 0$. The functions  
$$f,f^{ij}:\Omega\times\mathbb{R}\times\mathbb{R}^{n-1}\to\mathbb{R}$$  
have the same regularity than the metric with respect to variables $x$, $y$ and they are of class $C^{\infty}$ with respect to $z$.
We carry analogous calculations for the remaining $4$ terms of formula (\ref{h-norm1hyp}).
After these straightforward standard computations we obtain the following expression for the constant mean curvature equation of a normal graph based on a hypersurface
\begin{equation}
\left[ \frac{1}{W_u}g^{ij}(u,x )+l^{ij}(x,u,\nabla u)\right] u_{ij}=h=h_1 +h_2,
\end{equation}
where $h_1=H_{\nu}^{\partial T}-\frac{1}{W_u}H_{-\theta}^{{\partial B}^u}$ and $h_2 =h_2(x, u,\nabla u)$ satisfying $$||h_2||_{\infty}\rightarrow0,$$ when $||u||_{C^{1,\alpha}}\rightarrow 0$. Moreover
$$h_1,h_2:\Omega\times\mathbb{R}\times\mathbb{R}^{n-1}\to\mathbb{R}$$
have the same regularity as the Levi-Civita connection with respect to variables $x$, $y$ and are of class $C^{\infty}$ with respect to $z$. 
If $k\leq K^{M}\leq\delta$ (which is guaranteed by the fact that we are in a compact neighborhood of $B$ that is itself compact), then by Heintze-Karcher's theorem see (\ref{Eq:HKCorollary}) we get 
\begin{equation}
\left(c_{\delta}(u)-\beta s_{\delta}(u)\right)^2g_0 \leq g(u,x )\leq\left(c_{\delta}(u)+\beta s_{\delta}(u)\right)^2g_0, 
\end{equation}
and so
\begin{equation}
\frac{g_0 ^{-1}}{\left(c_{\delta}(u)-\beta s_{\delta}(u)\right)^2} \leq g(u,x )^{-1}\leq\frac{g_0 ^{-1}}{\left(c_{\delta}(u)+\beta s_{\delta}(u)\right)^2},
\end{equation}
where $g_0$ is the metric $g$ restricted to $\partial B$.
Consequently, there are $0<A_1\leq A_2$ for which
\begin{equation}
\frac{A_1 I_{n-1}}{\left(c_{\delta}(u)+\beta s_{\delta}(u)\right)^2}\leq g(u,x )^{-1}\leq\frac{A_2 I_{n-1}}{\left(c_{\delta}(u)-\beta s_{\delta}(u)\right)^2},
\end{equation}
hence the equation
$$Lu:=a^{ij}u_{ij}=\tilde{h}(x ),$$
with $a^{ij}(x):= \frac{1}{W_u}g^{ij}(u,x )+l^{ij}(x,u,\nabla u)$, $\tilde{h}(x )=h(x , u(x), \nabla u(x ))$ is uniformly elliptic as  the $l^{ij}\rightarrow 0$ when $||u||_{C^1}\searrow 0$ ($r\searrow0$). Using classical Schauder interior estimates for linear elliptic partial differential equations, e.g. Theorem $6.2$ and Corollary $6.3$ of \cite{GT} applied to a fixed covering by charts of $\partial B$ and taking as $\Lambda$ of $(6.13)$ of \cite{GT} a uniform fixed upper bound of $C(1,B,\partial^4g, \varepsilon)$ (a constant that is taken as an upper bound of $||u_T||_{C^{1,\alpha}}$ in the statement of Theorem $\ref{tr2}$) multiplied by a constant that depends on the diameter of $\partial B$ (i.e., $C^0$ on the metric $g$) it is not too hard to check, that 
\begin{eqnarray}\label{Eq:Thm3.1SchauderEstimates}
||u_T||_{C^{2,\alpha}_g(\partial B)} & \le & \tilde C(\partial B, g, \partial g, \partial^2 g, C(\bar{\varepsilon}))||u||_{C^{1,\alpha}_g(\partial B)}\\
 & \le & \tilde CC(1,B,\partial^4g, \varepsilon, \varepsilon_0)=:C(2,B,\partial^4g, \varepsilon, \varepsilon_0).
\end{eqnarray}
Thus we can choose $C^*(2,B,\partial^4g, \varepsilon, \varepsilon_0):=\tilde CC(1,B,\partial^4g, \varepsilon, \varepsilon_0)>0$ in the statement of Theorem $\ref{tr2}$. Now deriving equation \eqref{h-norm1hyp} and iterating a suitable number of times (in fact $m-1$ times) this Schauder interior estimates argument, we easily obtain our constants $C^*(m,B,\partial^4g, \varepsilon, \varepsilon_0, ||g||_{m, \alpha})$ for any $m$. This shows that for any $m$ the constants $C(m,B,\partial^4g, \varepsilon, \varepsilon_0)>0$ are small when $\varepsilon,\varepsilon_0$ are small. Now it is trivial to deduce the remaining parts of the statement of the theorem. We just point out that the fact that in case $B$ is the limit in flat norm of a sequence of isoperimetric regions implies that $B$ is also an isoperimetric region is due to fact that as a consequence of Theorem $2$ of \cite{FloresNardulli2016} $I_{(M,g)}$ is continuous when $(M,g)$ is of bounded geometry. With this last remark we complete the proof of the theorem.  
\end{Dem}
\subsection{Some refined mean curvature estimates}
In this last section we give an effective estimate of the difference of the mean curvature vector of $\partial T$ and $\partial B$. We give a more geometric characterizations of the explicit estimates of $$||H^{\partial T}-H^{\partial B}||_{C^{0}},$$ in terms of the geometric data of the isometric embedding of $\partial B$ into $(M,g)$ and the ambient metric $g$. This is not relevant for the sequel but it has an interest in itself; for this reason we included it here. From the preceding theorem we know that the interior normals to $\partial T$ converge to the interior normals of $\partial B$ by $C^1$ convergence and that the mean curvature vectors of $\partial T$ converge to the mean curvature vectors of $\partial B$ by $C^2$ convergence. We want to compare the mean curvature of $\partial T$ with the mean curvature of a touching inscribed equidistant hypersurface. This is possible only when the mean curvatures point in the same direction, and unfortunately when the mean curvature of $\partial B$ is $0$ or changes sign we are not able to do such a comparison. However when the mean curvature of $\partial B$ does not change direction and is not zero, Schauder estimates of Section \ref{preuve} show provided $\varepsilon$ is small enough, that $H^{\partial T}$ and $H^{\partial B_r}$ have the same direction at points of contact. So in particular for $\varepsilon$ small enough the mean curvature vector of $\partial T$ at a maximum point $x_0$ of $u_T$ and at a minimum point $x_1$ of $u_T$ point in the same direction of that of suitable circumscribed and inscribed tangent equidistant hypersurfaces of $\partial B$. 
We prove the following lemma.
\begin{Lemme}\label{Lemma:Schauderestimates} There exists $b_3 (s)$ such that whenever $y\in\partial B$, 
\begin{equation}
 |H_{\theta}^{{\partial B}^s }(y)-H_{\theta}^{\partial B}(y)|\leq b_3 (s),
\end{equation}
where $\partial B^s$ is the equidistant hypersurface at signed distance $s$ from $\partial B$.
\end{Lemme}
\begin{Dem} Let
$$b'_3 (s,y):=|\sum_{i=1}^{n-1}ctg_{\delta}(s+c_1 (y,\lambda_i (y)))-H_\theta ^{\partial B}(y)|,$$
$$b''_3 (s,y):=|a_k (s+c_2 (y,H_\theta ^{\partial B}(y)))-H_\theta ^{\partial B}(y)|,$$
 $$b_3 (s,y) :=Max\left\lbrace b'_3 (s,y),b''_3 (s,y)\right\rbrace,$$ 
where $ ctg_{\delta}(c_1(x,s))= s$, $c_1(x,s)\in]0,\frac{\pi}{\sqrt{\delta}}[ $, and $ctg_{k}(c_2(x,s))=s$,\\ if $ s>\sqrt{-k}$, $tg_{k}(c_2(x,s))=s $, if $s<\sqrt{-k}$ and $c_2(x,\sqrt{-k})=\sqrt{-k}$
$$a_k (s)=\left\lbrace \begin{array}{lll}
                                         ctg_k (s) & , &  s>\sqrt{-k}\\
                                         \sqrt{-k} & , & s=\sqrt{-k}\\
                                         tg_k (s) & , & s<\sqrt{-k}\\
            \end{array}\right.$$  
We find $b_3 (s):=||b_3 (s,y) ||_{\infty ,\partial B}$.          
\end{Dem}
\begin{Rem} $b_3 (s)\rightarrow 0$, when $s\rightarrow 0$. 
\end{Rem}

\begin{Thm}[The Comparison Principle for Mean Curvatures]\label{Thm:ComparisonPrinciple}
Let $B_1 $ and $B_2$ being two submanifolds with boundary, of dimension $n$ of $M$, $B_1\subseteq B_2$, with $\{x\}=\partial B_1 \cap\partial B_2$, for a single point $x\in M$, with the mean curvature vector that points in the same direction.\\
Then\\
$<H^{\partial B_1 }(x),\nu_{ext}>\leq <H^{\partial B_2 }(x),\nu_{ext}> $ 
\end{Thm}
\begin{Dem}
\cite{Alex}
\end{Dem}

\begin{Lemme}\label{lp} Let $\partial T_j$ be a sequence of normal graphs of $C^{2,\alpha}$ functions $u_j$ over $\partial B$. Assume that $u_j$ satisfies the constant mean curvature equation, $||u_j||_{\infty}$ converges to $0$ as $j\rightarrow+\infty$ and that $\partial T_j$ and $\partial B$ have mean curvature vectors such that $<H_{\partial B},\theta>$ and $<H_{\partial T_j},\theta>$ have the same sign. Then  
\begin{equation}\label{Eq:lp}
\left|H^{\partial T_j}-H_{\theta}^{{\partial B}^u}\right|\leq\max\{|b_3(u(x_1))|, |b_3(u(x_2))|\}\rightarrow 0,
\end{equation} 
when $j\rightarrow +\infty$. In particular, $\eqref{Eq:lp}$ holds if the sequence $(T_j)$ and $B$ satisfy the hypothesis of Theorem \ref{tr2}. 
\end{Lemme}
\begin{Dem}
Let $x_1 $, $x_2 \in \partial B$ be defined as $u(x_2 ):=Max_{x\in\partial B}\{u(x )\}$ and $u(x_1 ):=Min_{x\in\partial B}\{u(x )\}$.\\ Then
$$B^{u(x_1 )} \subseteq T\subseteq B^{u(x_2 )}$$
and $B^{u(x_1 )}$, $B^{u(x_2 )}$ have smooth boundary and are tangent to $\partial T$ at $p_1=(x_1 , u(x_1 ))$ and  
$p_2=(x_2 , u(x_2 ))$.  
We deduce then, by the comparison principle applied to $B^{u(x_1 )}$, $T$, $B^{u(x_2 )}$ that
\begin{equation}\label{ep1}
\left|H_{\nu}^{\partial T}(x )-H_\theta^{\partial B}\right|\leq\max\{|b_3(u(x_1))|, |b_3(u(x_2))|\}.
\end{equation}
\end{Dem}
\section{Proof of Theorem \ref{T4}: Normal Graph Theorem with variable metrics}\label{Section:VariableMetrics}
In this section we present the proof of our main Theorem \ref{T4}. We begin by summarizing results of Gromov \cite{Gr2}, p. 118 that we will need. We assume that the reader is familiar with the notions of fibration, vector bundle, jet bundle of a fibration, and partial differential operator. For these topics one can consult various basic texts such as \cite{Hirsch}. For a more advanced treatment relevant for our purposes we strongly recommend \cite{Eliashberg}, \cite{Spring}, and obviously the treatise \cite{Gr2}. We follow closely the treatment given in \cite{Gr2}.
\begin{Def}
Let $p:X^{n+q}\rightarrow V^n$ be a smooth fibration and let $\pi:G\rightarrow V$ be a smooth vector bundle. We denote by $\mathcal{X}^{\alpha}$ and $\mathcal{G}^{\alpha}$ the spaces of $C^{\alpha}$-sections of the fibrations $p$ and $\pi$ for $\alpha\in\mathbb{N}\mathring{\cup}\{\infty\}$, respectively. We say that $\mathcal{D}:\mathcal{X}^{\alpha}\to\mathcal{G}^{\alpha}$, is a \textbf{differential operator of order $r$}, if there exists $\Delta:X^{(r)}\rightarrow G$ $($here $X^{(r)}$ is the space of $r-jet$ of sections of $p$$)$ such that $\mathcal{D}(\sigma)=\Delta\circ J_\sigma^r$, for every $C^r$-section $\sigma$ of $p$. $\mathcal{D}$ is said \textbf{$C^{\alpha}$-smooth}, if $\Delta$ is $C^{\alpha}$-smooth. We assume in the sequel that $\mathcal{D}$ is $C^{\infty}$ and so the maps $\mathcal{D}:\mathcal{X}^{r+\alpha}\rightarrow\mathcal{G}^{\alpha}$ are continuous with respect to the usual compact-open and fine topologies. Here a typical neighborhood of a section $\sigma_1\in\Gamma(\xi)$ of a fibration $\xi:E\rightarrow V$ in the \textbf{$C^0$-fine topology} is of the form $\mathcal{U}_{\varepsilon}(\sigma_1):=\{\sigma_2\in\Gamma(\xi):\:d_E(\sigma_1(v),\sigma_2(v))<\varepsilon(v)\}$, where $\varepsilon(v)\in C^0(V,[0,+\infty[)$ and $d_E$ is a metric on $E$.  
\end{Def}
\begin{Rem} There are several equivalent definitions in the literature of our fine topology known also as the Whitney strong topology; for more details see $\cite{Hirsch}$ p. $59$ and the entire content of Chapter $2$ of the same book or $\cite{Spring}$ p. $9$.
\end{Rem}
\begin{Example}\label{Ex:VarMetric} Let $G$ be the bundle of symmetric bilinear forms over the manifold $V$, and let $(W, h)$ be a manifold endowed with a quadratic differential form $h$. We consider the trivial fibration $\xi:X=W\times V\rightarrow V$. As usual we identify every map $V\rightarrow W$  with a section of $\xi$. We obtain a first order partial differential operator $\mathcal{D}$ if we define $\mathcal{D}(\sigma):=\sigma^*(h)$. $\mathcal{D}$ is $C^{\infty}$ if $h$ is $C^{\infty}$. 
\end{Example}
\begin{Def}\label{Def:Petersen} For any $m\in\mathbb{N}$, $\alpha\in [0, 1]$, a sequence of pointed smooth complete Riemannian manifolds is said to \textbf{converge in the pointed $C^{m,\alpha}$, respectively $C^{m}$ topology to a smooth manifold $M$} (denoted $(M_i, p_i, g_i)\rightarrow (M,p,g)$), if for every $R > 0$ we can find a domain $\Omega_R$ with $B(p,R)\subseteq\Omega_R\subseteq M$, a natural number $\nu_R\in\mathbb{N}$, and $C^{m+1}$ embeddings $F_{i,R}:\Omega_R\rightarrow M_i$, for large $i\geq\nu_R$ such that $B(p_i,R)\subseteq F_{i,R} (\Omega_R)$ and $F_{i,R}^*(g_i)\rightarrow g$ on $\Omega_R$ in the $C^{m,\alpha}$, respectively $C^m$ topology. 
\end{Def}
\begin{Rem} 
As it easy to check when the manifolds are compact, pointed convergence is independent of the base point, so we can speak just of convergence without making any reference to the word pointed. 
\end{Rem}
We define now the fine topologies needed to state the continuity results with respect to the metric deducible from Nash's imbedding theorem. Following \cite{Hirsch} we give the following definition.
\begin{Def}[\cite{Gr2} page 18]\label{Def:WhitneyTopology} Let $(X, \tau_X)$ and $(Y,\tau_Y)$ be arbitrary topological spaces, denotes by $\tau_X\times\tau_Y$ the product space topology. Let $f\in C^0(X, Y)$, $\Gamma_f\subseteq X\times Y$ be the graph of $f$, $\mathcal{U}(f,W):=\{g\in C^0(X,Y):\:\Gamma_g\subseteq W\}$. The family $\{\mathcal{U}(f,W)\}$ with $(f,W)\in C^0(X,Y)\times(\tau_X\times\tau_Y)$ forms a base for a topology $\tau$ that is called in the literature \textbf{strong topology} or \textbf{fine topology} or \textbf{Whitney topology}. We chose here to call $\tau$ the fine topology. We define $C^0_S(X,Y):=(C^0(X,Y),\tau)$. The \textbf{fine-$C^r$-topology} in $\Gamma^r(V,X)$ is the relative topology induced from the fine topology in $\Gamma^r(V,X)\to C^0(V,X^{(r)})$ by the injection $f\mapsto J_f^r$ onto the space of holonomic sections $\Gamma(V, X^{(r)})$. 
\end{Def}
\begin{Def}[\cite{Gr2} page $8$] Let $f\in C^2(V^n,\R^q)$, and $x$ a chart centered at $p\in V$. Denote by $T_f^2(V,p)\leq T_{f(p)}(\R^q)$ the subspace spanned by the vectors $\frac{\partial f}{\partial x^i}(p), \frac{\partial^2 f}{\partial x^ix^j}(p)$, for $i,j\in\{1,...,n\}$. We say that $f$ is a \textbf{free immersion} if all the preceding vectors form a linearly independent set of vectors. 
\end{Def}
The following theorem is attributed in \cite{Gr2} page $116$ to John Nash.
\begin{Thm}[
compare \cite{Gr2}, page $116$]\label{Thm:Nash0} Let $\mathcal{D}$ be the operator of Example \ref{Ex:VarMetric}, with $(W,h)=(\mathbb{R}^q, \delta)$ where $\delta$ is the canonical Euclidean metric. Then over the space of free maps $V\rightarrow\R^q$, $\mathcal{D}$ admits an infinitesimal inversion $M$ of defect $d=2$ and of order $s=0$.
\end{Thm}
Roughly speaking this last theorem asserts that the differential operator $\mathcal{D}$ of degree $1$ of Example \ref{Ex:VarMetric} has its differential invertible on free immersions that are two ($d=2$) times differentiable with inverse a differential operator of order $s=0$. The optimal version of this theorem could be with $d=1$,  and this justifies the reason for the use of the word defect. For an account of the proof of this last result and for the rigorous definitions needed to understand its statement we refer the reader to the book \cite{Gr2} pages $116$-$117$. 
As a consequence of Theorem \ref{Thm:Nash0}, $(4)$ of Main Theorem pages $117$-$118$ 
  we have the following remarkable results, whose statement is just the statement of $(4)$ of \cite{Gr2} page $118$ specialized to the case of the differential operator $\mathcal{D}$ described in example \ref{Ex:VarMetric}, with $s=0$, $d=2$, $r=1$, $\sigma_0=\sigma_1=\eta_1=3$. 
\begin{Thm}[ $(4)$ \cite{Gr2} page $118$]\label{MainGromov}  Let $n<N$, and let $i_{g_{\infty}}:(M^n, g_{\infty})\to(\R^N,\delta)$ be a free $C^4$ isometric immersion. Then for every $\alpha\geq 4$, $\alpha\in\{0,1,..., \infty\}$ there exists a fine $C^{\alpha}$-neighborhood $\U^\alpha$ of $g_{\infty}$ such that for every $g\in\U$ there exists an isometric immersion $i_g:(M, g)\to(\R^N,\delta)$ of class $C^\sigma$, for any integer $\sigma<3$. Moreover such immersions can be chosen such that $i_g\to i_{\infty}$ in the $C^{\sigma}$-fine topology, as $g\to g_{\infty}$ in $C^{\alpha}$-topology. In particular $i_g\to i_{g_{\infty}}$ in $C^2$ topology when $g\to g_{\infty}$ in $C^{\alpha}$-topology, and so also the second fundamental forms $II_{i_g}\to II_{i_{g_{\infty}}}$, $g\to g_{\infty}$ in $C^{\alpha}$-topology.    
\end{Thm}
Observe that our choice of $r=1,\sigma_0=\sigma_1=\eta_1=3$ are the weakest possible in the range of integers, because it has to be  $\max\{d, 2r+s\}=:\bar{s}<\sigma_0\le\sigma_1\le\eta_1$. On the other hand if we take $(4)$ of page $118$ with $\sigma_0=3$, $\eta_1=\sigma_1=\infty$, $r=1$, $s=0$, $d=2$, and $\mathcal{D}$ the differential operator of Example \ref{Ex:VarMetric}, we have the following theorem.
\begin{Thm}[Nash in \cite{Nash}]\label{Thm:MainNash} If there exists a free $C^\infty$ isometric immersion $i_g:(M,g)\to(\R^N, \delta)$, then there exists a  $C^3$ neighborhood $\mathcal{U}^3$ of $g$ such that for every $\sigma\ge3$ and $h\in\mathcal{U}^3\cap \Gamma^\sigma(G)$, where $G$ is the bundle of symmetric bilinear forms on $M$ of Example \ref{Ex:VarMetric}, there exists a $C^\sigma$ isometric immersion $i_h:(M,g)\to(\R^N, \delta)$.
\end{Thm}
\begin{Thm}[Imbedding Theorem \cite{Gr2} page $223$]\label{Thm:VarMetricImbedding1} Every Riemannian $C^{\alpha}$-manifold $V^n$, $2<\alpha\leq\infty$, admits a free isometric $C^{\alpha}$-imbedding $i:V\to\R^N$, with $N=n^2+10n+3$.
\end{Thm}
Observe that our $i_h$ corresponds to the $\mathcal{D}_{i_g}^{-1}(h)$ in the notation of \cite{Gr2} and these corresponds to the same $i_g$ of Theorem \ref{MainGromov}.
Now we are ready to achieve the proof of our main result.

\begin{Dem}[of Theorem \ref{T4}]
Take the manifold $(M,g_{\infty})$ and apply Theorem \ref{Thm:VarMetricImbedding1} to $(M,g_{\infty})$ with $\alpha=\infty$ to obtain a free isometric $C^{\infty}$-imbedding  $i_{\infty}$ for $(M, g_{\infty})$ fixed.
Furthermore, an application of Theorem \ref{MainGromov} allows us to obtain $C^\infty$ free isometric embeddings  $i_{g_j}$ of $(M, g_j)$ into $(\mathbb{R}^N, \delta)$ close in the $C^2$ fine topology (see \cite{Gr2} p. 18) to $i_{g_{\infty}}$. If $M$ is compact the fine topology and the usual topology of convergence on compact sets are the same, so the $C^2$-fine-topology of Definition $\ref{Def:WhitneyTopology}$ is the same as the $C^2$-topology of Definition $\ref{Def:Petersen}$. If $M$ is not compact the explicit computation of the constants involved in Lemma \ref{tflat1} shows that they depend continuously on $Vol_g(B), H_{\partial B,g}, A_g(\partial B), inj_{B,g}$ and that when $j$ varies the numbers $n,k,v_0$ do not vary. Namely $Vol_g(B), A_g(\partial B), diam_g(B), |\nabla_g\xi|_g$ depend just on $g$ and they are continuous in $C^0$-topology, $H_{g,\partial B}$ depends on $g$ and the first derivatives of $g$, and moreover they depend on them continuously. For what concerns the injectivity radius of $(M,g)$ we have $inj_g\to inj_{g_{\infty}}$ in $C^2$-topology as is proved in the Theorem (there is no number in the paper of Sakai) of page $91$ of Sakai \cite{SakaiContinuityOfInjectivityRadius} and by Theorem \ref{MainGromov} $\beta_{i_g}\to \beta_{i_{g_\infty}}$. So the constants $\varepsilon^*_{7,j}$, and $\varepsilon'_{7,j}$ of Lemma \ref{Lemme:UniformlyBoundedDiameter} applied to the metrics $g_j$ satisfy $\varepsilon'_{7,j}\to\varepsilon_{7,\infty}$ and are obviously uniformly bounded below.
Hence we can put all the $T_j$ inside a big compact set $\tilde{B}\subseteq M$ such that $diam_{g_j}(\tilde{B}), diam_{g_\infty}(\tilde{B})\leq const.$ uniformly with respect to $j$. So we are reduced to the case when $M$ is a compact manifold and this just requires bounded geometry and $C^2$ convergence of the metrics. Now we are in position to apply our Theorem  $\ref{tr2}$ and obtain that $\partial T_j$ is a normal graph of a function $u_{T_j}\in C^{2,\alpha}(\partial B)$, moreover for every $\varepsilon>0$ there is $j_\varepsilon$ such that for every $j\ge j_\varepsilon$ it holds 
\begin{equation}\label{Eq:explained}
||u_{T_j}||_{C^{2,\alpha}}\le C=C(\varepsilon),
\end{equation} 
with $C$ depending just on $\varepsilon$ and satisfying 
\begin{equation}\label{Eq:explained0}
C(\varepsilon)\to0, \:\varepsilon\to0.
\end{equation} 
 To check the validity of \eqref{Eq:explained}, \eqref{Eq:explained0} observe that the explicit calculations made in the proof of Theorem  $\ref{tr2}$ the constants on which the estimates of Theorem $\ref{tr2}$ depend are divided into two disjoint finite sets $\mathcal{C}:=\mathcal{A}\mathring{\cup}\mathcal{B}$, satisfying the property that if $c\in\mathcal{A}$ then $c=c(B, \xi, g, \partial g, \partial^2g, i_g)=c(B, g, \partial g, \partial^2g,\partial^3g,\partial^4g)$, and if $c\in\mathcal{B}$ then $c=c(B, \xi, g, \partial g, \partial^2g)$ does not depends on $i_g$. Furthermore the constants $c\in\mathcal{A}$ depends continuously in $C^4$ topology on the metric and the constants $c\in\mathcal{B}$ depends continuously in $C^4$ topology on the metric. The dependence on $B$ and $\xi$ of the constants means the dependence on $B$ and $\xi$ differentiable objects independent of $g$. To differentiate between quantities that depend on the metric also we indicate it explicitly. By Theorem \ref{MainGromov} it follows easily that a typical constant $c(B, g, \partial g, \partial^2g, i_g)=c(B, g, \partial g, \partial^2g,\partial^3g,\partial^4g)$. This is the reason for requiring $C^4$ convergence. By the way this is just a temporarily technical obstacle due to the fact that the version of the Allard regularity theorem that we use, needs the Nash isometric embedding theorem. Now appears clear that $C^4$ convergence of the metric implies that the constants of Theorem $\ref{tr2}$ could be chosen independently of $j$. This last fact combined with \eqref{Eq:explained}, \eqref{Eq:explained0} readily yields $||u_{T_j}||_{C^{2,\alpha}}\to0$. Finally using the last part of Theorem $\ref{tr2}$ about higher order norm estimates we finish the proof of the theorem.  
\end{Dem}
\begin{Rem} We need this unpleasant $C^4$ convergence in Theorem $\ref{tr2}$ because of the dependence of the constants involved on the imbedding $i_g$ through the bounds on $\beta_g$ that are continuous with respect to $g$ only if the imbedding $i_g$ are continuous in $C^2$ topology. Unfortunately we can ensure the $C^2$ continuity of the embeddings $i_g$ only in case of $C^4$ convergence of the metrics $g_j$. However, it is still possible to drop the hypothesis of $C^4$ convergence and replace it by $C^2$ convergence if we are concerned just with $C^{2,\alpha}$ $($remember that we assumed that $M$ is a smooth differentiable manifold$)$ convergence of the $u_{T,j}$ as prescribed by the Allard's regularity Theorem. Consult Remark \ref{Rem:Allardintrinsic} on this last issue. 
\end{Rem}
\begin{Rem}\label{Rem:Allardintrinsic} As a final remark we expect that with a slight but cumbersome modification of the arguments contained in the proof of Theorem \ref{tr2}, Theorem \ref{T4} is true also if we replace $C^4$ convergence by $C^2$ convergence of the metrics $g_j$ to get $C^{2,\alpha}$ $($remember that we assumed that $M$ is a smooth differentiable manifold$)$ convergence of $u_{T_j}$. To achieve this goal one needs to write down carefully the dependence of all constants in the Euclidean proof of the regularity theorem of Allard to obtain a pure intrinsic Riemannian proof and then observe that indeed the constants involved depends just on the first and second derivatives of the metric and so they can be uniformly bounded over a sequence converging in $C^2$-topology in the sense of Definition $\ref{Def:Petersen}$, without any use of the Nash's isometric imbedding theorem. To be convinced of the $C^2$ dependence on the metric it is enough to remark that all these constants comes from a distortion of the metric due locally to the exponential map that is a bi-Lipschitz diffeomorphism. Hence the metric distortion depends on bounds on the sectional curvature of $M$, and so on the metric up to the second derivatives. The required details to make these arguments rigorous will be the object of a forthcoming paper.
\end{Rem}
      \newpage
      \markboth{Bibliography}{Bibliography}
      \bibliographystyle{alpha}
      \bibliography{these}

\begin{thebibliography}{DGCP72}

\bibitem[Ale62]{Alex}
A.~D. Alexandrov.
\newblock A characteristic property of spheres.
\newblock {\em Ann. Mat. Pura Appl.}, 58(4):303--315, 1962.

\bibitem[All72]{All}
William~K. Allard.
\newblock On the first variation of a varifold.
\newblock {\em Ann. of Math.}, 95:417--491, 1972.

\bibitem[Alm76]{Alm}
Frederick~J. Almgren.
\newblock {\em Existence and regularity almost everywhere of solutions to
  elliptic variational problems eith constraints}.
\newblock Number 165 in Memoirs of the American Mathematical Society. America
  Mathematical Society, 1976.

\bibitem[BBG85]{BessonGallotBerard}
P.~B\'erard, G.~Besson, and S.~Gallot.
\newblock Sur une in\'egalit\'e isop\'erim\'etrique qui g\'en\'eralise celle de
  {P}aul {L}\'evy-{G}romov.
\newblock {\em Invent. Math.}, 80(2):295--308, 1985.

\bibitem[BGG69]{BDeGG}
Enrico Bombieri, Ennio~De Giorgi, and Enrico Giusti.
\newblock Minimal cones and the bernstein problem.
\newblock {\em Invent. Math.}, 7:243--268, 1969.

\bibitem[BZ88]{BuragoZalgaller}
Yu.~D. Burago and V.~A. Zalgaller.
\newblock {\em Geometric inequalities}, volume 285 of {\em Grundlehren der
  Mathematischen Wissenschaften [Fundamental Principles of Mathematical
  Sciences]}.
\newblock Springer-Verlag, Berlin, 1988.
\newblock Translated from the Russian by A. B. Sosinski{\u\i}, Springer Series
  in Soviet Mathematics.

\bibitem[Cha06]{Chavel}
Isaac Chavel.
\newblock {\em Riemannian geometry}, volume~98 of {\em Cambridge Studies in
  Advanced Mathematics}.
\newblock Cambridge University Press, Cambridge, second edition, 2006.
\newblock A modern introduction.

\bibitem[DGCP72]{DeGiorgi}
E.~De~Giorgi, F.~Colombini, and L.~C. Piccinini.
\newblock {\em Frontiere orientate di misura minima e questioni collegate}.
\newblock Scuola Normale Superiore, Pisa, 1972.

\bibitem[EM02]{Eliashberg}
Y.~Eliashberg and N.~Mishachev.
\newblock {\em Introduction to the {$h$}-principle}, volume~48 of {\em Graduate
  Studies in Mathematics}.
\newblock American Mathematical Society, Providence, RI, 2002.

\bibitem[Fed69]{Fed}
Herbert Federer.
\newblock {\em Geometric Measure Theory}, volume 153 of {\em Grundelehren der
  mathematische wissenschaften}.
\newblock Springer Verlag, 1969.

\bibitem[GMT83]{GMTamanini}
E.~Gonzalez, U.~Massari, and I.~Tamanini.
\newblock On the regularity of boundaries of sets minimizing perimeter with a
  volume constraint.
\newblock {\em Indiana Univ. Math. J.}, 32(1):25--37, 1983.

\bibitem[GNP09]{GNP}
Renata Grimaldi, Stefano Nardulli, and Pierre Pansu.
\newblock Semianalyticity of isoperimetric profiles.
\newblock {\em Differ. Geom. Appl.}, 27(3):393--398, 2009.

\bibitem[God71]{God}
Claude Godbillon.
\newblock {\em \'{E}lments de Topologie Alg\'{e}brique}.
\newblock Hermann, 1971.

\bibitem[GR13]{GalliRitore}
Matteo Galli and Manuel Ritor\'e.
\newblock Existence of isoperimetric regions in contact sub-{R}iemannian
  manifolds.
\newblock {\em J. Math. Anal. Appl.}, 397(2):697--714, 2013.

\bibitem[Gra01]{Gra}
Alfred Gray.
\newblock {\em Tubes, Second Edition}, volume 221 of {\em Progress in
  Mathematics}.
\newblock Birk\"{a}user Verlag, 2001.

\bibitem[Gro86a]{Gr1}
Mikhael Gromov.
\newblock Isoperimetric inequalities in riemannian manifolds.
\newblock {\em Lectures notes in Mathematics}, 1200, 1986.

\bibitem[Gro86b]{Gr2}
Mikhael Gromov.
\newblock {\em Partial Differential Relations}.
\newblock Springer Verlag, 1986.

\bibitem[GT01]{GT}
David Gilbarg and Neil~S. Trudinger.
\newblock {\em Elliptic partial differential equations of second order}.
\newblock Classics in Mathematics. Springer-Verlag, Berlin, 2001.
\newblock Reprint of the 1998 edition.

\bibitem[Heb99]{Hebey}
Emmanuel Hebey.
\newblock {\em Nonlinear analysis on manifolds: {S}obolev spaces and
  inequalities}, volume~5 of {\em Courant Lecture Notes in Mathematics}.
\newblock New York University, Courant Institute of Mathematical Sciences, New
  York; American Mathematical Society, Providence, RI, 1999.

\bibitem[Hir94]{Hirsch}
Morris~W. Hirsch.
\newblock {\em Differential topology}, volume~33 of {\em Graduate Texts in
  Mathematics}.
\newblock Springer-Verlag, New York, 1994.
\newblock Corrected reprint of the 1976 original.

\bibitem[HK78]{HeintzeKarcher}
Ernst Heintze and Hermann Karcher.
\newblock A general comparison theorem with applications to volume estimates
  for submanifolds.
\newblock {\em Ann. Sci. \'Ecole Norm. Sup. (4)}, 11(4):451--470, 1978.

\bibitem[Mag12]{Maggi}
Francesco Maggi.
\newblock {\em Sets of finite perimeter and geometric variational problems},
  volume 135 of {\em Cambridge Studies in Advanced Mathematics}.
\newblock Cambridge University Press, Cambridge, 2012.
\newblock An introduction to geometric measure theory.

\bibitem[MJ00]{MJ}
Frank Morgan and David~L. Johnson.
\newblock Some sharp isoperimetric theorems for riemannian manifolds.
\newblock {\em Indiana Univ. Math. J.}, 49(2), 2000.

\bibitem[MN16]{FloresNardulli2016}
A.~{Mu{\~n}oz Flores} and S.~{Nardulli}.
\newblock {Local H\"{o}lder continuity of the isoperimetric profile in complete
  noncompact Riemannian manifolds with bounded geometry}.
\newblock {\em ArXiv e-prints}, June 2016.

\bibitem[Mor94]{Morgan94}
Frank Morgan.
\newblock Clusters minimizing area plus length of singular curves.
\newblock {\em Math. Ann.}, 299(4):697--714, 1994.

\bibitem[Mor03]{Morg1}
Frank Morgan.
\newblock Regularity of isoperimetric hypersurfaces in riemannian manifolds.
\newblock {\em Trans. Amer. Math. Soc.}, 355(12), 2003.

\bibitem[Nar09]{NarAnn}
Stefano Nardulli.
\newblock The isoperimetric profile of a smooth riemannian manifold for small
  volumes.
\newblock {\em Ann. Glob. Anal. Geom.}, 36(2):111--131, September 2009.

\bibitem[Nar14a]{NarAsian}
Stefano Nardulli.
\newblock Generalized existence of isoperimetric regions in non-compact
  {R}iemannian manifolds and applications to the isoperimetric profile.
\newblock {\em Asian J. Math.}, 18(1):1--28, 2014.

\bibitem[Nar14b]{NarCalcVar}
Stefano Nardulli.
\newblock The isoperimetric profile of a noncompact {R}iemannian manifold for
  small volumes.
\newblock {\em Calc. Var. Partial Differential Equations}, 49(1-2):173--195,
  2014.

\bibitem[Nas56]{Nash}
John Nash.
\newblock The imbedding problem for {R}iemannian manifolds.
\newblock {\em Ann. of Math. (2)}, 63:20--63, 1956.

\bibitem[NO16]{NardulliOsorio}
S.~{Nardulli} and L.~E. {Osorio Acevedo}.
\newblock {Sharp isoperimetric inequalities for small volumes in complete
  noncompact Riemannian manifolds of bounded geometry involving the scalar
  curvature}.
\newblock {\em ArXiv e-prints}, November 2016.

\bibitem[Ros05]{Ros}
Antonio Ros.
\newblock The isoperimetric problem.
\newblock In {\em Global theory of minimal surfaces}, volume~2 of {\em Clay
  Math. Proc.}, pages 175--209. Amer. Math. Soc., Providence, RI, 2005.

\bibitem[RR04]{RitoreRosales}
Manuel Ritor{\'e} and C{\'e}sar Rosales.
\newblock Existence and characterization of regions minimizing perimeter under
  a volume constraint inside {E}uclidean cones.
\newblock {\em Trans. Amer. Math. Soc.}, 356(11):4601--4622 (electronic), 2004.

\bibitem[Sak83]{SakaiContinuityOfInjectivityRadius}
Takashi Sakai.
\newblock On continuity of injectivity radius function.
\newblock {\em Math. J. Okayama Univ.}, 25(1):91--97, 1983.

\bibitem[Spr10]{Spring}
David Spring.
\newblock {\em Convex integration theory}.
\newblock Modern Birkh\"auser Classics. Birkh\"auser/Springer Basel AG, Basel,
  2010.
\newblock Solutions to the $h$-principle in geometry and topology, Reprint of
  the 1998 edition [MR1488424].

\end{thebibliography}
      \addcontentsline{toc}{section}{\numberline{}Bibliography}
      \emph{Stefano Nardulli\\ Departamento de Matem\'atica\\ Instituto de Matem\'atica\\ UFRJ-Universidade Federal do Rio de Janeiro, Brasil\\ email: nardulli@im.ufrj.br}  
\end{document}